\documentclass{sig-alternate-no-acm}
\usepackage{fullpage}
\usepackage{amsmath,mathrsfs,verbatim,color,lscape,xspace}
\allowdisplaybreaks
\usepackage{listings,mathtools}
\usepackage{graphicx}
\usepackage{secdot}
\usepackage{dsfont}
\usepackage{amsfonts}
\usepackage{amssymb}
\usepackage{enumerate,enumitem}
\usepackage{xspace,ifpdf}
\usepackage[tight]{subfigure}
\usepackage[shortcuts]{extdash}
\usepackage{hyperref}      %Use only when compiling pdf's, otherwise it creates errors
\usepackage{xcolor}\hypersetup{pdftex,breaklinks,citebordercolor=green,linkbordercolor=magenta}
\newcommand{\dltpts}[2][s]{\ensuremath{\widetilde{F}_{p}^{(#2)}(#1)}\xspace}        %for the derivative of the L-T of the phase-type service times

\newcommand{\lt}[1]{\ensuremath{\widetilde{#1}}}          %for L.T. of f, I have to type \LT{f}, argument not optional

\newcommand{\fun}[2][s]{\ensuremath{{#2}(#1)}}     %function of f with default parameter s,

\newcommand{\e}[1][]{\ensuremath{\mathbb{E}{#1}}\xspace}    %for the expectation
\newcommand{\pr}[1][]{\ensuremath{\mathbb{P}{#1}}\xspace}   %for the probability
\newcommand{\dm}[2][n]{\ensuremath{{#2}^{(#1)}}\xspace}   %for the derivative of matrices
\newcommand{\df}[2][]{\ensuremath{{#2}^{(#1)}}\xspace}   %for the derivative of functions
\newcommand{\ind}[3]{\ensuremath{{#1}^{#2}_{#3}}}           %for the indexes up and down
\newcommand{\vect}[1]{\ensuremath{\mathbf{#1}}\xspace}  %for the vector notation &&&&&&&&&&&&&&&&&&&&&&&&&&& CHECK IF I NEED IT &&&&&&&&&&&&&&&&&&&&
\newcommand{\diag}{\ensuremath{\text{diag}}\xspace}  %for the word diag for diagonal matrices
\newcommand{\n}[1][]{\ensuremath{\mathcal{N}}\xspace}    %for the calligraphic N
\newcommand{\rank}{\ensuremath{\textup{rank\/}}\xspace}  %for the word rank, for matrices
\newcommand{\mean}[1][]{\ensuremath{\mu_{#1}}\xspace}    %for the mean of the general service time distribution
\newcommand{\indfun}[1][]{\ensuremath{\mathds{1}_{#1}}\xspace}       %for the indicator function

\newcommand{\erlang}[2][]{\ensuremath{E_{#1}(#2)}\xspace}       %For the Erlang_k(l) distribution, when k=1 we mean the exponential distribution

\newcommand{\rootp}[1]{\ensuremath{s_{#1}}\xspace}    %for the roots with positive real part
\newcommand{\rootpmix}[1]{\ensuremath{s_{\epsilon,#1}}\xspace}    %for the roots with positive real part mixture model
\newcommand{\rootnden}[1]{\ensuremath{x_{#1}}\xspace}    %for the roots with negative real part on the denominator
\newcommand{\rootnnum}[1]{\ensuremath{y_{#1}}\xspace}    %for the roots with negative real part on the numerator
\newcommand{\workload}[1][]{\ensuremath{V_{#1}}\xspace}    %for the workload random variable
\newcommand{\workloadmix}[1][]{\ensuremath{V_{\epsilon #1}}\xspace}    %for the workload random variable
\newcommand{\epts}[1][]{\ensuremath{B^e_{#1}}\xspace}                 %for the stationary-excess phase-type service times
\newcommand{\ehts}[1][]{\ensuremath{C^e_{#1}}\xspace}                 %for the stationary-excess heavy-tailed service times

\newcommand{\las}[1]{\ensuremath{\boldsymbol\lambda^{#1}}\xspace}    %for the set of lambdas
\newcommand{\mxtrans}{\ensuremath{\mathbf{P}}\xspace}           %for the transition probability matrix P
\newcommand{\mxrates}{\ensuremath{\mathbf{\Lambda}}\xspace}           %for the rate \Lambda
\newcommand{\mxprobsdummy}{\ensuremath{\mathbf{Q}^{(1)}}\xspace}           %for the matrix Q of probabilities of dummy customers
\newcommand{\mxprobsreal}{\ensuremath{\mathbf{Q}^{(2)}}\xspace}           %for the matrix Q of probabilities of real customers
\newcommand{\ltmgs}[1][s]{\ensuremath{\widetilde{\mathbf{G}}(#1)}\xspace}     %for the matrix of the Laplace transforms of the service time distributions
\newcommand{\ltmgsmix}[1][s]{\ensuremath{\widetilde{\mathbf{G}}_\epsilon(#1)}\xspace}     %for the matrix of the Laplace transforms of the service time distributions of the mixture model
\newcommand{\uv}{\ensuremath{\mathbf{e}}\xspace}           %for the unit vector
\newcommand{\um}{\ensuremath{\mathbf{U}}\xspace}           %for the unit matrix, matrix with all units
\newcommand{\im}[1][]{\ensuremath{\mathbf{I}_{#1}}\xspace}           %for the Identity matrix
\newcommand{\vup}{\ensuremath{\mathbf{u}}\xspace}           %for the vector with unknown parameters
\newcommand{\vupmix}{\ensuremath{\mathbf{u}_\epsilon}\xspace}           %for the vector with unknown parameters in the mixture model
\newcommand{\va}[1][]{\ensuremath{\mathbf{a}_{#1}}\xspace}           %for the vector a
\newcommand{\vc}[1][]{\ensuremath{\mathbf{c}_{#1}}\xspace}           %for the vector c
\newcommand{\vd}[1][]{\ensuremath{\mathbf{d}_{#1}}\xspace}           %for the vector d
\newcommand{\mxadj}[1][]{\ensuremath{\mathbf{\mathcal{E}}_{#1}}\xspace}     %for the adjoint matrix of E(s)
\newcommand{\mxadjmix}[1][]{\ensuremath{\mathbf{\mathcal{E}}_{\epsilon #1}}\xspace}     %for the adjoint matrix of E(s) in the mixture model
\newcommand{\mxa}[1][]{\ensuremath{\mathbf{A}_{#1}}\xspace}     %for the matrix A{ij}
\newcommand{\mxb}[1][]{\ensuremath{\mathbf{B}_{#1}}\xspace}     %for the matrix B{ij}
\newcommand{\mxex}{\ensuremath{\mathbf{C}}\xspace}     %for the matrix \Gamma_{ij}, wich we use as example for introduction of new terminology or notation
\newcommand{\mxexs}{\ensuremath{\mathbf{D}}\xspace}     %for the matrix B_{ij}, wich we use as second example for introduction of new terminology or notation
\newcommand{\mxe}[1][]{\ensuremath{\mathbf{E}_{#1}}\xspace}     %for the matrix E_{ij}
\newcommand{\mxemix}[1][]{\ensuremath{\mathbf{E}_{\epsilon #1}}\xspace}     %for the matrix E_{ij} of the mixture model
\newcommand{\mxh}[1][]{\ensuremath{\mathbf{H}_{#1}}\xspace}     %for the matrix H_{ij}
\newcommand{\mxk}[1][]{\ensuremath{\mathbf{K}_{#1}}\xspace}     %for the matrix K_{ij}
\newcommand{\mxhmix}[1][]{\ensuremath{\mathbf{H}_{\epsilon #1}}\xspace}     %for the matrix H_{ij} of the mixture model
\newcommand{\mxmeans}[1][]{\ensuremath{\mathbf{M}_{#1}}\xspace}     %for the matrix with the mean service times
\newcommand{\mxmeansmix}[1][]{\ensuremath{\mathbf{M}_{\epsilon #1}}\xspace}     %for the matrix with the mean service times in the mixture model
\newcommand{\ltm}[2][s]{\ensuremath{#2(#1)}\xspace}             %for the L-T of a matrix
\newcommand{\ltpts}[1][s]{\ensuremath{\widetilde{F}_{p}(#1)}\xspace}        %for the L-T of the phase-type service times
\newcommand{\lthts}[1][s]{\ensuremath{\widetilde{F}_{h}(#1)}\xspace}        %for the L-T of the heavy-tailed service times
\newcommand{\ltptes}[1][s]{\ensuremath{\widetilde{F}^e_{p}(#1)}\xspace}        %for the Laplace transform of the phase-type stationary-excess service time distribution
\newcommand{\lthtes}[1][s]{\ensuremath{\widetilde{F}^e_{h}(#1)}\xspace}        %for the Laplace transform of the heavy-tailed stationary-excess service time distribution
\newcommand{\ltgs}[1][]{\ensuremath{\widetilde{G}^{#1}(s)}\xspace}        %for the L-T of the general service times
\newcommand{\ltgsc}[1]{\ensuremath{\widetilde{G}_{#1}(s)}\xspace}        %for the L-T of the general service times per class
\newcommand{\ltgsmix}[1][]{\ensuremath{\widetilde{G}^{#1}_\epsilon(s)}\xspace}        %for the L-T of the general service times of the mixture model
\newcommand{\ltgsdis}[1][]{\ensuremath{\widetilde{G}^{\bullet}_\epsilon(s)}\xspace}        %for the L-T of the general service times of the discard model
\newcommand{\vltw}[1][s]{\ensuremath{\widetilde{\mathbf{\Phi}}(#1)}\xspace}        %for the vector of the L-T with the workload distributions
\newcommand{\vltwmix}[1][s]{\ensuremath{\widetilde{\mathbf{\Phi}}_\epsilon(#1)}\xspace}        %for the vector of the L-T with the waiting time distributions of the mixture model
\newcommand{\vltwmixder}[1][s]{\ensuremath{\dm[1]{\widetilde{\mathbf{\Phi}}}_\epsilon(#1)}\xspace}        %for the vector of the first derivative L-T with the workload distributions of the mixture model
\newcommand{\ltwc}[2][s]{\ensuremath{\widetilde{\phi}_{#2}(#1)}\xspace}        %for the L-T of the waiting time distribution of each class
\newcommand{\ltwcmix}[2][s]{\ensuremath{\widetilde{\phi}_{\epsilon #2}(#1)}\xspace}        %for the L-T of the waiting time distribution of each class of the mixture model
\newcommand{\ptsd}[1][t]{\ensuremath{F_{p}(#1)}\xspace}        %for the phase-type service time distribution
\newcommand{\htsd}[1][t]{\ensuremath{F_{h}(#1)}\xspace}        %for the heavy-tailed service time distribution
\newcommand{\ptesd}[1][t]{\ensuremath{F^e_{p}(#1)}\xspace}        %for the phase-type stationary-excess service time distribution
\newcommand{\htesd}[1][t]{\ensuremath{F^e_{h}(#1)}\xspace}        %for the heavy-tailed stationary-excess service time distribution
\newcommand{\gsd}[1][t]{\ensuremath{G(#1)}\xspace}        %for the general service time distribution
\newcommand{\gsdmix}[1][t]{\ensuremath{G_\epsilon(#1)}\xspace}        %for the general service time distribution of the mixture model
\newcommand{\gsdc}[2][t]{\ensuremath{G_{#2}(#1)}\xspace}        %for the general service time distribution for different classes

\newcommand{\fef}[1][s]{\ensuremath{{n}(#1)}\xspace}        %function example n(s) first
\newcommand{\fes}[1][s]{\ensuremath{{d}(#1)}\xspace}        %function example d(s) second
\newcommand{\correpc}[2][t]{\ensuremath{\widehat{\varphi}_{r,\epsilon,#2}(#1)}\xspace}         %for the corrected replace approximation
\newcommand{\corsimrepc}[2][t]{\ensuremath{\widehat{\varphi}_{s.r,\epsilon,#2}(#1)}\xspace}         %for the simplified corrected replace approximation
\newcommand{\correp}[1][t]{\ensuremath{\widehat{\varphi}_{r,\epsilon}(#1)}\xspace}         %for the corrected replace approximation
\newcommand{\corsimrep}[1][t]{\ensuremath{\tilde{\varphi}_{s.r,\epsilon}(#1)}\xspace}         %for the simplified corrected replace approximation
\newcommand{\delay}[1][]{\ensuremath{W^{#1}}\xspace}    %for the workload random variable
\newcommand{\delaymix}[1][]{\ensuremath{W^{#1}_{\epsilon}}\xspace}    %for the workload random variable mixture model

\newcommand{\ltd}[1][s]{\ensuremath{\widetilde{w}(#1)}\xspace}        %for the L-T of the delay distribution
\newcommand{\ltdmix}[1][s]{\ensuremath{\widetilde{w}_\epsilon(#1)}\xspace}        %for the L-T of the delay distribution in the mixture model

\newcommand{\vweights}{\ensuremath{\boldsymbol \omega}\xspace}
\newcommand{\weight}[1]{\ensuremath{\omega_{#1}}\xspace}    %for the weights in front of the LT worklaod for the LI waiting time

\newtheorem{dummy}{Dummy}[section]
\newtheorem{proposition}[dummy]{Proposition}
\newtheorem{lemma}[dummy]{Lemma}

\newdef{definition}{Definition}
\newtheorem{theorem}[dummy]{Theorem}
\newtheorem{corollary}[dummy]{Corollary}
\newdef{remark}{Remark}
\newdef{approximation}{Approximation}
\newenvironment{toybeginning}{\paragraph{\textbf{\textsl{Running example}}}}{\hfill$\blacksquare$}
\newenvironment{toy}{\paragraph{\textbf{\textsl{Running example (continued)}}}}{\hfill$\blacksquare$}

\begin{document}

%\conferenceinfo{MAM8}{2014 Calicut, Kerala India}
%\CopyrightYear{2007} % Allows default copyright year (20XX) to be over-ridden - IF NEED BE.
%\crdata{0-12345-67-8/90/01}  % Allows default copyright data (0-89791-88-6/97/05) to be over-ridden - IF NEED BE.
% --- End of Author Metadata ---

\title{Corrected phase\-/type approximations of heavy\-/tailed queueing models in a Markovian environment}
\author{
\begin{tabular}{c@{\hspace*{2em}}c@{\hspace*{2em}}c@{\hspace*{2em}}c}
Eleni Vatamidou$^\star$                & Ivo Adan$^{\star,\dag}$       & Maria Vlasiou$^{\star,\ddagger}$             & Bert Zwart$^{\star,\ddagger}$ \\
\texttt{e.vatamidou@tue.nl}    & \texttt{i.j.b.f.adan@tue.nl}  & \texttt{m.vlasiou@tue.nl}  & \texttt{Bert.Zwart@cwi.nl}\\
\end{tabular}\vspace{0.25cm}\\
%%%
%%% Now we list the addresses horizontally:
%%%
{\footnotesize $^\star$\textsc{Eurandom} and Department of Mathematics \& Computer Science, Eindhoven University of Technology,}\\
{\footnotesize  P.O. Box 513, 5600 MB Eindhoven, The Netherlands}\\
{\footnotesize $^\dag$Centrum Wiskunde \& Informatica (CWI), P.O. Box 94079, 1090 GB Amsterdam, The Netherlands}\\
{\footnotesize $^\ddagger$Department of Mechanical Engineering, Eindhoven University of Technology,}\\
{\footnotesize  P.O. Box 513, 5600 MB Eindhoven, The Netherlands}
}

\maketitle
\begin{abstract}
Significant correlations between arrivals of load\-/generating events make the numerical evaluation of the workload of a system a challenging problem. In this paper, we construct highly accurate approximations of the workload distribution of the MAP/G/1 queue that capture the tail behavior of the exact workload distribution and provide a bounded relative error. Motivated by statistical analysis, we consider the service times as a mixture of a phase\-/type and a heavy\-/tailed distribution. With the aid of perturbation analysis, we derive our approximations as a sum of the workload distribution of the MAP/PH/1 queue and a heavy\-/tailed component that depends on the perturbation parameter. We refer to our approximations as {\it corrected phase\-/type approximations}, and we exhibit their performance with a numerical study.
\end{abstract}


\section{Introduction}\label{s3.Introduction}
%\blue{what is the problem and what we do}
The evaluation of performance measures in stochastic models is a key problem that has been widely studied in the literature \cite{abate94a,asmussen-RP,knessl87a,wu99}. In this paper, we focus on the evaluation of the workload distribution of a single server queue where customers arrive according to a Markovian Arrival Process (MAP) \cite{asmussen93,lucantoni90} and their service times follow some general distribution. Under the presence of heavy\-/tailed service times, such evaluations become more challenging and sometimes even problematic \cite{ahn12,asmussen05}. In such cases, it is necessary to construct approximations. In this study, we propose to modify existing approximations by adding a small refinement term, which can serve two purposes. On the one hand, the refinement term helps in constructing approximations not only with a small absolute error, but also with a small relative error. On the other hand, it gives information on the accuracy of the approximation without the modification: the smaller the refinement term, the better the pre\-/modified approximation.

%\blue{definition of a MAP}
An important generalization of the Poisson point process is the MAP. In a MAP, the arrivals are not homogenous in time, but they are determined by a Markov process $\{J_t\}_{t\geq 0}$ with a finite state space. The class of MAPs is a very rich class of point processes, containing many well\-/known arrival processes as special cases. A special case of a MAP is the Markov\-/modulated Poisson process (MMPP), which is a popular model for bursty arrivals \cite{fisher93}. The class of MAPs contains also the class of phase\-/type renewal processes, i.e.\ renewal processes with phase\-/type interarrivals \cite{neuts78}. %The intensity matrix governing $\{J_t\}$ can be decomposed as the sum of two matrices, where the first matrix is related to phase transitions not associated with arrivals, while transitions in the second matrix are related to arrivals. %More precisely, a transition from state $i$ to state $j$ is accompanied by an arrival of a customer with rate that is equal to the $(i,j)$ element of the second matrix, while with rate that is equal to the $(i,j)$ element of the first matrix it is not accompanied by any arrival.

%\blue{importance of MAPs}
%It has been proven that the stationary MAPs are dense in the family of all stationary point processes \cite{asmussen93}. Due to their importance, many algorithms have been developed to fit any general arrival process to a MAP \cite{breuer02,buchholz10,okamura11,ryden96}. As a consequence of the fact that MAPs can be used to describe any arrival process, the workload of a MAP/G/1 queue can be used as an approximation of the workload of a G/G/1 queue.

%\blue{closed form expressions for the workload}
It has been shown that the Laplace transform of the workload of a MAP/G/1 queue has a matrix expression analogous to the Pollazceck\-/Khinchine equation of an M/G/1 queue \cite{neuts-SSM,ramaswami80}. However, these closed\-/form expressions are only practical in case of phase\-/type service times \cite{asmussen00a,asmussen-APQ}, where the workload distribution has a phase\-/type representation \cite{ramaswami90} in a form which is explicit up to the solution of a matrix functional equation.%fixed point problem.

%\blue{phase\-/type approximations, characteristics, drawbacks}
Since the class of phase\-/type distributions is dense in the class of all distributions on $(0,\infty)$ \cite{asmussen00a}, a common approach to approximate the workload is by approximating the service time distribution with a phase\-/type one; see e.g.\ \cite{feldmann98,starobinski00}. We refer to these methods as {\it phase\-/type approximations}. There are many algorithms for phase\-/type approximations, which provide highly accurate approximations for the workload distribution when the service times are light\-/tailed. However, in many cases, a heavy\-/tailed distribution is most appropriate to model the service times \cite{embrechts-MEE,rolski-SPIF}. In these cases, the exponential decay of phase\-/type approximations gives a big relative error at the tail and the evaluation of the workload becomes more complicated. Since heavy\-/tailed distributions have cumbersome expressions for their Laplace transform, this prevents the usage of techniques that require transform expressions, such as \cite{gail96}.

%\blue{The gap we are filling in., inspiration for our approach}
In this paper, we develop approximations of the workload distribution for heavy\-/tailed service times that maintain the computational tractability of phase\-/type approximations, capture the correct tail behavior and provide small absolute and relative errors. In order to achieve these desirable characteristics, our key idea is to use a mixture model for the service times. The idea of our approach stems from fitting procedures of the service time distribution to data. Heavy\-/tailed statistical analysis suggests that only a small fraction of the upper\-/order statistics of a sample is relevant for estimating tail probabilities \cite{resnick-HTP}. The remaining data set may be used to fit the bulk of the distribution, where a natural choice is to fit a phase\-/type distribution to the remaining data set \cite{asmussen96b}. As a result, a mixture model for the service times is a natural assumption.

%\blue{Explanations for the technique of our approximation.}
We now briefly explain how to derive our approximations when the service time distribution is a mixture of a phase\-/type distribution and a heavy\-/tailed one. We show that if the service time distribution is such a mixture, then the workload can also be written as a mixture, in the sense that it involves the workload of a model with purely phase\-/type service times and some additional terms related to the heavy\-/tailed distribution of our mixture model. Consequently, we first need to compute the workload in a MAP/PH/1 queue and afterwards use this as a base to calculate the rest of the terms involving the heavy\-/tailed distribution.

As a first step to derive our approximations, we write the service time distribution as perturbation of the phase\-/type distribution by a function that contains the heavy\-/tailed component. By ignoring the perturbation term and by taking the service time distribution equal to the phase\-/type distribution, we find the workload of a resulting simpler MAP/PH/1 queue, which is a phase\-/type approximation of the workload. By applying perturbation analysis to all parameters that depend on the service time distribution, we can write the workload as a series expansion, where the constant term is the workload of the MAP/PH/1 queue used as base and all other terms contain the heavy\-/tailed component.

%cause ruin

Large deviations theory suggests that a single catastrophic event, i.e.\ a stationary heavy\-/tailed service time, is sufficient to give a non\-/zero tail probability for the workload \cite{embrechts-MEE}. As we will see in Section~\ref{s3.workload distribution of the perturbed model}, the second term of the series expansion of the workload can be expressed in terms of such a catastrophic event. Thus, we define our approximations as the sum of the first two terms of the series expansion of the workload, and we show that the addition of the second term leads to improved approximations when compared to their phase\-/type counterparts. In other words, the second term makes the phase\-/type approximation more robust so that the relative error at the tail does not explode. Therefore, we call this term correction term, and inspired by the terminology {\it corrected heavy traffic approximations} \cite{asmussen-APQ} we refer to our approximations as {\it corrected phase\-/type approximations}. In a previous study \cite{vatamidou13}, we applied this approach to Poisson arrivals.

%\blue{connection with other models, possible extensions}
The connection between the stationary workload distribution of a MAP/G/1 queue and ruin probabilities for a risk process in a Markovian environment, where the claim sizes in the risk model correspond to the service times and the arrival process of claims is the time\-/reversed MAP of the queueing model, is well known \cite{asmussen-APQ,asmussen-RP}. Thus, the corrected phase\-/type approximations can also be used to estimate the ruin probabilities of the above mentioned risk model. Finally, our technique can be applied to more general queueing models, i.e.\ queuing models with dependencies between interarrival and service times \cite{boxma01,smits04}, and also to models that allow for customers to arrive in batches (the arrival process is called Batch Markovial Arrival Process) \cite{lucantoni91,lucantoni93-BMAP,lucantoni94}.
%, namely the BMAP/G/1 queues

%\blue{closely related work of Ivo, differences, applicability. Explain what is done there, what we do here. Simplifications, generalizations of his work.}
A closely related work is Adan and Kulkarni \cite{adan03}. They consider a single server queue, where the interarrival times and the service times depend on a common discrete Markov Chain. In addition, they assume that a customer arrives in each phase transition, and they find a closed form expression for the waiting time distribution under general service time distributions. However, when there exist also phase transitions not related to arrivals of customers, their results remain valid for the evaluation of the workload. This can be seen by using the standard technique of including {\it dummy} customers in the model; namely customers with zero service time.

%\blue{Setup of the rest paper.}
The rest of the paper is organized as follows. In Section~\ref{s3.presentation of the model}, we introduce the model under consideration without assuming any special form for the service time distribution, and in Section~\ref{s3.preliminaries} we find the general expressions for the Laplace transforms of the workload prior to a transition from each state. In Section~\ref{s3.Construction of the corrected phase-type approximations}, we consider service time distributions that are a mixture of a phase\-/type distribution and a heavy\-/tailed one, and we explain the idea to construct our approximations. Later in Section~\ref{s3.replace base model}, we specialize the results of Section~\ref{s3.preliminaries} for phase\-/type service times. We use as base model the phase\-/type model of Section~\ref{s3.replace base model}, and we apply perturbation analysis to find in Section~\ref{s3.parameters of the perturbed system} the perturbed parameters and in Section~\ref{s3.workload distribution of the perturbed model} the desired Laplace transforms of the workload in the mixture model. Using the latter results, we construct in Section~\ref{s3.properties corrected replace approximation} the approximations and we discuss their properties. Finally, in Section~\ref{s3.numerics}, we use a specific mixture service time distribution for which the exact workload distribution can be calculated and we exhibit the accuracy of our approximations through numerical experiments. Finally, in the Appendix, we give the proofs of all theorems, the necessary theory on perturbation analysis, and other related results. Due to the complexity of the formulas, we use a simple running example in order to explain the idea behind the calculations.

\section{Presentation of the model}\label{s3.presentation of the model}
We consider a single server queue with FIFO discipline, where customers arrive according to a Markovian Arrival Process (MAP). The arrivals are regulated by a Markov process $\{J_t\}_{t\geq 0}$ with a finite state space \n, say with $N$ states. We assume that the service time distribution of a customer is independent of the state of $\{J_t\}$ upon his arrival. For this model, we are interested in finding accurate approximations for the workload distribution.

The intensity matrix \vect{D} governing $\{J_t\}$ is denoted by the decomposition $ \vect{D}=\vect{D}^{(1)}+\vect{D}^{(2)}$, where the matrix $\vect{D}^{(1)}$ is related to arrivals of {\it dummy} customers, while transitions in $\vect{D}^{(2)}$ are related to arrivals of {\it real} customers. Note that the diagonal elements of the matrix $\vect{D}^{(2)}$ may not be identically equal to zero. This means that if $d_{ii}^{(2)}>0$, then a real customer arrives with rate $d_{ii}^{(2)}$ and we have a transition from state $i$ to itself. However, phase transitions not associated with arrivals (dummy customers) from any state to itself are not allowed. Since the matrix \vect{D} is an intensity matrix, its rows sum up to zero. Therefore, the diagonal elements of the matrix $\vect{D}^{(1)}$ are negative and they are defined as $d_{ii}^{(1)} = -\sum_{k\neq i}d_{ik}^{(1)} - \sum_{k=1}^Nd_{ik}^{(2)}$.

%The class of MAPs is a very rich class of point processes, containing many well\-/known arrival processes as special cases. A special case of a MAP is the Markov\-/modulated Poisson process (MMPP), which is a popular model for bursty arrivals \cite{fisher93}. This is defined in terms of a background Markov process $\{J_t\}_{t\geq 0}$ with a finite number of states, such that the arrival intensity is $\beta_i$ on time intervals where $J_t=i$. In the setup of MAPs, the matrix $\vect{D}^{(2)}$ is equal to the diagonal matrix $\diag(\beta_1,\dots,\beta_N)$, which means that arrivals of real customers occur only when there is a transition from a state to itself. The class of MAPs contains also the class of phase\-/type renewal processes, i.e.\ renewal processes with phase\-/type interarrivals \cite{neuts78}. Consider for example a renewal process with interarrival distribution of phase type with representation $(\boldsymbol \alpha,\mathbf{T})$ (the corresponding exit rate vector is $\mathbf{t}=-\mathbf{T} \uv$, where \uv is the column vector with appropriate dimensions and all elements equal to 1). The intensity matrix of the background Markov process is then $\vect{D}= \mathbf{T} + \mathbf{t} \boldsymbol \alpha$, where the matrices $\mathbf{T}$ and $\mathbf{t} \boldsymbol \alpha$ correspond to state changes without and with arrivals, respectively.

In this paper, we are interested in modeling heavy\-/tailed service times. As stated earlier, motivated by statistical analysis, we assume that the service time distribution of a real customer is a mixture of a phase\-/type distribution, \ptsd, and a heavy\-/tailed one, \htsd. Namely, the service time distribution of a real customer has the form
\begin{equation}\label{e3.mixture service time distribution}
  \gsdmix = (1-\epsilon)\ptsd + \epsilon \htsd, %\qquad \epsilon \in [0,1).
\end{equation}
where $\epsilon$ is typically small.

Our goal is to find the workload distribution for this mixture model. Towards this direction, we present in the next section existing results \cite{adan03} for the evaluation of the workload distribution under the assumption of generally distributed service times. Ultimately, we wish to specialize these results to service times of the aforementioned form \eqref{e3.mixture service time distribution}. %As we will explain in the sequel, although the expressions provided in \cite{adan03} are compact, they are not practical when the service time distribution of a real customers is some heavy\-/tailed distribution. In this paper, we explain how to overcome this problem by suggesting a new approach and by providing an algorithm that constructs highly accurate approximations for the measure under consideration.

\subsection{Preliminaries}\label{s3.preliminaries}
Since the results of this section are valid for any service time distribution, we suppress the index $\epsilon$ and we use the notation \gsd for the service time distribution of a real customer. We consider now the embedded Markov chain $\{Z_n\}_{n\geq 0}$ on the arrival epochs of customers (real and dummy) and we denote by \mxtrans the transition probability matrix of the regulating Markov chain $\{Z_n\}$, which we assume to be irreducible. If $\lambda_i$ is the exponential exit rate from state $i$, i.e.
\begin{equation}\label{e3.definition rates}
  \lambda_i = \sum_{k\neq i}d_{ik}^{(1)} + \sum_{k=1}^Nd_{ik}^{(2)},
\end{equation}
the transition probabilities can be calculated by
\begin{equation}\label{e3.definition transition probabilities}
  p_{ij} = \frac{d_{ij}^{(1)}(1-\delta_{ij})+d_{ij}^{(2)}}{\lambda_i},
\end{equation}
where $\delta_{ij}$ is the Kronecker delta ($\delta_{ij}=0$ when $i\neq j$ and $\delta_{ij}=1$ when $i=j$). In addition, an arriving customer at a transition from state $i$ to state $j$ is tagged $i$. If $p_{ij}>0$, then we define the probability %Recall that the elements $d_{ii}^{(1)}$ are negative and make the sum of their corresponding row of the intensity matrix equal to zero, which explains their exclusion both from the exit rate $\lambda_i$ and the $p_{ii}$; namely the transition probabilities to the same state.
\begin{equation}\label{e3.definition conditional probabilities of dummy customer}
  q^{(1)}_{ij} = \frac{d_{ij}^{(1)}(1-\delta_{ij})}{d_{ij}^{(1)}(1-\delta_{ij})+d_{ij}^{(2)}},
\end{equation}
which is the probability of an arriving customer to be dummy conditioned on the event that there is a phase transition from state $i$ to $j$. Similarly, conditioned on the event that there is a phase transition from $i$ to $j$, the arriving customer is real with probability
\begin{equation}\label{e3.definition conditional probabilities of real customer}
  q^{(2)}_{ij} = \frac{d_{ij}^{(2)}}{d_{ij}^{(1)}(1-\delta_{ij})+d_{ij}^{(2)}}.
\end{equation}
If $p_{ij}=0$, then we define $q^{(1)}_{ij}=q^{(2)}_{ij}=0$. Consequently, the conditional service time distribution of an arriving customer at a transition from $i$ to $j$ is $\gsdc{ij} = q^{(1)}_{ij} + q^{(2)}_{ij} \gsd$, and its Laplace\-/Stieltjes transform (LST) is $\ltgsc{ij} = q^{(1)}_{ij}+ q^{(2)}_{ij} \ltgs $, $i,j=1,\dots,N$, where $\ltgs$ is the LST of the service time distribution \gsd of a real customer. In matrix form, the above quantities can be written as
\begin{align}
  \mxrates      &= \diag(\lambda_1,\dots,\lambda_N), \label{e3.definition rate matrix}\\
  \mxprobsdummy &= [q^{(1)}_{ij}], \label{e3.definition matrix of probabilities of dummy customers}\\
  \mxprobsreal  &= [q^{(2)}_{ij}], \label{e3.definition matrix of probabilities of real customers}\\
  \ltmgs        &= \mxprobsdummy +\ltgs \mxprobsreal. \label{e3.definition l-t general service times matrix}
\end{align}
Let now $\circ$ denote the Hadamard product between two matrices of same dimensions; i.e.\ if $\vect{B}=(b_{ij})$ and $\vect{C}=(c_{ij})$ are $m \times n$ matrices, then the $(i,j)$ element of the $m \times n$ matrix $\vect{B} \circ \vect{C}$ is equal to $b_{ij}c_{ij}$. We also define the matrix
\begin{equation}\label{e3.matrix h}
    \ltm{\mxh} = \ltmgs \circ \mxtrans \mxrates,
\end{equation}
which we will need later. Finally, let $\boldsymbol\pi = [\pi_1,\dots,\pi_N]$ be the stationary distribution of $\{Z_n\}_{n\geq 0}$, and \mean be the mean of the service time distribution \gsd. Then the system is stable if the mean service time of a customer is less than the mean inter\-/arrival times between two consecutive customers in steady state. Namely,
\begin{equation}\label{e3.stability condition}
  \boldsymbol\pi \left(\mxrates^{-1} - \mxmeans\right) \uv >0,
\end{equation}
where $\mxmeans=\mean \mxprobsreal\circ \mxtrans $ and \uv is the column vector with appropriate dimensions and all elements equal to 1. Note that the $(i,j)$ element of the matrix $\mxprobsreal\circ \mxtrans $ is the unconditional probability that a real customer arrives at a transition from $i$ to $j$.

From this point on, we use a simple running example so that we display the involved parameters and the derived formulas. The running example evolves progressively, which means that its parameters are introduced only once and the reader should consult a previous block of the example to recall the notation.

\begin{toybeginning}
   For our running example, we consider a MAP with Erlang-2 distributed interarrival times, where the exponential phases have both rate $\lambda$ ($N=2$). Therefore, the matrices $\vect{D}^{(1)}$  and $\vect{D}^{(2)}$ are given as follows:
  \begin{equation*}
   \vect{D}^{(1)} =  \left(\begin{array}{rr}
                        -\lambda    &\lambda \\
                        0           &-\lambda
                   \end{array} \right) \qquad \text {and} \qquad
   \vect{D}^{(2)} =  \left(\begin{array}{cc}
                        0    &0 \\
                        \lambda           &0
                   \end{array} \right).
 \end{equation*}
 In this case, we have that $\lambda_1=\lambda_2=\lambda$, $p_{ij}=1-\delta_{ij}$, $q^{(1)}_{12} = q^{(2)}_{21}=1$, and all other elements of the matrices \mxprobsdummy and \mxprobsreal are equal to zero. Observe that we only have transitions from state~1 to state~2 and from state~2 to state~1. Therefore, in state~1 we always have arrivals of dummy customers while in state~2 we only have arrivals of real customers. Thus, only the diagonal elements of the matrix \ltmgs are not equal to zero, so that $\ltgsc{11}=1$ and $\ltgsc{22}=\ltgs$. Finally, the stability condition takes its known form $\lambda \mean/2 <1$.
\end{toybeginning}

Let now \workload denote the steady\-/state workload of the system just prior to an arrival of a customer. If the arriving customer is real, then the workload just prior to its arrival equals the waiting time of the customer in the queue, which we denote by \delay. In terms of Laplace transforms, the steady\-/state workload of the system just prior to an arrival of a customer in state $i$ is found as
\begin{equation*}
  \ltwc{i} = \e{(e^{-s\workload};Z=i)}, \qquad \Re(s) \geq 0, \quad i =1,\dots,N,
\end{equation*}
where $Z$ is the steady\-/state limit of $Z_n$.\ Gathering all the above Laplace transforms \ltwc{i}, $i=1,\dots,N$, we construct the transform vector
\begin{equation}\label{e3.transform vector}
  \vltw = [\ltwc{1},\dots,\ltwc{N}].
\end{equation}
%Then, the Laplace transform of the workload \workload of the system just prior to an arrival\footnote{With the inclusion of dummy customers in the model, every phase transition is accompanied by an arrival of a customer.} is
%\begin{align}
%  \ltw   =& \e{\big[e^{-s\workload}\big]} = \sum_{i=1}^N \e{\big[e^{-s\workload};Z=i\big]} \notag \\
%         =& \sum_{i=1}^N \ltwc{i}
%         = \vltw \uv, \qquad \Re(s) \geq 0.\label{e3.total workload}
%\end{align}

We first provide some general theorems for the transform vector \vltw, which we later on refine in order to provide more detailed information regarding the form of the elements \ltwc{i}, $i=1,\dots,N$. In the following, \im stands for the identity matrix, with appropriate dimensions.

\begin{theorem}\label{t3.equations for the transform vector}
  Provided that the stability condition \eqref{e3.stability condition} is satisfied, the transform vector \vltw satisfies
  \begin{align}
    \vltw \big(\ltm{\mxh} +s \im -\mxrates \big) &= s \vup,\label{e3.eq1 for trasform vector}\\
    \vltw[0] \uv &= 1, \label{e3.eq2 for transform vector}
  \end{align}
  where $\vup=\left[u_1,\dots,u_N\right]$ is a vector with $N$ unknown parameters that needs to be determined.
\end{theorem}

Note that the above theorem is similar to Theorem~3.1 in \cite{adan03} and so does its proof. Therefore, we omit here the proof and we refer the reader to Theorem~3.1 of \cite{adan03} for more details.

\begin{remark}\label{r3.waiting times}
  Let \vweights be a column vector of dimension $N$, such that $\vweights = \mxrates^{-1} \vect{D}^{(2)} \uv/\boldsymbol\pi \mxrates^{-1} \vect{D}^{(2)} \uv$. Then, it can easily be verified that $\ltd = \vltw \vweights$ is the Laplace transform of the waiting time of a real customer. If, however, $\vweights = \uv$, then $\vltw \uv$ is the Laplace transform of the workload just prior to an arrival of a customer. Thus, for the study of our system it is sufficient to determine the transform vector \vltw.
\end{remark}

If $\det\big(\ltm{\mxh} +s\im -\mxrates\big)$ denotes the determinant of the square matrix $\ltm{\mxh} +s\im -\mxrates$, then for the determination of the unknown vector \vup, we have the following theorem.
\begin{theorem}\label{t3.solution for vector u}
  The next two statements hold:
  \begin{enumerate}
    \item The equation $\det\big(\ltm{\mxh} +s\im -\mxrates\big)=0$ has exactly $N$ solutions $\rootp{1},\dots,\rootp{N}$, with $\rootp{1}=0$ and $\Re(\rootp{i})>0$ for $i=2,\dots,N$.\label{part1}
    \item Suppose that the stability condition \eqref{e3.stability condition} is satisfied and that the above mentioned $N-1$ solutions $\rootp{2},\dots,\rootp{N}$ are distinct. Let $\va[i]$ be a non\-/zero column vector satisfying
        \begin{equation*}
          \big(\ltm[\rootp{i}]{\mxh} +\rootp{i} \im -\mxrates \big)\va[i] = 0, \quad i=2,\dots,N.
        \end{equation*}
        Then \vup is given by the unique solution to the following $N$ linear equations:
        \begin{align}
          \vup \mxrates^{-1} \uv &= \boldsymbol\pi\left(\mxrates^{-1} - \mxmeans\right) \uv, \label{e3.equation 1 vector u}\\
          \vup \va[i] &= 0, \qquad i=2,\dots,N.\label{e3.equation 2 vector u}
        \end{align}
  \end{enumerate}
\end{theorem}

Again, Theorem~\ref{t3.solution for vector u} is similar to Theorems~3.2 \& 3.3 in \cite{adan03}, and therefore, its proof is omitted.

Theorem~\ref{t3.solution for vector u} on one hand provides us with an algorithm to calculate the vector \vup and on the other hand it guarantees that all elements of the transform vector \vltw are well\-/defined on the positive half\-/plane. To understand the latter remark observe the following. For simplicity, we set
\begin{equation}\label{e3.matrix E}
  \ltm{\mxe} = \ltm{\mxh} +s\im -\mxrates.
\end{equation}
Let \ltm{\mxadj} be the adjoint matrix of \ltm{\mxe}, so $\ltm{\mxe} \cdot \ltm{\mxadj} = \det \ltm{\mxe} \im$. Post\-/multiplying Eq.~\eqref{e3.eq1 for trasform vector} with \ltm{\mxadj}, we have that
$\vltw \det \ltm{\mxe} = s \vup \ltm{\mxadj}$, and consequently
\begin{equation}\label{e3.Laplace transform waiting time as a fraction}
  \vltw = \frac{1}{\det \ltm{\mxe}} s \vup \ltm{\mxadj}.
\end{equation}
The first statement of Theorem~\ref{t3.solution for vector u} says that the determinant $\det \ltm{\mxe}$ has the factors $s-\rootp{i}$, $i=1,\dots,N$, in its expression. This means that the transform vector \vltw has $N$ potential singularities on the positive half plane, as the determinant appears at the denominator. However, the second statement of Theorem~\ref{t3.solution for vector u} explains that the vector \vup is such that these problematic factors are canceled out.

Observe that Theorem~\ref{t3.solution for vector u} does not give us any information about the form of the elements of the transform vector \vltw, which is the stepping stone for the construction of our approximations. For this reason, we proceed by finding an analytic expression for the aforementioned elements. It is apparent from Eq.~\eqref{e3.Laplace transform waiting time as a fraction} that for the evaluation of \vltw we only need $\det \ltm{\mxe}$ and the adjoint matrix \ltm{\mxadj}. For the determination of these quantities, we introduce the following notation:
\begin{itemize}
  \item As before, we denote the set of all states of the Markov process $\{J_t\}$ as $\n =\{1,\dots,N\}$ .
  \item If $S\subset \Omega$, for some set $\Omega \subset \n$, then $S^c$ is the complementary set of $S$ with respect to $\Omega$. Observe that all subset relations will be used locally and that the symbol ``$\subset$'' does not imply strict subsets. The number of elements in a set $S$ is denoted as $\left| S\right|$.
  \item For a subset $S$ of \n we define $\las{S} = \prod_{i\in S}\lambda_i$ and $\fun{\zeta^{S}} = \prod_{i\in S}(s-\lambda_i)$. We also define $\las{\emptyset} = \fun{\zeta^{\emptyset}} = 1$.
  \item Suppose that $U,W\subset \n$ and that \vect{A} is a square matrix of dimension $N$. Then $\vect{A}^W_U$ is the submatrix of \vect{A} if we keep the rows in $U$ and the columns in $W$. Whenever the notation becomes very complicated, to avoid any confusion with the indices, we will denote the $i$th column and row of matrix \vect{A} with $\vect{A}_{\bullet i}$ and $\vect{A}_{i\bullet}$, respectively. We also define $\det\vect{A}_\emptyset^\emptyset=1$.
  \item Suppose that $S$ is a subset of $\Omega$, for some set $\Omega \subset \n$, and that it follows some properties, i.e.\ ``Property~1", etc. If we want to sum with respect to $S$, then we write under the symbol of summation first $S\subset \Omega$, followed by the properties. Namely, we write $\sum_{\substack{S\subset \Omega\\ \text{Property~1} \\ \text{etc}}}$. In some cases, to avoid lengthy expressions we will write instead of $\sum_{\substack{S\subset \Omega\\ \text{Properties of } S}}\sum_{\substack{R\subset \Omega_1\\ \text{Properties of } R}}$ the double sum $\sum_{\substack{S\subset \Omega\\ \text{Properties of } S; \\ R\subset \Omega_1\\ \text{Properties of } R}}$, where $R$ is a subset of $\Omega_1$, for some set $\Omega_1 \subset \n$. We apply the same rule also for multiple sums.
  \item Suppose that \vect{A} and \vect{B} are two square matrices of dimension $N$, and that $U$ and $W$ are two disjoint subsets of \n. For all $\Omega\subset\n$, we use the notation $\vect{A}^U_\Omega\Join \vect{B}^W_\Omega$ for the matrix that has  columns the union of the columns $V$ of matrix \vect{A} and the columns $W$ of matrix \vect{B}, ordered according to the index set $U\cup W$; e.g.\ if $\Omega=\n=\{1,\dots,5\}$, $U=\{1,2,4\}$, and $W=\{3,5\}$, then $\vect{A}^{\{1,2,4\}}_{\n}\Join~\vect{B}^{\{3,5\}}_{\n} = (\vect{A}_{\bullet 1},\vect{A}_{\bullet 2},\vect{B}_{\bullet 3},\vect{A}_{\bullet 4},\vect{B}_{\bullet 5})$.
\end{itemize}

Using the above notation, we proceed with refining the desired quantities. More precisely, we first find $\det \ltm{\mxe}$, then the adjoint matrix \ltm{\mxadj}, and finally the vector $s\vup \ltm{\mxadj}$ that appears in the numerator of the transform vector \vltw (see Eq.~\eqref{e3.Laplace transform waiting time as a fraction}). Combining these results, one can easily derive \vltw. We start by finding the determinant of matrix \ltm{\mxe} (see Eq.~\eqref{e3.matrix E}).

\begin{theorem}\label{t3.determinant matrix E}
  The determinant of matrix \ltm{\mxe} can be explicitly calculated as follows:
  \begin{align*}
    \det \ltm{\mxe} =& \sum_{S\subset \n}\las{S}\fun{\zeta^{S^c}}\det \big(\mxprobsdummy\circ\mxtrans\big)^S_S \\
                    +& \sum_{k=1}^N \ltgs[k] \sum_{\substack{\Gamma\subset \n\\  \mid \Gamma\mid =k}}  \sum_{\substack{S\subset \n\\ S\supset \Gamma}}
                      \las{S}  \fun{\zeta^{S^c}} \\
                    & \times \det \Big( \big(\mxprobsdummy\circ\mxtrans\big)^{S\setminus \Gamma}_S  \Join \big(\mxprobsreal\circ\mxtrans\big)^\Gamma_S\Big).
  \end{align*}
\end{theorem}
\begin{proof}
    See Appendix~\ref{appendix B}.
\end{proof}

Observe that the determinant $\det \ltm{\mxe}$ is an at most $N$ degree polynomial with respect to the LST of the service time distribution \ltgs of a real customer. Moreover, the coefficients of this polynomial are all polynomials with respect to $s$. Therefore, in case \ltgs is a rational function in $s$, then $\det \ltm{\mxe}$ is also a rational function in $s$ and its eigenvalues can be easily calculated. Furthermore, the subset $\Gamma$ of \n that appears in the second summand has at least one element, thus in the formula of $\det \ltm{\mxe}$ it always holds that $\Gamma \ne \emptyset$.

\begin{toy}
  The matrix \ltm{\mxe} has elements $\ltm{\mxe[ii]}= s-\lambda$, $i=1,2$, $\ltm{\mxe[12]}= \lambda$, and $\ltm{\mxe[21]}= \lambda \ltgs$. We calculate its determinant using Theorem~\ref{t3.determinant matrix E}. It holds that $\det \big(\mxprobsdummy\circ~\mxtrans\big)^S_S=0$ for all subsets $S$ of \n, except for $S=\emptyset$. Since $\Gamma \ne \emptyset$, it is evident that $\det \Big( \big(\mxprobsdummy\circ~\mxtrans\big)^{S\setminus \Gamma}_S \allowbreak\Join \big(\mxprobsreal\circ\mxtrans\big)^\Gamma_S \Big)\neq 0$ only for $\Gamma=\{1\}$ and $S=\n$, because the $1$st column of the matrix \mxprobsdummy and the $2$nd column of the matrix \mxprobsreal are zero. Combining all these we obtain
  \begin{align*}
   \det \ltm{\mxe} =& \las{\emptyset}\fun{\zeta^{\n}}\det \big(\mxprobsdummy\circ\mxtrans\big)^\emptyset_\emptyset
                    + \ltgs \las{\n} \fun{\zeta^{\emptyset}} \\
                    & \times \det \Big( \big(\mxprobsdummy\circ\mxtrans\big)^{\{2\}}_{\n}  \Join \big(\mxprobsreal\circ~\mxtrans\big)^{\{1\}}_{\n} \Big)\\
                   =& (s-\lambda)^2 - \lambda^2 \ltgs.
  \end{align*}
\end{toy}

In a similar manner, we find the explicit form of the adjoint matrix \ltm{\mxadj} in the following theorem.

\begin{theorem}\label{t3.adjoint matrix of E}
  The adjoint matrix \ltm{\mxadj} has elements
  \begin{small}
  \begin{align*}
    &\ltm{\mxadj[ij]}\\
                    &= \begin{cases}
                            \sum_{k=0}^{N-1} \ltgs[k] \sum_{\substack{\Gamma\subset \n\setminus\{i\}\\  \mid \Gamma\mid =k; \\ S\subset \n\setminus\{i\}\\ S\supset \Gamma}}
                                \las{S}  \fun{\zeta^{S^c}} \\
                            \times \det \Big( \big(\mxprobsdummy\circ\mxtrans\big)^{S\setminus \Gamma}_S \Join \big(\mxprobsreal\circ\mxtrans\big)^\Gamma_S  \Big),
                                            & i=j,\\
                            (-1)^{i+j} \sum_{k=1}^{N-1} \ltgs[k] \sum_{\substack{\Gamma \subset \n\setminus\{i,j\}\\  \mid \Gamma\mid =k-1}} \\
                            \times \sum_{\substack{S\subset \n\setminus\{i,j\}\\ S\supset \Gamma;\\ R\subset S\cap T_{ij}}} (-1)^{\mid R\mid} \las{S\cup \{j\}} \fun{\zeta^{S^c}} \\
                            \times \det \Big( \big(\mxprobsdummy\circ\mxtrans\big)^{S\setminus \Gamma}_{S\cup \{i\}} \Join \big(\mxprobsreal\circ\mxtrans\big)^{\Gamma \cup \{j\}}_{S\cup \{i\}}  \Big)\\
                            + (-1)^{i+j} \sum_{k=0}^{N-2} \ltgs[k] \sum_{\substack{\Gamma \subset \n\setminus\{i,j\}\\  \mid \Gamma\mid =k}} \\
                            \times \sum_{\substack{S\subset \n\setminus\{i,j\}\\ S\supset \Gamma;\\ R\subset S\cap T_{ij}}} (-1)^{\mid R\mid} \las{S\cup \{j\}} \fun{\zeta^{S^c}} \\
                            \times \det \Big( \big(\mxprobsdummy\circ\mxtrans\big)^{(S\setminus \Gamma)\cup \{j\}}_{S\cup \{i\}} \Join \big(\mxprobsreal\circ\mxtrans\big)^{\Gamma}_{S\cup \{i\}}  \Big),
                                 & i\neq j,
                      \end{cases}
  \end{align*}
  \end{small}
  where $m_{ij}=\min\{i,j\}$, $M_{ij}=\max\{i,j\}$, and $T_{ij}=\{m_{ij}+1,\dots,M_{ij}-1\}$.
\end{theorem}
\begin{proof}
    See Appendix~\ref{appendix B}.
\end{proof}

%\\                            \times \comqus{(S\cup \{i\})\setminus \Gamma}

The adjoint matrix \ltm{\mxadj} is equal to the transpose of the cofactor matrix of \ltm{\mxe}. Therefore, similarly to $\det \ltm{\mxe}$, each element of \ltm{\mxadj} is an at most $N-1$ degree polynomial with respect to \ltgs. This observation explains also the similarity between the formula of $\det \ltm{\mxe}$ and the diagonal elements of \ltm{\mxadj}.

\begin{toy}
  Using the same arguments as for the evaluation of the determinant, we have for the adjoint matrix
  \begin{align*}
    \ltm{\mxadj[ii]} =& \ltgs[0]\las{\emptyset} \fun{\zeta^{\n \setminus \{i\}}} \\
                      & \times \det \Big( \big(\mxprobsdummy\circ\mxtrans\big)^\emptyset_\emptyset \Join \big(\mxprobsreal\circ\mxtrans\big)^\emptyset_\emptyset  \Big) \\
                     =& s-\lambda, \qquad \qquad \qquad i=1,2, \notag \\
    \ltm{\mxadj[12]} =& (-1)^{1+2} (-1)^{\mid \emptyset\mid} \las{\emptyset \cup \{2\}} \fun{\zeta^{\emptyset}} \\
                      & \times \det \Big( \big(\mxprobsdummy\circ\mxtrans\big)^{\{2\}}_{\{1\}} \Join \big(\mxprobsreal\circ\mxtrans\big)^{\emptyset}_{\{1\}}  \Big)\\
                     =& -\lambda \notag \\
    \ltm{\mxadj[21]} =& (-1)^{2+1} \ltgs  (-1)^{\mid \emptyset\mid} \las{\emptyset \cup \{1\}} \fun{\zeta^{\emptyset}}\\
                      & \times \det \Big( \big(\mxprobsdummy\circ\mxtrans\big)^{\emptyset}_{\{2\}} \Join \big(\mxprobsreal\circ\mxtrans\big)^{\{1\}}_{\{2\}}  \Big)\\
                     =& -\lambda \ltgs. \notag
  \end{align*}
\end{toy}

Observe that the elements of the transform vector \vltw are defined as $\ltwc{i}= s \vup \ltm{\mxadj}\uv_i / \det\ltm{\mxe}$ (see Eq.~\eqref{e3.Laplace transform waiting time as a fraction}), where $\uv_i$ is a column vector with element equal to 1 in position $i$ and all other elements zero. The outcome of $s\vup\ltm{\mxadj}\uv_i$ is the inner product of the vector $s\vup$ with the $i$th column of matrix $\ltm{\mxadj}$. Therefore, we have the following theorem.
\begin{theorem}\label{t3.numerator LT workload state i}
  The numerator of the $i$th element of the transform vector \vltw takes the form
  \begin{small}
  \begin{align*}
  s\vup& \ltm{\mxadj}\uv_i = su_i \sum_{k=0}^{N-1} \ltgs[k] \sum_{\substack{\Gamma\subset \n\setminus\{i\}\\  \mid \Gamma\mid =k; \\ S\subset \n\setminus\{i\}\\ S\supset \Gamma}}
                                \las{S}  \fun{\zeta^{S^c}} \\
                        & \times \det \Big( \big(\mxprobsdummy\circ\mxtrans\big)^{S\setminus \Gamma}_S \Join \big(\mxprobsreal\circ\mxtrans\big)^\Gamma_S  \Big)\\
                       +& s\sum_{\substack{l=1\\ l\neq i}}^N u_l (-1)^{l+i} \sum_{k=1}^{N-1} \ltgs[k] \sum_{\substack{\Gamma\subset \n\setminus\{l,i\}\\  \mid \Gamma\mid =k-1}}
                          \sum_{\substack{S\subset \n\setminus\{l,i\}\\ S\supset \Gamma;\\ R\subset S\cap T_{li}}} (-1)^{\mid R\mid} \\
                        & \times \las{S\cup \{i\}} \fun{\zeta^{S^c}}
                          \det \Big( \big(\mxprobsdummy\circ\mxtrans\big)^{S\setminus \Gamma}_{S\cup \{l\}} \Join \big(\mxprobsreal\circ\mxtrans\big)^{\Gamma \cup \{i\}}_{S\cup \{l\}}  \Big)
                         \notag\\
                       +& s\sum_{\substack{l=1\\ l\neq i}}^N u_l (-1)^{l+i} \sum_{k=0}^{N-2} \ltgs[k] \sum_{\substack{\Gamma\subset \n\setminus\{l,i\}\\  \mid \Gamma\mid =k}}
                          \sum_{\substack{S\subset \n\setminus\{l,i\}\\ S\supset \Gamma;\\ R\subset S\cap T_{li}}} (-1)^{\mid R\mid} \las{S\cup \{i\}} \\
                        & \times \fun{\zeta^{S^c}}
                          \det \Big( \big(\mxprobsdummy\circ\mxtrans\big)^{(S\setminus \Gamma)\cup \{i\}}_{S\cup \{l\}} \Join \big(\mxprobsreal\circ\mxtrans\big)^{\Gamma}_{S\cup \{l\}}  \Big).
  \end{align*}
  \end{small}
\end{theorem}
\begin{proof}
  See Appendix~\ref{appendix B}.
\end{proof}

Combining now the results of the Theorems~\ref{t3.determinant matrix E} and \ref{t3.numerator LT workload state i} by using Eq.~\eqref{e3.Laplace transform waiting time as a fraction}, one can find the transform vector \vltw.

\begin{toy}
  For each state we have
  \begin{align*}
    \vup \ltm{\mxadj} \uv_1 =& u_1 \ltm{\mxadj[11]} + u_2 \ltm{\mxadj[21]} = u_1 (s-\lambda) - u_2 \lambda \ltgs, \notag \\
    \vup \ltm{\mxadj} \uv_2 =& u_1 \ltm{\mxadj[12]} + u_2 \ltm{\mxadj[22]} = - u_1 \lambda + u_2 (s-\lambda).  \notag
  \end{align*}
  The transform vector \vltw is then
 \begin{equation*}
   \vltw =  \left[\frac{s u_1 (s-\lambda) - s u_2 \lambda \ltgs}{(s-\lambda)^2 - \lambda^2 \ltgs}, \frac{-s u_1 \lambda + s u_2 (s-\lambda)}{(s-\lambda)^2 - \lambda^2 \ltgs} \right].
 \end{equation*}
\end{toy}

The following remark connects the system of equations that is required for the evaluation of \vup, which was introduced in Theorem~\ref{t3.solution for vector u}, to the adjoint matrix $\ltm{\mxadj}$.

\begin{remark}\label{r3.eigenvectors a}
  The second statement of Theorem~\ref{t3.solution for vector u} practically says that each \rootp{i}, $i=2,\dots,N$, is a simple eigenvalue of the matrix $\ltm{\mxh} +s \im -\mxrates$.  Therefore, the column vector $\va[i]$ belongs to the null space of the matrix $\ltm[\rootp{i}]{\mxh} +\rootp{i} \im -\mxrates$. Combining the results of Theorem~\ref{t3.eigenvectors that correspond to the eigenvalue zero}, Remark~\ref{r3.proportional eigenvectors} and Corollary~\ref{c3.eigenvectors of E} (see Appendix~\ref{appendix A}), which provide some general results with respect to the form of the null space of a singular matrix, without loss of generality we can assume that the vector $\va[i]$ is any non\-/zero column of the matrix $\ltm[\rootp{i}]{\mxadj}$. Namely, if the $m$th column of $\ltm[\rootp{i}]{\mxadj}$ is such a column, then
  \begin{equation}\label{e3.column vectors a}
    \va[i] := \fun[\rootp{i}]{\va[i]} = \big(\ltm[\rootp{i}]{\mxadj}\big)^{\{m\}}_{\n}, \qquad i=2,\dots,N.
  \end{equation}
  This observation is very useful, because it allows us to calculate in a straightforward way the desired system of equations and find closed form expressions for the vector \vup.  In addition, since the vectors $\va[i]$, $i=2,\dots,N$, are matrix functions evaluated at the point $s=\rootp{i}$ we define the derivative of each $\va[i]$ as
  \begin{equation*}
    \dm[1]{\va[i]} = \left.\frac{d}{ds}\fun{\va[i]}\right|_{s=\rootp{i}}, \qquad i=2,\dots,N.
  \end{equation*}
  The usefulness of the latter definition will be apparent in Section~\ref{s3.parameters of the perturbed system}, where we provide an extension of Theorem~\ref{t3.solution for vector u} that helps us to calculate our approximations.
\end{remark}

\begin{toy}
  If \rootp{2} is the only positive (and real) root of the equation $\det \ltm{\mxe}  =0$, the vector \vup satisfies the system of equations \eqref{e3.equation 1 vector u}--\eqref{e3.equation 2 vector u}
  \begin{align*}
    \frac1\lambda u_1 + \frac1\lambda u_2 &= \frac1{\lambda}-\frac{\mean}{2},\\
    -\lambda u_1 + (\rootp{2}-\lambda) u_2 &=0,
  \end{align*}
  where for the derivation of the second equation we used the second column of the matrix \ltm{\mxadj}. Namely, we used $\va[2] = \big(\ltm[\rootp{2}]{\mxadj}\big)^{\{2\}}_{\n}$. It is easy to verify that the solution to the above system is given by
  \begin{equation*}
    \vup = \left( \Big(1-\frac{\lambda}{\rootp{2}}\Big)\Big(1-\frac{\lambda \mean}{2}\Big)  ,\  \frac{\lambda}{\rootp{2}}\Big(1-\frac{\lambda \mean}{2}\Big) \right).
  \end{equation*}
\end{toy}

Although Theorems~\ref{t3.adjoint matrix of E} and \ref{t3.numerator LT workload state i} provide explicit expressions for the transform vector, they may not be practical in cases where the LST of the service time distribution of a real customer \ltgs, which is involved in the formulas, does not have a closed form; i.e.\ Pareto distribution. In such cases, one would have to either consort to a numerical evaluation of \ltgs or approximate the transform vector \vltw in some other fashion. This paper focuses on the latter approach, which we work out in detail in the following section by taking as starting point a mixture model for the service time distribution of a real customer.

%\subsection{Perturbation of the service time distribution}\label{s3.perturbation of the service time distribution}

\subsection{Construction of the corrected phase\-/type approximations}\label{s3.Construction of the corrected phase-type approximations}
We assume now that the service time distribution of a real customer is \gsdmix, which was defined in Eq.~\eqref{e3.mixture service time distribution} as a mixture of a phase\-/type distribution and a heavy\-/tailed one. We will eventually show that the workload can be written also as a mixture, in the sense that we can identify the workload of a model with purely phase\-/type service times and some additional terms that involve the heavy\-/tailed service times.  As a result, in order to derive our approximations, we first need to compute the workload in a MAP/PH/1 queue and afterwards use this as a base to further develop our approximations involving a heavy\-/tailed component. In the sequel, we give a more detailed description of our technique.

In terms of Laplace transforms we get for our mixture service time distribution $\ltgsmix = (1-\epsilon)\ltpts + \epsilon \lthts$. As observed in Section~\ref{s3.preliminaries}, when the service time distribution of a real customer is of phase type, then the determinant $\det \ltm{\mxe}$ and the elements of the adjoint matrix \ltm{\mxadj} are all rational functions in $s$. Therefore, after the cancelation of the problematic factors $s-\rootp{i}$, $i=1,\dots,N$ (see the analysis below Theorem~\ref{t3.solution for vector u}), the elements of the transform vector \vltw are also rational functions in $s$ and they can easily be inverted to find the workload distribution. %Thus, to overcome the difficulty that occurs from the LST \lthts, a phase\-/type approximation would suggest to first remove the heavy\-/tailed customers, and then find an explicit phase\-/type representation for the transform vector \vltw of the resulting simpler MAP/PH/1 queue.

Note now that the LST of the service time distribution of a real customer \ltgsmix can be written in the following way:
\begin{equation*}
  \ltgsmix = \ltpts + \epsilon \big(\lthts - \ltpts\big).
\end{equation*}
In this formula, \ltgsmix can be seen as perturbation of the phase\-/type distribution \ltpts by the term $\epsilon \big(\lthts - \ltpts\big)$. The index $\epsilon$ is interpreted as the perturbation parameter and it used for all parameters of the system that depend on it. By setting $\lthts \equiv \ltpts$\footnote{In other words, we assume that all of the customers come from the same phase\-/type distribution or equivalently that we replace all the heavy\-/tailed customers with phase\-/type ones.} in the formula, one can find with $\ltgsmix = \ltpts$ the workload of a simpler MAP/PH/1 queue, by specializing the formulas of Section~\ref{s3.preliminaries} to phase\-/type service times. As a next step, we find all the parameters of the mixture model as perturbation of the simpler phase\-/type model, which we use as base. Then, we write the workload of the mixture model in a series expansion in $\epsilon$, where the constant term is the workload of the MAP/PH/1 queue we used as base and all other terms contain the heavy\-/tailed service times.
%Since the latter service time distribution is a phase\-/type approximation of the original mixture service time distribution, the resulting workload will be also a phase\-/type approximation of the workload of the mixture model.%As stated earlier, the idea is to compute first the workload in a MAP/PH/1 queue and then built upon it to derive our approximations. Next, we explain how this interpretation of the formula can lead, with the aid of perturbation analysis, to a phase\-/type base models and consequently to the approximations for the workload of the mixture model.

%The first step to find a phase\-/type approximation for the workload is to remove all the heavy\-/tailed customers from the system. We will do this here, by replacing them with phase\-/type customers. In other words, we assume that all of the customers come from the same phase\-/type distribution and we set $\lthts \equiv \ltpts$.

We define our approximation by taking the first two terms of the aforementioned series, namely the up to $\epsilon$-order terms. We call this approximation {\it corrected replace approximation}. The characterization ``corrected'' comes from the fact that the $\epsilon$-order term corrects the tail behavior of the constant term, which as a phase\-/type approximation of the workload is incapable of capturing the correct tail behavior. Finally, the characterization ``replace" is due to the phase\-/type base model we used. We give analytically all the steps to derive the corrected replace approximation in Section~\ref{s3.corrected replace approximation}.

\section{Corrected replace approximation}\label{s3.corrected replace approximation}

In this section, we construct the corrected replace approximation. First, we calculate the workload for the phase\-/type model that appears when we replace all the heavy\-/tailed customers with phase\-/type ones in Section~\ref{s3.replace base model}; i.e.\ we specialize the results of Section~\ref{s3.preliminaries} to phase\-/type service times. Later, in Section~\ref{s3.parameters of the perturbed system}, we calculate the parameters of the mixture model with service time distribution \ltgsmix given by Eq.~\eqref{e3.mixture service time distribution} as perturbation of the parameters of the corresponding phase\-/type model, with perturbation parameter $\epsilon$. In Section~\ref{s3.workload distribution of the perturbed model}, we find a series expansion in $\epsilon$ of the workload in the mixture model with constant term the workload in the phase\-/type base model and all higher terms involving the heavy\-/tailed services. Finally, in Section~\ref{s3.corrected replace approximation}, we construct the corrected replace approximation by keeping only the first two terms of the aforementioned series. We start in the next section with the analysis of the replace base model; i.e.\ the one containing only phase\-/type service times.

\subsection{Replace base model}\label{s3.replace base model}
When we replace the heavy\-/tailed customers with phase\-/type ones, we consider the service time distribution $\ltgsmix \allowbreak = \ltpts$ for our phase\-/type base model. Observe that this service time distribution is independent of the parameter $\epsilon$, and so will be all the other parameters of this simpler model. Thus, from a mathematical point of view, the action of replacing the heavy\-/tailed claim sizes with phase\-/type ones is equivalent to setting $\epsilon=0$ in the mixture model.

To avoid overloading the notation, we omit the subscript ``0" (which is a consequence of the fact that $\epsilon=0$) from the parameters of the replace phase\-/type model and we assume that the service time distribution of a real customer is some phase\-/type distribution with LST $\ltgs:=\ltpts = \fun{q}/\fun{p}$, where \fun{q} and \fun{p} are appropriate polynomials without common roots. The degree of \fun{p} is $M$, and without loss of generality, we choose the coefficient of its highest order term to be equal to 1. Finally, the degree of the polynomial \fun{q} is less than or equal to $M-1$. Define
\begin{align}\label{e3.definition kappa}
  K &= \max_{k\neq 0}\left\{ \max_{\Gamma \subset \n} \left\{ \rank \Big( \big(\mxprobsdummy\circ\mxtrans\big)^{S\setminus \Gamma}_S \Join \right. \right. \notag\\
     & \left. \left.\big(\mxprobsreal\circ\mxtrans\big)^\Gamma_S\Big)\right\} : k=\mid\Gamma\mid, \text{ and }\Gamma \subset S \subset \n \right\}.
\end{align}
Then, the following result holds.

%Moreover, since $\ltpts[0] =1$, the constant terms of the polynomial are also equal.

\begin{proposition}\label{p3.LT replace workload state i}
  There exist \rootnden{j}, with $\Re(\rootnden{j})>0$, $j=1,\dots,rM$, and for each state $i \in \n$, there exist \rootnnum{i,j} with $\Re(\rootnnum{i,j})>0$, $j=1,\dots,rM$, such that the Laplace transform \ltwc{i} takes the form
  \begin{equation*}
    \ltwc{i} = \frac{u_i \prod_{j=1}^{rM} (s+\rootnnum{i,j})}{\prod_{j=1}^{rM} (s+\rootnden{j})},
  \end{equation*}
  where $u_i$ is the $i$th element of the vector \vup that can be calculated according to Theorem~\ref{t3.solution for vector u} with the LST of the service times being equal to \ltpts, and $r$ is some positive integer less than or equal to $K$ defined by \eqref{e3.definition kappa}.
\end{proposition}
\begin{proof}
  See Appendix~\ref{appendix B}.
\end{proof}

The formula of \ltwc{i} is a rational function that corresponds to a phase\-/type distribution. Applying Laplace inversion to \ltwc{i}, we can find the exact tail probabilities of the workload prior to an arrival of a customer in state $i$; namely we can find $\pr(\workload[i] >t)$.

\begin{toy}
  Here, instead of calculating the transform vector \vltw for phase\-/type customers, we deal with the Laplace transform of the waiting time of a real customer in the queue; namely $\ltd = \vltw \vweights$ (see Remark~\ref{r3.waiting times}) with $\vweights^T = (0,2)$, where superscript $T$ denotes the transpose of a vector or a matrix. In our example, $K=1$ and consequently, $r=1$. Thus, \ltd under phase\-/type service times is
  \begin{align*}
  \ltd =& 2\ \ltwc{2}
       = 2 \frac{s^2 u_2 - s \lambda (u_1+u_2)}{(s-\lambda)^2 - \lambda^2 \ltpts}\\
       =& 2 \frac{s^2 \fun{p} u_2 - s \fun{p} \lambda (u_1+u_2)}{(s-\lambda)^2 \fun{p} - \lambda^2 \fun{q}}.
  \intertext{Observe that both the numerator and the denominator of \ltd are polynomials of degree $M+2$. Moreover, Theorem~\ref{t3.solution for vector u} guarantees that 0 and \rootp{2} are common roots of them. If $-\rootnnum{j}$ and $-\rootnden{j}$, $j=1,\dots,M$, $\Re(\rootnden{j}),\Re(\rootnnum{j})>0$, are the remaining roots of the numerator and the denominator, respectively, the Laplace transform of the waiting time can be written as}
  \ltd =& \frac{2 u_2  s (s-\rootp{2}) \prod_{j=1}^{M} (s+\rootnnum{j})}{s (s-\rootp{2}) \prod_{j=1}^{M} (s+\rootnden{j})}
       = \frac{2 u_2 \prod_{j=1}^{M} (s+\rootnnum{j})}{\prod_{j=1}^{M} (s+\rootnden{j})}.
  \end{align*}
\end{toy}

As pointed out in Section~\ref{s3.Construction of the corrected phase-type approximations}, the LST of the service time distribution \ltgsmix (see Eq.~\eqref{e3.mixture service time distribution}) can be seen as perturbation of \ltpts by the term $ \epsilon \big(\lthts - \ltpts\big)$. In the next section we write the parameters of the mixture model as perturbation of the parameters of the replace base model.

\subsection{Perturbation of the parameters of the replace base model}\label{s3.parameters of the perturbed system}
In order to find the workload in the mixture model as a series expansion in $\epsilon$ with constant term the workload in the replace base model, we apply perturbation analysis to the parameters of the mixture model that depend on $\epsilon$. Thus, we first check which of the parameters in the mixture model depend on $\epsilon$ and then we represent them as perturbation of the parameters of the replace base model.

Since the matrices \mxtrans, \mxprobsdummy, \mxprobsreal, and \mxrates (see Section~\ref{s3.preliminaries}) depend only on the arrival process, they are invariant under any perturbation of the service time distribution. However, the matrix \ltmgsmix,  and consequently \ltm{\mxhmix} change, and so does the stability condition (see Eqs.~\eqref{e3.definition l-t general service times matrix}--\eqref{e3.stability condition}). Let now \ltptes and \lthtes be the LSTs of the stationary\-/excess service time distributions \ptesd and \htesd, and \mean[p] and \mean[h] be the finite means of the phase\-/type and heavy\-/tailed service times, respectively. Then, we obtain
\begin{align}
  \ltmgsmix     =& \ltmgs + \epsilon s \big(\mean[p]\ltptes - \mean[h]\lthtes\big) \mxprobsreal, \label{e3.definition l-t general service times matrix mixture model}
  \intertext{and}
  \ltm{\mxhmix} =& \ltmgsmix \circ \mxtrans \mxrates \notag \\
                =& \ltm{\mxh} + \epsilon s \big(\mean[p]\ltptes - \mean[h]\lthtes\big)\mxprobsreal \circ\mxtrans \mxrates.\label{e3.definition matrix h mixture model}
\end{align}
Finally, the stability condition takes the form
\begin{equation}\label{e3.stability condition mixture model}
  \boldsymbol\pi (\mxrates^{-1} - \mxmeansmix)\uv >0,
\end{equation}
where $ \mxmeansmix = \mxmeans + \epsilon s \big(\mean[h] - \mean[p]\big) \mxprobsreal\circ\mxtrans$.

Under the stability condition \eqref{e3.stability condition mixture model}, Theorem~\ref{t3.equations for the transform vector} holds for the transform vector \vltwmix, for some row vector \vupmix. More precisely, there exist a unique vector \vupmix such that the transform vector \vltwmix satisfies the system of equations:
\begin{align}
 \vltwmix \big(\ltm{\mxhmix} +s \im -\mxrates \big) &= s \vupmix, \label{e3.eq1 for trasform vector mixture model}\\
 \vltwmix[0] \uv &= 1, \label{e3.eq2 for trasform vector mixture model}
\end{align}
where the vector \vupmix is calculated according to Theorem~\ref{t3.solution for vector u}.

Recall that the evaluation of \vupmix goes through the evaluation of the positive eigenvalues of the matrix
\begin{align}
  \ltm{\mxemix} =& \ltm{\mxhmix} +s \im -\mxrates  \notag\\
                =& \ltm{\mxe} + \epsilon s \big(\mean[p]\ltptes - \mean[h]\lthtes\big)\mxprobsreal \circ\mxtrans \mxrates.\label{e3.definition perturbed matrix e}
\end{align}
Observe that the above representation of the matrix \ltm{\mxemix} is a linear perturbation in $\epsilon$ of the matrix \ltm{\mxe} of the base model. Thus, according to results on perturbation of analytic matrix functions \cite{deteran11,lancaster03}, we have that the positive eigenvalues of the matrix \ltm{\mxemix} and their corresponding eigenvectors are analytic functions in $\epsilon$. Consequently, one can find a series representation in $\epsilon$ for all the involved quantities that are needed for the evaluation of the vector \vupmix (see Theorem~\ref{t3.solution for vector u}). By using these parameters, we can find a complete series representation for the transform vector \vltwmix and by applying Laplace inversion to each term of this series we can find a formal expression for the workload that is a series expansion in $\epsilon$. As we stated earlier, we only need the first two terms of the latter series to define the corrected replace approximation. Therefore, in our analysis, we keep only the terms up to order $\epsilon$ of each involved perturbed parameter.

In the next theorem, we provide an algorithm to calculate the first order approximation in $\epsilon$ of the vector \vupmix, given that we have already calculated the vector \vup of the replace base model, by specializing Theorem~\ref{t3.solution for vector u} to phase\-/type service times. We denote by \um the square matrix of appropriate dimensions with all its elements equal to one.

\begin{theorem}\label{t3.vector u in the mixture model}
  Let \vup be the unique solution to the Eqs.~\eqref{e3.equation 1 vector u}--\eqref{e3.equation 2 vector u} for the replace base model. If the roots $\rootp{2},\dots,\rootp{N}$ of\/ $\det\big(\ltm{\mxh} +s\im -\mxrates\big)=0$ with positive real part are simple, then
  \begin{enumerate}
    \item the equation $\det\big( \ltm{\mxhmix} +s\im -\mxrates \big) = 0$ has exactly $N$ non\-/negative solutions $\rootpmix{1},\dots,\rootpmix{N}$, with $\rootpmix{1}=0$ and $\rootpmix{i}=\rootp{i} - \epsilon \delta_i + O(\epsilon^2)$ for $i=2,\dots,N$, where
        \begin{small}
        \begin{align*}
          \delta_i:&=\delta(\rootp{i})\\
                   &= \frac{\sum_{j=1}^N \det\big(\ltm[\rootp{i}]{\mxe}_{\bullet 1},\dots,\ltm[\rootp{i}]{\mxk}_{\bullet j},\dots,
                        \ltm[\rootp{i}]{\mxe}_{\bullet N}\big)} {\sum_{j=1}^N
                        \det\big(\ltm[\rootp{i}]{\mxe}_{\bullet 1},\dots,\ltm[\rootp{i}]{\dm[1]{\mxe}}_{\bullet j},
                        \dots,\ltm[\rootp{i}]{\mxe}_{\bullet N}\big)},
        \end{align*}\label{part 1p}
        \end{small}
        and $\ltm{\mxk} = s \big(\mean[p]\ltptes-\mean[h]\lthtes\big)\mxprobsreal \circ\mxtrans \mxrates$.
    \item We set $\mxa =\big(\mxrates^{-1} \uv, \va[2],\dots,\va[N]\big)$ (see Eq.~\eqref{e3.column vectors a}) and $\vc  = \big(\boldsymbol{\pi}(\mxrates^{-1} - \mxmeans) \uv,0,\dots,0\big)$, and we assume that the stability condition \eqref{e3.stability condition mixture model} is satisfied. Then, the vector \vupmix is the unique solution to the system of $N$ linear equations
        \begin{equation}\label{e3.system of equations mixture model}
          \vupmix \big(\mxa-\epsilon \mxb + O(\epsilon^2 \um)\big)=\vc + \epsilon \vd,
        \end{equation}
        where $\mxb= \big(\vect{0},\delta_2\dm[1]{\va[2]} - \vect{k}_2,\dots,\delta_N\dm[1]{\va[N]} - \vect{k}_N\big)$ and $\vd = \big((\mean[p]-\mean[h])\boldsymbol{\pi}\mxprobsreal\circ\mxtrans\uv,0,\dots,0\big)$, with $\vect{k}_i$, $i=2,\dots,N$, being a column vector with coordinates
        \begin{small}
        \begin{align*}
          &k_{i,j} = (-1)^{m+j}\sum_{k=1}^{N-1}\det \Bigg(\Big(\ind{\ltm[\rootp{i}]{\mxe}}{\n\setminus\{j\}}{\n\setminus\{m\}}\Big)_{\bullet 1},\dots, \\
          &\Big(\ind{\ltm[\rootp{i}]{\mxk}}{\n\setminus\{j\}}{\n\setminus\{m\}}\Big)_{\bullet k},\dots,
          \Big(\ind{\ltm[\rootp{i}]{\mxe}}{\n\setminus\{j\}}{\n\setminus\{m\}}\Big)_{\bullet N-1} \Bigg),
          \quad j \in \n,
        \end{align*}
        \end{small}
        and the choice of $m$ explained in Remark~\ref{r3.eigenvectors a}.
  \end{enumerate}
\end{theorem}
\begin{proof}
  See Appendix~\ref{appendix B}.
\end{proof}

\begin{remark}\label{r3.vector k}
  When the number of states is $N=2$, the column vector $\vect{k}_2$ of Theorem~\ref{t3.vector u in the mixture model} is equal to
  \begin{align*}
    \vect{k}_2 &= \big( \ltm[\rootp{2}]{\mxk[22]}, -\ltm[\rootp{2}]{\mxk[21]} \big)^T \qquad \text{or} \\
    \vect{k}_2 &= \big( -\ltm[\rootp{2}]{\mxk[12]}, \ltm[\rootp{2}]{\mxk[11]} \big)^T
  \end{align*}
  depending on whether $m=1$ or $m=2$, respectively, where superscript $T$ denotes the transpose of a vector or matrix. The case $N=1$ has been treated earlier by the authors; see \cite{vatamidou13}.
\end{remark}

\begin{toy}
  In order to evaluate the vector \vupmix, we first need to calculate the perturbed root \rootpmix{2}, and more precisely the term $\delta_2$. Observe that in our case only the element $\ltm{\mxk[21]} = s \lambda \big(\mean[p]\ltptes-\mean[h]\lthtes\big)$ of the matrix \ltm{\mxk} is not equal to zero. Then, the numerator of $\delta_2$ becomes
  \begin{align*}
    \det &\big(\ltm[\rootp{2}]{\mxe}^{\{1\}}_{\n},\ltm[\rootp{2}]{\mxk}^{\{2\}}_{\n}\big)
    + \det\big(\ltm[\rootp{2}]{\mxk}^{\{1\}}_{\n},\ltm[\rootp{2}]{\mxe}^{\{2\}}_{\n}\big) \\
    &= - \rootp{2} \lambda^2 \big(\mean[p]\ltptes[\rootp{2}]-\mean[h]\lthtes[\rootp{2}]\big),
  \end{align*}
  and its denominator takes the form
  \begin{align*}
    \det & \big(\ltm[\rootp{2}]{\mxe}^{\{1\}}_{\n},\ltm[\rootp{2}]{\dm[1]{\mxe}}^{\{2\}}_{\n}\big)
    + \det\big(\ltm[\rootp{2}]{\dm[1]{\mxe}}^{\{1\}}_{\n},\ltm[\rootp{2}]{\mxe}^{\{2\}}_{\n}\big) \\
    &= 2(\rootp{2} - \lambda) - \lambda^2 \dltpts[\rootp{2}]{1},
  \end{align*}
  because the first derivative of the matrix \ltm{\mxe} is
  \begin{equation*}
   \ltm{\dm[1]{\mxe}} =  \left(\begin{array}{cc}
                        1           &0 \\
                        \lambda \dltpts{1}       &1
                   \end{array} \right).
  \end{equation*}
  Combining the above we have
  \begin{equation*}
    \delta_2 = \frac{- \rootp{2} \lambda^2 \big(\mean[p]\ltptes[\rootp{2}]-\mean[h]\lthtes[\rootp{2}]\big)}{2(\rootp{2} - \lambda) - \lambda^2 \dltpts[\rootp{2}]{1}}.
  \end{equation*}
  Recall that for the determination of the vector \va[2] we had used the second column of the adjoint matrix, namely we had chosen $m=2$. Thus, according to Remark~\ref{r3.vector k} the vector $\vect{k}_2$ is a zero column vector of dimension 2. Since $\dm[1]{\va[2]}$ is the second column of the matrix \ltm{\dm[1]{\mxe}}, it holds that $\mxb[22]=\delta_2$ and all other elements of \mxb are equal to zero. Finally, $\vd = \big( \frac12 (\mean[p]-\mean[h]),0 \big)$.
\end{toy}

By matching the coefficients of $\epsilon$ on the left and right side of Eq.~\eqref{e3.system of equations mixture model}, we can write the vector of unknown parameters \vupmix as $\vupmix = \vup + \epsilon \vect{z} + O(\epsilon^2 \uv)$. The exact form of the vector \vect{z} is given in the following lemma, which we give without proof.

\begin{lemma}\label{l3.vector u mixture model}
  The vector \vupmix can be written in the form
  \begin{equation*}
    \vupmix = \vup + \epsilon \vect{z} + O(\epsilon^2 \uv),
  \end{equation*}
  where
  \begin{equation*}
    \vect{z} =  \left(\vc\mxa^{-1}\mxb + \vd\right)\mxa^{-1}.
  \end{equation*}
\end{lemma}

\begin{toy}
For the evaluation of \vect{z} we need to find the inverse of matrix \mxa, namely we need
  \begin{equation}
    \mxa^{-1} = \frac{\lambda}{\rootp{2}}\left( \begin{array}{cc}
    \rootp{2}-\lambda       &\lambda \\
    -\frac1{\lambda}       &\frac1{\lambda}
  \end{array} \right).
  \end{equation}
  By observing that $\vc \mxa^{-1} = \vup$ and by following the calculations of Lemma~\ref{l3.vector u mixture model} we obtain
  \begin{align*}
    \vect{z} = \frac{\lambda}{\rootp{2}} \Bigg[ &\frac12(\mean[p]-\mean[h])(\rootp{2}-\lambda)
    - \frac1{\rootp{2}}\Big(1-\frac{\lambda \mean[p]}{2}\Big) \delta_2,\\
    &\frac{\lambda}2(\mean[p]-\mean[h])
    + \frac1{\rootp{2}}\Big(1-\frac{\lambda \mean[p]}{2}\Big) \delta_2  \Bigg].
  \end{align*}
\end{toy}

In our analysis, we used first order perturbation with respect to the parameter $\epsilon$. The exact same procedure can be followed if higher order terms of $\epsilon$ are desired. However, this would result to the increase of the complexity of the formulas. In the next section, we provide the formulas for the evaluation of the perturbed transform vector \vltwmix.

\subsection{Workload distribution of the perturbed model}\label{s3.workload distribution of the perturbed model}
If \ltm{\mxadjmix} is the adjoint matrix of \ltm{\mxemix} (see Eq.~\eqref{e3.definition perturbed matrix e}), then the $i$th element of the transform vector \vltwmix is defined as
\begin{equation}\label{e3.lt mixture workload in each state}
  \ltwcmix{,i} = \frac{s\vupmix \ltm{\mxadjmix}\uv_i}{\det\ltm{\mxemix}}.
\end{equation}
Therefore, to find the exact formula of  \ltwcmix{,i} we need to find $\det\ltm{\mxemix}$ and $s\vupmix \ltm{\mxadjmix}\uv_i$. By using the binomial identity and by omitting higher order powers of $\epsilon$, we have that $\Big(\ltpts + \epsilon s \big(\mean[p]\ltptes-\mean[h]\lthtes\big)\Big)^k = \big(\ltpts\big)^k + \epsilon k \big(\ltpts\big)^{k-1} s\big(\mean[p]\ltptes - \mean[h]\lthtes\big) + O(\epsilon^2)$. We give the following lemmas without proof. The first one gives the formula for the evaluation of the denominator of the desired quantity.
\begin{lemma}\label{l3.determinant matrix E replace}
  If $\det \ltm{\mxe}$ is evaluated according to Theorem~\ref{t3.determinant matrix E} with $\ltgs = \ltpts$, then $\det \ltm{\mxemix}$ can be written as perturbation of $\det \ltm{\mxe}$ as follows
  \begin{align*}
    \det \ltm{\mxemix} = & \det \ltm{\mxe} + \epsilon s \big(\mean[p]\ltptes - \mean[h]\lthtes\big) \\
                         & \times \sum_{k=1}^N k \big(\ltpts\big)^{k-1} \sum_{\substack{\Gamma\subset \n\\  \mid \Gamma\mid =k}}
                                \sum_{\substack{S\subset \n\\ S\supset \Gamma}} \las{S}  \fun{\zeta^{S^c}}\\
                         & \times \det \Big( \big(\mxprobsdummy\circ\mxtrans\big)^{S\setminus \Gamma}_S  \Join \big(\mxprobsreal\circ\mxtrans\big)^\Gamma_S\Big) \\
                        +& O(\epsilon^2).
  \end{align*}
\end{lemma}

\begin{toy}
  Only the combination $k=1$ with  $\Gamma=\{1\}$, and  $S= \n $ gives a non\-/zero coefficient for $\epsilon$. Therefore,
  \begin{align*}
    \det \ltm{\mxemix} =& \det \ltm{\mxe}
                          + \epsilon s \big(\mean[p]\ltptes - \mean[h]\lthtes\big)  \las{\n} \\
                         &\times \fun{\zeta^{\emptyset}} \det \Big( \big(\mxprobsdummy\circ\mxtrans\big)^{\{2\}}_{\n} \Join \big(\mxprobsreal\circ\mxtrans\big)^{\{1\}}_{\n}\Big)\\
                       =&\det \ltm{\mxe} - \epsilon \lambda^2 s \big(\mean[p]\ltptes - \mean[h]\lthtes\big).
  \end{align*}
\end{toy}

The next lemma gives the numerator of each \ltwcmix{,i}, $i \in \n$.

\begin{lemma}\label{l3.numerator LT workload state i replace}
  If $s \vup \ltm{\mxadj}\uv_i$ is evaluated according to Theorem~\ref{t3.numerator LT workload state i} with $\ltgs = \ltpts$, then $s\vupmix \ltm{\mxadjmix}\uv_i$ can be written as perturbation of $s\vup \ltm{\mxadj}\uv_i$ as follows
  \begin{align}
  s\vupmix  \ltm{\mxadjmix} & \uv_i = s\vup \ltm{\mxadj}\uv_i + \epsilon s \Bigg[ z_i \sum_{k=1}^{N-1} \big(\ltpts\big)^k
                                \sum_{\substack{\Gamma\subset \n\setminus\{i\}\\  \mid \Gamma\mid =k;\\ S\subset \n\setminus\{i\}\\ S\supset \Gamma}} \las{S}  \notag \\
                            & \times \fun{\zeta^{S^c}} \det \Big( \big(\mxprobsdummy\circ\mxtrans\big)^{S\setminus \Gamma}_S \Join \big(\mxprobsreal\circ\mxtrans\big)^\Gamma_S  \Big) \notag \\
                           +& z_i \sum_{\substack{S\subset \n\setminus\{i\}}} \las{S} \fun{\zeta^{S^c}} \det \big(\mxprobsdummy\circ\mxtrans\big)^{S}_S    \notag \\
                           +&\sum_{\substack{l=1\\ l\neq i}}^N z_l (-1)^{l+i} \sum_{k=1}^{N-1} \big(\ltpts\big)^k
                                \sum_{\substack{\Gamma\subset \n\setminus\{l,i\}\\  \mid \Gamma\mid =k-1;\\ S\subset \n\setminus\{l,i\}\\ S\supset \Gamma;\\ R\subset S\cap T_{li}}}
                                (-1)^{\mid R\mid}  \las{S\cup \{i\}}  \notag\\
                            & \times  \fun{\zeta^{S^c}}   \det \Big( \big(\mxprobsdummy\circ\mxtrans\big)^{S\setminus \Gamma}_{S\cup \{l\}}
                                \Join \big(\mxprobsreal\circ\mxtrans\big)^{\Gamma \cup \{i\}}_{S\cup \{l\}}  \Big)     \notag \\
                           +&\sum_{\substack{l=1\\ l\neq i}}^N z_l (-1)^{l+i} \sum_{k=0}^{N-2} \big(\ltpts\big)^k \notag \\
                            & \times  \sum_{\substack{\Gamma\subset \n\setminus\{l,i\}\\  \mid \Gamma\mid =k}}
                                \sum_{\substack{S\subset \n\setminus\{l,i\}\\ S\supset \Gamma;\\ R\subset S\cap T_{li}}} (-1)^{\mid R\mid} \las{S\cup \{i\}}  \fun{\zeta^{S^c}} \notag \\
                            & \times \det \Big( \big(\mxprobsdummy\circ\mxtrans\big)^{(S\setminus \Gamma)\cup \{i\}}_{S\cup \{l\}}
                                \Join \big(\mxprobsreal\circ\mxtrans\big)^{\Gamma}_{S\cup \{l\}}  \Big) \notag \\
                           +& s\big(\mean[p]\ltptes - \mean[h]\lthtes\big)\Bigg(u_i \sum_{k=1}^{N-1} k \big(\ltpts\big)^{k-1} \notag \\
                            & \times \sum_{\substack{\Gamma\subset \n\setminus\{i\}\\  \mid \Gamma\mid =k}}
                                \sum_{\substack{S\subset \n\setminus\{i\}\\ S\supset \Gamma}} \las{S} \fun{\zeta^{S^c}} \notag \\
                            & \times \det \Big( \big(\mxprobsdummy\circ\mxtrans\big)^{S\setminus \Gamma}_S \Join \big(\mxprobsreal\circ\mxtrans\big)^\Gamma_S  \Big) \notag\\
                           +& \sum_{\substack{l=1\\ l\neq i}}^N u_l (-1)^{l+i} \sum_{k=1}^{N-1} k \big(\ltpts\big)^{k-1} \notag \\
                            & \times \sum_{\substack{\Gamma\subset \n\setminus\{l,i\}\\  \mid \Gamma\mid =k-1}}
                                \sum_{\substack{S\subset \n\setminus\{l,i\}\\ S\supset \Gamma;\\ R\subset S\cap T_{li}}} (-1)^{\mid R\mid} \las{S\cup \{i\}} \fun{\zeta^{S^c}}\notag\\
                            & \times \det \Big( \big(\mxprobsdummy\circ\mxtrans\big)^{S\setminus \Gamma}_{S\cup \{l\}}
                                \Join \big(\mxprobsreal\circ\mxtrans\big)^{\Gamma \cup \{i\}}_{S\cup \{l\}}  \Big) \notag \\
                           +& \sum_{\substack{l=1\\ l\neq i}}^N u_l (-1)^{l+i} \sum_{k=1}^{N-2} k \big(\ltpts\big)^{k-1} \notag \\
                            & \times  \sum_{\substack{\Gamma\subset \n\setminus\{l,i\}\\  \mid \Gamma\mid =k}}
                                \sum_{\substack{S\subset \n\setminus\{l,i\}\\ S\supset \Gamma;\\ R\subset S\cap T_{li}}} (-1)^{\mid R\mid} \las{S\cup \{i\}} \fun{\zeta^{S^c}} \notag\\
                            & \times \det \Big( \big(\mxprobsdummy\circ\mxtrans\big)^{(S\setminus \Gamma)\cup \{i\}}_{S\cup \{l\}}
                                \Join \big(\mxprobsreal\circ\mxtrans\big)^{\Gamma}_{S\cup \{l\}}  \Big)  \Bigg)  \Bigg]  \notag \\
                           +& O(\epsilon^2), \notag
  \end{align}
  where $z_i$, $i \in \n$, are the coordinates of the vector \vect{z} given in Lemma~\ref{l3.vector u mixture model}.
\end{lemma}

\begin{toy}
  By doing the calculations for each state without taking into account terms that are equal to zero, we obtain:
  \begin{small}
  \begin{align*}
  s\vupmix \ltm{\mxadjmix}\uv_1 =& s\vup \ltm{\mxadj}\uv_1 + \epsilon s\Bigg[  z_1  \las{\emptyset}\fun{\zeta^{\{2\}}} \det \big(\mxprobsdummy\circ\mxtrans\big)^\emptyset_\emptyset\notag \\
                               +& z_2 (-1)^{2+1} \ltpts (-1)^{\mid \emptyset\mid} \las{\emptyset\cup \{1\}} \fun{\zeta^{\emptyset}} q_{21}^{(2)}p_{21} \\
                               +& s\big(\mean[p]\ltptes - \mean[h]\lthtes\big)\Bigg( u_2 (-1)^{2+1}  (-1)^{\mid \emptyset\mid}\\
                                & \times \las{\emptyset\cup \{1\}} \fun{\zeta^{\emptyset}} q_{21}^{(2)} p_{21} \Bigg)
                              \Bigg] + O(\epsilon^2)\\
                               =& s\vup \ltm{\mxadj}\uv_1 + \epsilon s\Big(  z_1  (s-\lambda) - z_2 \lambda \ltpts \\
                                &+ s\big(\mean[p]\ltptes - \mean[h]\lthtes\big) ( - \lambda u_2 )  \Big) + O(\epsilon^2),
  \end{align*}
  \end{small}
  and
  \begin{small}
  \begin{align*}
  s\vupmix \ltm{\mxadjmix}\uv_2 =& s\vup \ltm{\mxadj}\uv_2 + \epsilon s\Bigg[  z_2 \las{\emptyset} \fun{\zeta^{\{1\}}} \det \big(\mxprobsdummy\circ\mxtrans\big)^\emptyset_\emptyset  \notag \\
                              +&z_1 (-1)^{1+2} (-1)^{\mid \emptyset\mid} \las{\emptyset\cup \{2\}} \fun{\zeta^{\emptyset}} q_{12}^{(1)} p_{12} \Bigg] + O(\epsilon^2)\\
                              =&  s\vup \ltm{\mxadj}\uv_2 + \epsilon s\Big(- z_1 \lambda + z_2  (s-\lambda)
                               \Big) + O(\epsilon^2).
  \end{align*}
  \end{small}
\end{toy}

Combining the results of Lemmas~\ref{l3.determinant matrix E replace}--\ref{l3.numerator LT workload state i replace}, we have the following proposition for the transform vector \vltwmix.

\begin{proposition}\label{p3.series expansion of the corrected replace transform vector}
  If \ltwc{i} is calculated according to Proposition~\ref{p3.LT replace workload state i} for the replace base model, then there exist unique coefficients $\beta$, $\gamma$, $\alpha_{i,k}$, $\beta_{i,k}$, $\gamma_{i,k}$, $k=2,\dots,N$, and $\alpha''_{i,j,l}$, $\beta''_{i,j,l}$ and $\gamma''_{i,j,l}$, $j=1,\dots,\sigma$, $l=1,\dots,r_{i,j}$, such that
  \begin{align*}
    \ltwcmix{,i} =& \ltwc{i} + \epsilon \frac1{u_i}\ltwc{i} \Bigg[ \Bigg(z_i + \sum_{k=2}^N \frac{\alpha_{i,k}}{s-\rootp{k}} \\
                  &+ \sum_{j=1}^{\sigma}\sum_{l=1}^{r_{i,j}} \frac{\alpha''_{i,j,l} \cdot (\rootnnum{i,j})^{r_{i,j}-l+1}}{(s+\rootnnum{i,j})^{r_{i,j}-l+1}}\Bigg) \\
                 +& \big(\mean[p]\ltptes - \mean[h]\lthtes\big)\Bigg( \beta
                   + \sum_{k=2}^N \frac{\beta_{i,k}}{s-\rootp{k}} \\
                  &+ \sum_{j=1}^{\sigma}\sum_{l=1}^{r_{i,j}} \frac{\beta''_{i,j,l} \cdot (\rootnnum{i,j})^{r_{i,j}-l+1}}{(s+\rootnnum{i,j})^{r_{i,j}-l+1}}\Bigg) \\
                 -& \big(\mean[p]\ltptes - \mean[h]\lthtes\big)\ltwc{i} \Bigg( \gamma
                   + \sum_{k=2}^N \frac{\gamma_{i,k}}{s-\rootp{k}} \\
                  & + \sum_{j=1}^{\sigma}\sum_{l=1}^{r_{i,j}} \frac{\gamma''_{i,j,l} \cdot (\rootnnum{i,j})^{r_{i,j}-l+1}}{(s+\rootnnum{i,j})^{r_{i,j}-l+1}} \Bigg) \Bigg]
                  + O(\epsilon^2),
  \end{align*}
where $z_i$, $i \in \n$, are the coordinates of the vector \vect{z} given in Lemma~\ref{l3.vector u mixture model}.
\end{proposition}
\begin{proof}
  See Appendix~\ref{appendix B}.
\end{proof}

Before we evaluate the Laplace transform of the waiting time of a real customer \ltdmix in our running example, we apply Laplace inversion to the coefficient of $\epsilon$ in the series expansion of \ltwcmix{,i}. We denote by \erlang[k]{\lambda} the r.v.\ that follows an Erlang distribution with $k$ phases and rate $\lambda$. For simplicity, we write \erlang{\lambda} for the exponential r.v.\ with rate $\lambda$. Finally, let \epts and \ehts be the generic stationary excess phase\-/type and heavy\-/tailed service times, respectively.

\begin{theorem}\label{t3.series expansion of mixture workload state i}
If \fun{\lt{\theta}_i} is the coefficient of $\epsilon$ in the series expansion of \ltwcmix{,i} in Proposition~\ref{p3.series expansion of the corrected replace transform vector}, its Laplace inversion $\fun[t]{\Theta_i}=\mathcal{L}^{-1}\{\fun{\lt{\theta}_i}\}$ is given as follows
\begin{small}
\begin{align*}
  \fun[t]{\Theta_i} &=\frac1{u_i}\Bigg[\Big(z_i - \sum_{k=2}^N \frac{\alpha_{i,k}}{\rootp{k}} \Big) \pr (\workload[i] >t ) \\
                      +& \Big(\beta - \sum_{k=2}^N \frac{\beta_{i,k}}{\rootp{k}} \Big)
                              \Big( \mean[p] \pr (\workload[i] + \epts >t ) -\mean[h] \pr (\workload[i] + \ehts >t ) \Big) \notag \\
                      -& \Big(\gamma - \sum_{k=2}^N \frac{\gamma_{i,k}}{\rootp{k}} \Big) \Big( \mean[p] \pr (\workload[i] + \workload[i]' + \epts >t ) \\
                       &- \mean[h] \pr (\workload[i] + \workload[i]' + \ehts >t ) \Big) \notag \\
  -& \sum_{k=2}^N \frac1{\rootp{k}} \Bigg( \gamma_{i,k} \Big( \mean[p] \pr \big(t < \workload[i] + \workload[i]' + \epts < t + \erlang{\rootp{k}}\big) \\
                       &- \mean[h] \pr \big(t<\workload[i] + \workload[i]' + \ehts < t + \erlang{\rootp{k}} \big) \Big) \notag \\
                       &- \beta_{i,k} \Big( \mean[p] \pr\big(t<\workload[i] + \epts < t + \erlang{\rootp{k}} \big) \\
                       &-\mean[h] \pr\big(t<\workload[i] + \ehts < t + \erlang{\rootp{k}}\big) \Big) \\
                       &- \alpha_{i,k} \pr\big(t<\workload[i] < t + \erlang{\rootp{k}} \big) \Bigg) \notag \\
  -& \sum_{j=1}^{\sigma}\sum_{l=1}^{r_{i,j}} \Bigg( \gamma''_{i,j,l} \Big( \mean[p] \pr\big(\workload[i]+ \workload[i]'+ \epts + \erlang[r_{i,j}-l+1]{\rootnnum{i,j}}>t\big) \\
                       &- \mean[h] \pr\big(\workload[i]+ \workload[i]'+ \ehts + \erlang[r_{i,j}-l+1]{\rootnnum{i,j}}>t \big) \Big) \notag \\
                       &- \beta''_{i,j,l} \Big(  \mean[p] \pr\big(\workload[i]+ \epts + \erlang[r_{i,j}-l+1]{\rootnnum{i,j}}>t\big) \\
                       &- \mean[h] \pr\big(\workload[i]+ \ehts + \erlang[r_{i,j}-l+1]{\rootnnum{i,j}}>t\big) \Big) \notag \\
                             &-\alpha''_{i,j,l} \pr\big(\workload[i]+ \erlang[r_{i,j}-l+1]{\rootnnum{i,j}}>t\big)  \Bigg) \Bigg],
\end{align*}
\end{small}
where $\workload[i]'$ is independent and follows the same distribution of \workload[i].
\end{theorem}
\begin{proof}
  See Appendix~\ref{appendix B}.
\end{proof}

\begin{remark}\label{r3.real roots}
  Note that an \erlang[k]{\lambda} distribution ($k\geq 1$) is defined for a non\-/negative real valued rate $\lambda$. To state Theorem~\ref{t3.series expansion of mixture workload state i}, we assumed that all the roots \rootp{k}, $k=2,\dots,N$, and $-\rootnnum{i,j}$, $j=1,\dots,rM$, are real-valued for all $i \in \n$. In most systems, this assumption in not always true. Recall that the previously mentioned roots are roots of a polynomial with real coefficients (see analysis above Eq.~\eqref{e3.numerator as a product replace}). Therefore, from the Complex Conjugate Root Theorem
  it holds that if e.g.\ \rootp{2} is complex, then its complex conjugate $\overline{\rootp{2}}$ is also a root. Thus, we write $E_{Re(\rootp{2})}$ instead of $E_{\rootp{2}}$ and $E_{\overline{\rootp{2}}}$, because every parameter or function that depends on $\overline{\rootp{2}}$ appears as a complex conjugate of the corresponding quantity that depends on $\rootp{2}$, and their imaginary parts cancel out. The same result holds for all other roots.
\end{remark}

\begin{toy}
  For the evaluation of the Laplace transform $\ltdmix = \vltwmix \vweights$ of the waiting time of a real customer \delaymix, we follow similar steps as in the proof of Proposition~\ref{p3.series expansion of the corrected replace transform vector}. Recall that in our example, $r=1$, and assume that only $\sigma$ of the roots $-\rootnnum{j}$ are distinct and that the multiplicity of each of them is $r_{j}$, such that $\sum_{j=1}^\sigma r_{j} = M$.

   Therefore, we first find $\fun{p} \det \ltm{\mxemix}$ and $\fun{p} s \vupmix \ltm{\mxadjmix}\vweights$. If we set $\fun{\xi}= -\lambda^2 \fun{p}$, $\fun{\xi'_1}= -2\lambda \fun{p}$, and $\fun{\xi'_2} = 2(s-\lambda) \fun{p}$,
%  \begin{align*}
%    \fun{\xi}     &= -\lambda^2 \fun{p}, \\ %\label{e3.function xi toy}
%    \fun{\xi'_1}  &= -2\lambda \fun{p}, \\ %\label{e3.function xi prime 1 toy}
%    \fun{\xi'_2}  &= 2(s-\lambda) \fun{p}, %\label{e3.function xi prime 2 toy}
%  \end{align*}
  then we obtain
  \begin{align*}
    \fun{p} &\det \ltm{\mxemix}             =  \fun{p}\det \ltm{\mxe} \\
                                           &+ \epsilon s \big(\mean[p]\ltptes - \mean[h]\lthtes\big) \fun{\xi} + O(\epsilon^2),\notag \\
    \fun{p} & s \vupmix \ltm{\mxadjmix}\vweights  = \fun{p} s \vup \ltm{\mxadj} \vweights \\
                                           &+ \epsilon s \sum_{l=1}^2 z_l \fun{\xi'_{l}}
                                                  + O(\epsilon^2). \notag
  \end{align*}
  In our case, we define the functions \fes and \fef (see Eqs.~\eqref{e3.function d mixture model} and \eqref{e3.function n mixture model} respectively) as
  \begin{align*}
  \fes =& \frac{\big(\mean[p]\ltptes - \mean[h]\lthtes\big) \fun{\xi}\ltd}{\vup \vweights (s-\rootp{2}) \prod_{j=1}^{\sigma} (s+\rootnnum{j})^{r_j}}
              - \frac{\delta_2}{s - \rootp{2}}, \\ % \label{e3.function d toy}
  \fef =& \frac{\sum_{l=1}^2 z_l \fun{\xi'_{l}}}{\vup \vweights (s-\rootp{2}) \prod_{j=1}^{\sigma} (s+\rootnnum{j})^{r_j}}
         - \frac{\delta_2}{s - \rootp{2}}, %\label{e3.function n toy}
  \end{align*}
  where the two equivalent definitions of $\delta_2$ (see Eqs.~\eqref{e3.definition 1 for dk} and \eqref{e3.definition 2 for dk}) take the form
  \begin{align*}%\label{e3.definition 1 for dk}
  \delta_2  =& \frac{\big(\mean[p]\ltptes[\rootp{2}] - \mean[h]\lthtes[\rootp{2}]\big) \fun[\rootp{2}]{\xi}\ltd[\rootp{2}]}{\vup \vweights \prod_{j=1}^{\sigma} (\rootp{2}+\rootnnum{j})^{r_j}} \\
            =& \frac{\sum_{l=1}^2 z_l \fun[\rootp{2}]{\xi'_{l}}}{\vup \vweights \prod_{j=1}^{\sigma} (\rootp{2}+\rootnnum{j})^{r_j}}.
  \end{align*}
  Following the calculations after Eq.~\eqref{e3.numerator mixture} we get that
  \begin{align}
  \ltdmix  &= \ltd + \epsilon \frac1{\vup \vweights }\ltd \Bigg( \frac{\sum_{l=1}^2 z_l \fun{\xi'_{l}}}{(s-\rootp{2}) \prod_{j=1}^{\sigma} (s+\rootnnum{j})^{r_j}} \notag \\
          -& \big(\mean[p]\ltptes - \mean[h]\lthtes\big)\ltd \frac{\fun{\xi}}{(s-\rootp{2}) \prod_{j=1}^{\sigma} (s+\rootnnum{j})^{r_j}} \Bigg) \notag \\
          +& O(\epsilon^2). \label{e3.laplace transform mixture workload toy}
  \end{align}
  Now, we apply simple fraction decomposition to the rational functions
  \begin{gather*}
    \frac{\sum_{l=1}^2 z_l \fun{\xi'_{l}}}{(s-\rootp{2}) \prod_{j=1}^{\sigma} (s+\rootnnum{j})^{r_j}}, \qquad
    \frac{\fun{\xi}}{(s-\rootp{2}) \prod_{j=1}^{\sigma} (s+\rootnnum{j})^{r_j}}.
  \end{gather*}
  Thus, we calculate
  \begin{gather*}
  \alpha_{2}  = \frac{\sum_{l=1}^2 z_l \fun[\rootp{2}]{\xi'_{l}}}{\prod_{j=1}^{\sigma} (\rootp{2}+\rootnnum{j})^{r_j}}, \qquad
  \gamma_{2}  = \frac{\fun[\rootp{2}]{\xi}}{\prod_{j=1}^{\sigma} (\rootp{2}+\rootnnum{j})^{r_j}}, %\label{e3.coeff toy}
  \end{gather*}
  and for $j=1,\dots,\sigma$, $p=1,\dots,r_{j}$, the coefficients $\alpha''_{j,p}$ and $\gamma''_{j,p}$, are respectively the unique solutions to the following two linear systems of $r_{j}$ equations
  \begin{small}
  \begin{gather*}
    \frac{d}{ds^n}\left.\Bigg[\sum_{l=1}^2 z_l \fun{\xi'_{l}}\Bigg]\right|_{s=-\rootnnum{j}}
                = \\
                \frac{d}{ds^n}\left.\Bigg[ (s-\rootp{2}) \prod_{\substack{l=1\\ l\neq j}}^{\sigma}(s+\rootnnum{l})^{r_{l}}
                        \sum_{p=1}^{r_{j}} \alpha''_{j,p}(\rootnnum{j})^{r_{j}-p+1}(s+\rootnnum{j})^{p-1} \Bigg]\right|_{s=-\rootnnum{j}},\\%\quad n=0,\dots,r_{j},
    \frac{d}{ds^n}\left.\Bigg[\fun{\xi}\Bigg]\right|_{s=-\rootnnum{j}}
                = \\
                \frac{d}{ds^n}\left.\Bigg[ (s-\rootp{2}) \prod_{\substack{l=1\\ l\neq j}}^{\sigma}(s+\rootnnum{l})^{r_{l}}
                        \sum_{p=1}^{r_{j}} \gamma''_{j,p} (\rootnnum{j})^{r_{j}-p+1}(s+\rootnnum{j})^{p-1} \Bigg]\right|_{s=-\rootnnum{j}}, %\quad n=0,\dots,r_{j}.
  \end{gather*}
  \end{small}
  $n=0,\dots,r_{j}$. In addition, the polynomial $\fun{\xi}$ is of degree $M$, and the polynomial $\sum_{l=1}^2 z_l \fun{\xi'_{l}}$ is of degree $M+1$ with the coefficient of $s^{M+1}$ equal to $2 z_2$. Combining all these, we write Eq.~\eqref{e3.laplace transform mixture workload toy} as
  \begin{align*}
  \ltdmix =& \ltd + \epsilon \frac1{2 u_2}\ltd \Bigg[ \Bigg(2 z_2+  \frac{\alpha_{2}}{s-\rootp{2}} \\
           &+ \sum_{j=1}^{\sigma}\sum_{l=1}^{r_{j}} \frac{\alpha''_{j,l} \cdot (\rootnnum{j})^{r_{j}-l+1}}{(s+\rootnnum{j})^{r_{j}-l+1}}\Bigg) \\
          -& \big(\mean[p]\ltptes - \mean[h]\lthtes\big)\ltd \Bigg( \frac{\gamma_{2}}{s-\rootp{2}} \\
           & + \sum_{j=1}^{\sigma}\sum_{l=1}^{r_{j}} \frac{\gamma''_{j,l} \cdot (\rootnnum{j})^{r_{j}-l+1}}{(s+\rootnnum{j})^{r_{j}-l+1}} \Bigg) \Bigg]
                + O(\epsilon^2).
  \end{align*}
  Comparing the above formula with the one in Proposition~\ref{p3.series expansion of the corrected replace transform vector}, we see that here we used \ltdmix, $2u_2$, and $2z_2$ instead of \ltwcmix{,i}, $u_i$, and $z_i$, respectively. Moreover, we found that in our case $\gamma=0$ and all $\beta$ coefficients are also equal to zero. Thus, if \fun{\lt{\theta}} is the coefficient of $\epsilon$ in the series expansion of \ltdmix, we apply Theorem~\ref{t3.series expansion of mixture workload state i} to find its Laplace inversion as
  \begin{small}
  \begin{align*}
    \fun[t]{\Theta}&  = \frac1{2 u_2}
                        \Bigg[\Big(2 z_2  - \frac{\alpha_{2}}{\rootp{2}} \Big) \pr (\delay >t ) \\
                  +&  \frac{\gamma_{2}}{\rootp{2}}  \Big( \mean[p] \pr (\delay + \delay{'} + \epts >t ) %\\
                    - \mean[h] \pr (\delay + \delay{'} + \ehts >t ) \Big) \notag \\
  -& \frac1{\rootp{2}} \Bigg( \gamma_{2} \Big( \mean[p] \pr \big(t < \delay + \delay{'} + \epts < t + \erlang{\rootp{2}}\big) \\
                   &- \mean[h] \pr \big(t<\delay + \delay{'} + \ehts < t + \erlang{\rootp{2}} \big) \Big) \notag \\
                   &- \alpha_{2} \pr\big(t<\delay < t + \erlang{\rootp{2}} \big) \Bigg) \notag \\
  -& \sum_{j=1}^{\sigma}\sum_{l=1}^{r_{j}} \Bigg( \gamma''_{j,l} \Big( \mean[p] \pr\big(\delay+ \delay{'}+ \epts + \erlang[r_{j}-l+1]{\rootnnum{j}}>t\big) \\
                   &- \mean[h] \pr\big(\delay+ \delay{'}+ \ehts + \erlang[r_{j}-l+1]{\rootnnum{j}}>t \big) \Big) \notag \\
                   &-\alpha''_{j,l} \pr\big(\delay+ \erlang[r_{j}-l+1]{\rootnnum{j}}>t\big)  \Bigg) \Bigg],
  \end{align*}
  \end{small}
  where $\delay{'}$ is independent and follows the same distribution of \delay.
\end{toy}

By applying Laplace inversion to the first two terms of the series expansion in $\epsilon$ of the workload just prior to an arrival of a customer in each state, we obtain that the first term is a phase\-/type approximation of the aforementioned workload that results from the replace base model (see Section~\ref{s3.replace base model}). In addition, the second term, which we refer to as correction term and is found explicitly in Theorem~\ref{t3.series expansion of mixture workload state i}, involves linear combinations of terms that have probabilistic interpretation. More precisely, these terms with probabilistic interpretation are either tail probabilities of convoluted r.v.\ or probabilities for some of the aforementioned convoluted r.v.\ to lie between a fixed value $t$ and the same value $t$ shifted by an exponential time. Finally, observe that these convoluted r.v.\ involve the heavy\-/tailed stationary\-/excess service time r.v.\ \ehts in a maximum appearance of one. Combining the results of Proposition~\ref{p3.series expansion of the corrected replace transform vector} and Theorem~\ref{t3.series expansion of mixture workload state i}, in the next section we define our approximations.

\subsection{Corrected replace approximation}\label{s3.properties corrected replace approximation}
The goal of this section is to provide approximations that maintain the numerical tractability but improve the accuracy of the phase\-/type approximations and that are able to capture the tail behavior of the exact workload distribution. As we pointed out in the introduction, a single appearance of a stationary excess heavy\-/tailed service time \ehts is sufficient to capture the correct tail behavior of the exact workload. As we observed in Section~\ref{s3.workload distribution of the perturbed model}, the correction term contains terms with single appearances of \ehts. For this reason, the proposed approximation for the workload is constructed by the first two terms of its respective series expansion. We propose the following approximation:

\begin{approximation}\label{d3.corrected replace approximation}
  The corrected replace approximation of the survival function $\pr(\workloadmix[,i]>t)$ of the exact workload prior to an arrival of a customer in state $i$, $i \in \n$, is defined as
  \begin{small}
  \begin{align*}
  \correpc{i} &:= \pr(\workload[i] > t)  + \epsilon\frac1{u_i}
                            \Bigg[\Big(z_i - \sum_{k=2}^N \frac{\alpha_{i,k}}{\rootp{k}} \Big) \pr (\workload[i] >t ) \\
                      +& \Big(\beta - \sum_{k=2}^N \frac{\beta_{i,k}}{\rootp{k}} \Big)
                              \Big( \mean[p] \pr (\workload[i] + \epts >t ) -\mean[h] \pr (\workload[i] + \ehts >t ) \Big) \notag \\
                      -& \Big(\gamma - \sum_{k=2}^N \frac{\gamma_{i,k}}{\rootp{k}} \Big) \Big( \mean[p] \pr (\workload[i] + \workload[i]' + \epts >t ) \\
                       &- \mean[h] \pr (\workload[i] + \workload[i]' + \ehts >t ) \Big) \notag \\
  -& \sum_{k=2}^N \frac1{\rootp{k}} \Bigg( \gamma_{i,k} \Big( \mean[p] \pr \big(t < \workload[i] + \workload[i]' + \epts < t + \erlang{\rootp{k}}\big) \\
                       &- \mean[h] \pr \big(t<\workload[i] + \workload[i]' + \ehts < t + \erlang{\rootp{k}} \big) \Big) \notag \\
                       &- \beta_{i,k} \Big( \mean[p] \pr\big(t<\workload[i] + \epts < t + \erlang{\rootp{k}} \big) \\
                       &-\mean[h] \pr\big(t<\workload[i] + \ehts < t + \erlang{\rootp{k}}\big) \Big) \\
                       &- \alpha_{i,k} \pr\big(t<\workload[i] < t + \erlang{\rootp{k}} \big) \Bigg) \notag \\
  -& \sum_{j=1}^{\sigma}\sum_{l=1}^{r_{i,j}} \Bigg( \gamma''_{i,j,l} \Big( \mean[p] \pr\big(\workload[i]+ \workload[i]'+ \epts + \erlang[r_{i,j}-l+1]{\rootnnum{i,j}}>t\big) \\
                       &- \mean[h] \pr\big(\workload[i]+ \workload[i]'+ \ehts + \erlang[r_{i,j}-l+1]{\rootnnum{i,j}}>t \big) \Big) \notag \\
                       &- \beta''_{i,j,l} \Big(  \mean[p] \pr\big(\workload[i]+ \epts + \erlang[r_{i,j}-l+1]{\rootnnum{i,j}}>t\big) \\
                       &- \mean[h] \pr\big(\workload[i]+ \ehts + \erlang[r_{i,j}-l+1]{\rootnnum{i,j}}>t\big) \Big) \notag \\
                             &-\alpha''_{i,j,l} \pr\big(\workload[i]+ \erlang[r_{i,j}-l+1]{\rootnnum{i,j}}>t\big)  \Bigg) \Bigg],
  \end{align*}
  \end{small}
  where $\pr(\workload[i] > t)$ is the replace phase\-/type approximation of $\pr(\workloadmix[,i]>t)$, $\workload[i]'$ is independent and follows the same distribution of \workload[i], and the coefficients $\beta$, $\gamma$, $\alpha_{i,k}$, $\beta_{i,k}$, $\gamma_{i,k}$, $k=2,\dots,N$, and $\alpha''_{i,j,l}$, $\beta''_{i,j,l}$ and $\gamma''_{i,j,l}$, $j=1,\dots,\sigma$, $l=1,\dots,r_{i,j}$, are calculated according to Proposition~\ref{p3.series expansion of the corrected replace transform vector}.
\end{approximation}
The following result shows that the corrected replace approximation makes sense rigorously.
\begin{proposition}\label{p3.extended continuity theorem}
  If $\pr(\workload[i] > t)$ is the replace approximation of the exact workload just prior to an arrival of a customer in state $i$, $i \in \n$, $\pr(\workloadmix[,i] > t)$, then as $\epsilon \rightarrow 0$, it holds that
  \begin{equation*}
    \frac{\pr(\workloadmix[,i] > t) - \pr(\workload[i] > t)}{\epsilon} \rightarrow \fun[t]{\Theta_i},
  \end{equation*}
  where \fun[t]{\Theta_i} is given in Theorem~\ref{t3.series expansion of mixture workload state i}.
\end{proposition}
\begin{proof}
  See Appendix~\ref{appendix B}.
\end{proof}

Although Approximation~\ref{d3.corrected replace approximation} gives an approximation of the workload that can be calculated explicitly and is computationally tractable, it involves the evaluation of many terms. Therefore, to simplify the formula of the approximation, it makes sense to ignore terms that do not contribute significantly to the accuracy of the corrected replace approximation. Such terms seem to be the probabilities of convoluted r.v.\ that lie between a fixed value $t$ and the same value $t$ shifted by an exponential time. Therefore, we define the {\it simplified corrected replace approximation} as follows.

\begin{approximation}\label{d3.simplified corrected replace approximation}
  The simplified corrected replace approximation of the survival function $\pr(\workloadmix[,i]>t)$ of the exact workload just prior to an arrival of a customer in state $i$, $i \in \n$, is defined as
  \begin{small}
  \begin{align*}
  \corsimrepc{i} &:= \pr(\workload[i] > t)  + \epsilon\frac1{u_i}
                            \Bigg[\Big(z_i - \sum_{k=2}^N \frac{\alpha_{i,k}}{\rootp{k}} \Big) \pr (\workload[i] >t ) \\
                      +& \Big(\beta - \sum_{k=2}^N \frac{\beta_{i,k}}{\rootp{k}} \Big)
                              \Big( \mean[p] \pr (\workload[i] + \epts >t ) -\mean[h] \pr (\workload[i] + \ehts >t ) \Big) \notag \\
                      -& \Big(\gamma - \sum_{k=2}^N \frac{\gamma_{i,k}}{\rootp{k}} \Big) \Big( \mean[p] \pr (\workload[i] + \workload[i]' + \epts >t ) \\
                       &- \mean[h] \pr (\workload[i] + \workload[i]' + \ehts >t ) \Big) \notag \\
  -& \sum_{j=1}^{\sigma}\sum_{l=1}^{r_{i,j}} \Bigg( \gamma''_{i,j,l} \Big( \mean[p] \pr\big(\workload[i]+ \workload[i]'+ \epts + \erlang[r_{i,j}-l+1]{\rootnnum{i,j}}>t\big) \\
                       &- \mean[h] \pr\big(\workload[i]+ \workload[i]'+ \ehts + \erlang[r_{i,j}-l+1]{\rootnnum{i,j}}>t \big) \Big) \notag \\
                       &- \beta''_{i,j,l} \Big(  \mean[p] \pr\big(\workload[i]+ \epts + \erlang[r_{i,j}-l+1]{\rootnnum{i,j}}>t\big) \\
                       &- \mean[h] \pr\big(\workload[i]+ \ehts + \erlang[r_{i,j}-l+1]{\rootnnum{i,j}}>t\big) \Big) \notag \\
                             &-\alpha''_{i,j,l} \pr\big(\workload[i]+ \erlang[r_{i,j}-l+1]{\rootnnum{i,j}}>t\big)  \Bigg) \Bigg],
  \end{align*}
  \end{small}
  where $\pr(\workload[i] > t)$ is the replace phase\-/type approximation of $\pr(\workloadmix[,i]>t)$, $\workload[i]'$ is independent and follows the same distribution of \workload[i], and the coefficients $\beta$, $\gamma$, $\alpha_{i,k}$, $\beta_{i,k}$, $\gamma_{i,k}$, $k=2,\dots,N$, and $\alpha''_{i,j,l}$, $\beta''_{i,j,l}$ and $\gamma''_{i,j,l}$, $j=1,\dots,\sigma$, $l=1,\dots,r_{i,j}$, are calculated according to Proposition~\ref{p3.series expansion of the corrected replace transform vector}.
\end{approximation}

\begin{remark}\label{r3.definition corrected replace approximation for total workload}
  One way to define the corrected replace approximations for the survival function $\pr(\delaymix > t)$ of the waiting time of a real customer is to follow the steps in the running example. An alternative way is to define the approximations of $\pr(\delaymix > t)$ as the weighted sum of the approximations of the survival functions $\pr(\workloadmix[,i]~>~t)$, $i \in \n$. More precisely, the corrected replace approximation of the waiting time is defined as $\correp = \sum_{i \in \n}\weight{i}\allowbreak \times\correpc{i}$, and the simplified corrected replace approximation is defined as $\corsimrep = \sum_{i \in \n}\weight{i}\corsimrepc{i}$. Both approaches lead to the same result.
\end{remark}

In the next section, we perform numerical experiments to check the accuracy of the corrected replace and the simplified corrected replace approximations. In addition, we show that indeed the corrected replace approximation does not differ significantly from its simplified version.

\section{Numerical experiment}\label{s3.numerics}
In Section~\ref{s3.workload distribution of the perturbed model}, we pointed out that the first term of the corrected replace expansion is already a phase\-/type approximation of the workload. In this section we show that adding the correction term leads to improved approximations that are significantly more accurate than their phase\-/type counterpart. Therefore, we check here the accuracy of the corrected replace approximations (see Definitions~\ref{d3.corrected replace approximation} and \ref{d3.simplified corrected replace approximation}) by comparing them with the exact workload distribution and their corresponding phase\-/type approximation.

For the MAP arrival process of customers we choose a MMPP with two states. Since it is more meaningful to compare approximations with exact results than with simulation outcomes, we choose the service time distribution such that we can find an exact formula for the workload.

As service time distribution we use a mixture of an exponential distribution with rate $\nu$ and a heavy\-/tailed one that belongs to a class of long-tailed distributions introduced in \cite{abate99a}. The Laplace transform of the latter distribution is $\lthts = 1-\frac{s}{(\kappa+\sqrt{s})(1+\sqrt{s})}$, where $\e C=\kappa^{-1}$ and all higher moments are infinite. Furthermore, the Laplace transform of the stationary heavy\-/tailed claim size distribution is
\begin{equation*}
  \lthtes = \frac{\kappa}{(\kappa+\sqrt{s})(1+\sqrt{s})},
\end{equation*}
which for $\kappa\neq 1$ can take the form
\begin{equation*}
  \lthtes = \left(\frac{\kappa}{1-\kappa}\right) \left(\frac{1}{\kappa+\sqrt{s}} - \frac{1}{1+\sqrt{s}}\right).
\end{equation*}
For this combination of service time distributions, the survival workload can be found explicitly, by following same ideas as in Theorem~9 of \cite{vatamidou13}.

What is left now is to fix values for the parameters of the mixture model and perform our numerical experiments. Thus, for the MMPP arrival process we choose the parameters such that $\lambda_1=7$, $\lambda_2=1/2$, $p_{11} = 8/9$, and $p_{22}=3/{100}$ (the rest of the parameters can be calculated using the formulas \eqref{e3.definition rates}--\eqref{e3.definition conditional probabilities of real customer}). Although we do not have any restrictions for the parameters of the involved service time distributions, from a modeling point of view, it is counterintuitive to fit a heavy\-/tailed claim size distribution with a mean smaller than the mean of the phase\-/type claim size distribution. For this reason, we select $\kappa=2$ and $\nu=3$.

Finally, note that we performed extensive numerical experiments for various values of the perturbation parameter $\epsilon$ in the interval $[0.001,0.1]$. We chose to present only one example for $\epsilon=0.01$, since the qualitative conclusions for all other values of $\epsilon$ are similar to those presented in this section. The load of the system for this choice of parameters is then equal to 0.852.

\begin{figure}
 \centering
 \includegraphics[scale=0.85]{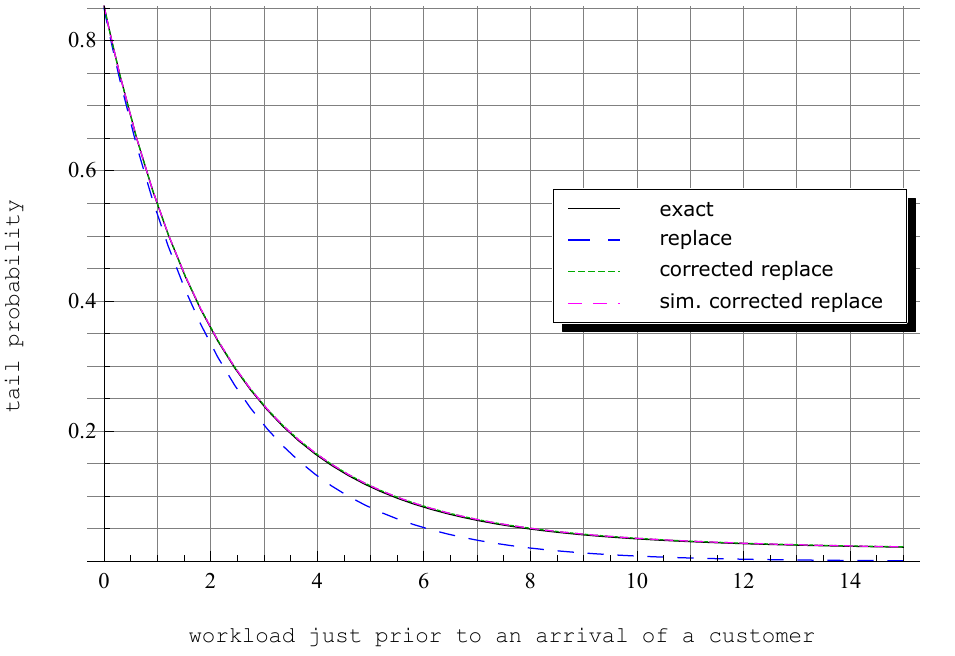}
\caption{Exact workload with phase-type, corrected replace and simplified corrected replace approximations for perturbation parameter 0.01 and load of the system 0.8527}\label{f3.approximations total workload}
\end{figure}

As we observe from Figure~\ref{f3.approximations total workload}, the replace phase\-/type approximation gives accurate estimates for small values of the workload, while it is incapable of capturing the correct tail behavior of the exact survival function of the workload. Contrary, both corrected replace approximations are highly accurate and give a small relative error at the tail. More precisely, we can observe the following:
\begin{itemize}
  \item The corrected replace approximation does not differ significantly from its simplified version. The maximum observed absolute error between the two approximations is approximately equal to $0.00073$.
  \item We found that the absolute error between the exact workload and the corrected replace approximation lies in the interval $[0.00045,0.00047]$, and the absolute error between the exact and the simplified corrected approximation lies in the interval $[0.0008,0.0009]$.
  \item Finally, we found that the relative error at the tail for both corrected replace approximations is smaller than $0.04$.
\end{itemize}

\section{Conclusions}\label{s3.conclusions}
To conclude, both corrected replace approximations are highly accurate and there is no significant difference between them. For this reason, it makes sense to use only the simplified corrected replace approximation to obtain reliable estimates for the workload. In addition, the corrected replace approximations give a small relative error at the tail. More precisely, the relative error at the tail is $O(\epsilon)$ and we are currently writing a rigorous proof for this statement.

\section{Acknowledgments}
The work of Maria Vlasiou and Eleni Vatamidou is supported by Netherlands Organisation for Scientific Research (NWO) through project number 613.001.006. The work of Bert Zwart is supported by an NWO VIDI grant and an IBM faculty award. Bert Zwart is also affiliated with VU University Amsterdam, and the Georgia Institute of Technology.

%
%\bibliographystyle{apt}
%\bibliography{D:/elena}

\appendix
\section{Results on perturbation theory}\label{appendix A}
In this section, we provide some preliminary results on linear algebra, matrix functions, and perturbation theory that are needed in our analysis. We introduce an $N\times N$ matrix function \ltm{\mxe} with a single parameter $s>0$. We say that the matrix function \ltm{\mxe} is {\it regular} if $\det \ltm{\mxe}$ is not identically zero as a function of $s$. In addition, if \ltm{\mxe} is regular (we denote it as  $\det \ltm{\mxe} \not\equiv0$), then the eigenvalues of \ltm{\mxe} are the solutions of the equations $\det \ltm{\mxe} =0$ \cite{deteran11}. Throughout our analysis, we assume that the matrix \ltm{\mxe} is regular and that $r$ is a simple eigenvalues of it. In addition, we assume that the matrix \ltm{\mxe} is analytic in the neighborhood of $r$.  We use the notation \ltm{\dm{\mxe}} for the $n$th derivative of the matrix function \ltm{\mxe}. Thus, \ltm{\mxe} can be written as a Taylor series in the following form:
\begin{equation}\label{e3.taylor series matrix E}
  \ltm{\mxe} = \ltm[r]{\dm[0]{\mxe}} + (s-r) \ltm[r]{\dm[1]{\mxe}} + \dots = \sum_{n=0}^\infty \frac{(s-r)^n}{n!} \ltm[r]{\dm{\mxe}}.
\end{equation}
To avoid redundant notation, in the forthcoming analysis we use the conventions that $\mxe = \ltm[r]{\dm[0]{\mxe}}=\ltm[r]{\mxe}$ and $\dm{\mxe} = \ltm[r]{\dm{\mxe}}$.

 As a consequence of the fact that the multiplicity of the eigenvalue $r$ is one, the dimension of the nullspace of \mxe is equal to one. Our first goal is to find the form of the eigenvectors of the nullspace of matrix \mxe. The following theorem gives us exactly the form of these eigenvectors.

\begin{theorem}\label{t3.eigenvectors that correspond to the eigenvalue zero}
  If \mxex is an $N\times N$ matrix with determinant equal to zero, i.e.\ $\det{\mxex}=0$, and nullspace of dimension one, then a right $N\times 1$ eigenvector that corresponds to the simple eigenvalue zero is \vect{t} with coordinates $t_j =(-1)^{1+j} \det \mxex_{\n\setminus\{1\}}^{\n\setminus\{j\}}$, $j \in \n$.
\end{theorem}
\begin{proof}
  We need to prove that the inner product of every row of \mxex with \vect{t} is equal to zero. More precisely, if $\vect{c}_i$ denotes the $i$th row of matrix \mxex, we need to show that
  \begin{equation*}
    \vect{c}_i\vect{t} =0, \qquad i \in \n.
  \end{equation*}
  If $c_{ij}$ is the $(i,j)$ element of matrix \mxex, for the first row we have
  \begin{equation*}
    \vect{c}_1\vect{t} = \sum_{j=1}^N c_{1j} (-1)^{1+j} \det \mxex_{\n\setminus\{1\}}^{\n\setminus\{j\}}
                        \stackrel{\text{def.}}{=} \det \mxex = 0.
  \end{equation*}
  For an arbitrary row $i=2,\dots,N$, we have
  \begin{align*}
    \vect{c}_i\vect{t} &= \sum_{j=1}^N c_{ij} (-1)^{1+j} \det \mxex_{\n\setminus\{1\}}^{\n\setminus\{j\}}.
  \end{align*}
  We expand the determinant of each matrix $\mxex_{\n\setminus\{1\}}^{\n\setminus\{j\}}$, $j \in \n$, in minors of the $i$th row of matrix \mxex. Observe that the $i$th row of the initial matrix is indexed $i-1$ in every matrix $\mxex_{\n\setminus\{1\}}^{\n\setminus\{j\}}$, due to the removal of the first row of \mxex. Note also that, every column $k$ placed to the right of the $j$th column of matrix \mxex, after the removal of the $j$th column is shifted one position to the left, therefore it is indexed as $k-1$. Using the notation \indfun for the indicator function, after the above observations, we have
  \begin{align*}
    \vect{c}_i\vect{t} &= \sum_{j=1}^N c_{ij} (-1)^{1+j} \det \mxex_{\n\setminus\{1\}}^{\n\setminus\{j\}} \\
                       &= \sum_{j=1}^N c_{ij} (-1)^{1+j} \sum_{k\neq j} c_{ik} (-1)^{i-1+k-\indfun[\{k>j\}]}\det \mxex_{\n\setminus\{1,i\}}^{\n\setminus\{j,k\}}\\
                       &=  (-1)^i \sum_{j=1}^N \sum_{k\neq j} c_{ij} c_{ik}
                        (-1)^{j+k-\indfun[\{k>j\}]}\det \mxex_{\n\setminus\{1\}}^{\n\setminus\{j\}} =0,
  \end{align*}
  because for any two arbitrary columns $m$ and $l$, with $m>l$, only the summands
  \begin{align*}
    & c_{il} c_{im}(-1)^{l+m-1}\det \mxex_{\n\setminus\{1,i\}}^{\n\setminus\{l,m\}},\quad \text{and}\\
    & c_{im} c_{il}(-1)^{m+l}\det \mxex_{\n\setminus\{1,i\}}^{\n\setminus\{l,m\}},
  \end{align*}
  appear in the expression of $\vect{c}_i\vect{t}$ and they cancel out with one another. Since, all summands of the above double sum are coupled and canceled out, the double sum is equal to zero.
  Thus, we have proven that the inner product of any column of \mxex with \vect{t} is equal to zero. Consequently, \vect{t} is an eigenvector of matrix \mxex that corresponds to its eigenvalue zero.
\end{proof}

\begin{remark}\label{r3.observation on the eigenvectors}
  If the nullspace of an $N\times N$ matrix \mxex has dimension one, then $\rank \mxex = N-1$. Therefore, there exists at least one submatrix of \mxex such that its determinant is not equal to zero. More precisely, there exists at least one combination of row-column $(m,n)$ with $\det \mxex_{\n\setminus\{m\}}^{\n\setminus\{n\}} \neq 0$. Thus, if all determinants $\det \ind{\mxex}{\n\setminus\{j\}}{\n\setminus\{1\}}$, $j \in \n$, are equal to zero, we can choose the coordinates of the right eigenvector \vect{t}, which corresponds to the eigenvalue zero, as $t_j =(-1)^{m+j} \det \mxex_{\n\setminus\{m\}}^{\n\setminus\{j\}}$, $j \in \n$.
\end{remark}

\begin{remark}\label{r3.proportional eigenvectors}
   If \vect{t} is an arbitrary eigenvector that belongs to the nullspace of \mxex, then any other eigenvector \vect{z} that belongs to the same nullspace is proportional to \vect{t}. Namely, there exists $\sigma \in \mathds{R}$ such that $\vect{z} = \sigma \vect{t}$.
\end{remark}

From Theorem~\ref{t3.eigenvectors that correspond to the eigenvalue zero} and Remark~\ref{r3.observation on the eigenvectors}, we have as consequence the following corollary for the right eigenvectors of the matrix \mxe.
\begin{corollary}\label{c3.eigenvectors of E}
  If $m\in \n$ is such that $\det \ind{\mxe}{\n\setminus\{j\}}{\n\setminus\{m\}} \neq 0$ for at least one $j\in \n$, a right eigenvector \vect{t} of the nullspace of\/ \mxe has coordinates
  \begin{equation*}
    t_j =(-1)^{m+j} \det \ind{\mxe}{\n\setminus\{j\}}{\n\setminus\{m\}}, \qquad j \in \n.
  \end{equation*}
\end{corollary}

We now perturb the matrix function \ltm{\mxe} by $\epsilon \ltm{\mxk}$. Namely, we consider the matrix $\ltm{\mxe}+\epsilon \ltm{\mxk}$, where we assume that the matrix \ltm{\mxk} is analytic in the neighborhood of $r$. If \dm[n]{\mxk} is the $n$th derivative of the matrix function \ltm{\mxk} at $s=r$, the Taylor series of matrix \ltm{\mxk} around $r$ is:
\begin{equation}\label{e3.taylor series matrix K}
  \ltm{\mxk} = \mxk + (s-r) \dm[1]{\mxk} + \dots = \sum_{n=0}^\infty \frac{(s-r)^n}{n!} \dm{\mxk},
\end{equation}
where $\dm{\mxk} = \ltm[r]{\dm{\mxk}}$ and $\mxk=\dm[0]{\mxk}$. Our goal is to find the form of the eigenvectors of the nullspace of $\ltm{\mxe}+\epsilon \ltm{\mxk}$. Thus, as a first step we find the roots of the solution
\begin{equation}\label{e3.determinant equation perturbed}
  \det \big(\ltm{\mxe}+\epsilon \ltm{\mxk}\big)=0.
\end{equation}

At this point, we need the following result from perturbation theory, which gives us the root of a function $\fun{f}$ when it it perturbed by a small amount.

\begin{theorem}\label{t3.method to find the perturbed roots}
  Let $r$ be a simple root of an analytic function $\fun{f}$. For some function $\fun[s,\epsilon]{h}$ and for all small real values $\epsilon$, we define the perturbed function
  \begin{equation}\label{e3.perturbed function}
    F(s,\epsilon) = \fun{f} + \fun[s,\epsilon]{h}.
  \end{equation}
  If $\fun[s,\epsilon]{h}$ is analytic in $s$ and $\epsilon$ near $(r,0)$, then $F(s,\epsilon)$ has a unique simple root $(x(\epsilon),\epsilon)$ near $(r,0)$ for all small values of $\epsilon$. Moreover, $x(\epsilon)$ is an analytic function in $\epsilon$, and if $\frac{\partial}{\partial s^n}\fun[s,0]{h} \equiv 0$, $n=0,1,\dots$, then it holds
  \begin{equation}\label{perturbed root}
    \fun[\epsilon]{x}=r - \epsilon \frac{\frac{\partial}{\partial \epsilon}\fun[r,0]{h}}{\fun[r]{\df[1]{f}}} +O(\epsilon^2).
  \end{equation}
\end{theorem}
\begin{proof}
  From the Implicit function theorem \cite{amman-HCE}, we know that there exist a unique function $x$, with $\fun[0]{x}=r$, such that for all small values of $\epsilon$, it holds that $F\big(\fun[\epsilon]{x},\epsilon\big)=0$ close to $(r,0)$. Moreover, the function $x$ is analytic in $\epsilon$. To find the linear Taylor polynomial approximation of \fun[\epsilon]{x}, which is defined as
  \begin{equation*}
    \fun[\epsilon]{x} = \fun[0]{x} + \epsilon \fun[0]{\df[1]{x}} + O(\epsilon^2),
  \end{equation*}
  we differentiate the function $F\big(\fun[\epsilon]{x},\epsilon\big)=0$ as a function of $\epsilon$, and by using the chain rule we obtain
  \begin{gather*}
    \frac{\partial}{\partial \fun[\epsilon]{x}} F\big(\fun[\epsilon]{x},\epsilon\big) \fun[\epsilon]{\df[1]{x}} + \frac{\partial}{\partial \epsilon} F\big(\fun[\epsilon]{x},\epsilon\big) = 0 \\
    \Rightarrow \\
    \big(\fun[\fun[\epsilon]{\df[1]{f}}]{x}  + \frac{\partial}{\partial \fun[\epsilon]{x}} \fun[\fun[\epsilon]{h},\epsilon]{x} \big) \fun[\epsilon]{\df[1]{x}} +  \frac{\partial}{\partial \epsilon}\fun[\fun[\epsilon]{h},\epsilon]{x}  = 0.
  \end{gather*}
  In the latter equation, we substitute $\epsilon=0$ and we solve it with respect to \fun[0]{\df[1]{x}}. Since $r$ is a simple root the function $f$, it holds that $\fun[r]{\df[1]{f}}\neq 0$ \cite{krantz-HCV}. Thus, we have
  \begin{equation*}
    \fun[r]{\df[1]{f}}\fun[0]{\df[1]{x}} + \frac{\partial}{\partial \epsilon}\fun[r,0]{h} = 0 \
     \Rightarrow \
     \fun[0]{\df[1]{x}} = -\frac{\frac{\partial}{\partial \epsilon}\fun[r,0]{h}}{\fun[r]{\df[1]{f}}},
  \end{equation*}
 which completes the proof.
\end{proof}

From Theorem~\ref{t3.method to find the perturbed roots}, we have the following lemma.
\begin{lemma}\label{l3.perturbed function}
  If the functions \fun{f} and \fun[s,\epsilon]{h} satisfy the assumptions of Theorem~\ref{t3.method to find the perturbed roots}, and \fun{g} is an analytic function with $\fun[r]{g}\neq 0$, then the perturbed function
  \begin{equation*}
    G(s,\epsilon) = \fun{f}\fun{g} + \fun[s,\epsilon]{h}\fun{g},
  \end{equation*}
  has the same unique simple root $(x(\epsilon),\epsilon)$ near $(r,0)$ for all small values of $\epsilon$ with the perturbed function $F(s,\epsilon) = \fun{f} + \fun[s,\epsilon]{h}$. Namely
  \begin{equation*}
    \fun[\epsilon]{x}=r - \epsilon \frac{\frac{\partial}{\partial \epsilon}\fun[r,0]{h}}{\fun[r]{\df[1]{f}}} +O(\epsilon^2).
  \end{equation*}
\end{lemma}
\begin{proof}
  According to Theorem~\ref{t3.method to find the perturbed roots}, the unique simple root $\fun[\epsilon]{x}$ of $G(s,\epsilon)$ near $(r,0)$ for all small values of $\epsilon$ satisfies
  \begin{align*}
    \fun[\epsilon]{x} & = r - \epsilon
            \frac{\frac{\partial}{\partial \epsilon} \left. \big( \fun[s,\epsilon]{h}\fun{g}\big) \right|_{(r,0)}}{\left.\frac{\partial}{\partial s}\big(\fun[r]{f}\fun[r]{g}\big)\right|_{s=r}} +O(\epsilon^2)\\
                & = r - \epsilon \frac{\frac{\partial}{\partial \epsilon}\fun[r,0]{h}\fun[r]{g}}{\fun[r]{\df[1]{f}} \fun[r]{g} + \fun[r]{f}\fun[r]{\df[1]{g}}} +O(\epsilon^2)\\
                & = r - \epsilon \frac{\frac{\partial}{\partial \epsilon}\fun[r,0]{h}}{\fun[r]{\df[1]{f}}} +O(\epsilon^2),
  \end{align*}
  because $\fun[r]{f}=0$.
\end{proof}

We also need the following property for the determinant of a square matrix.
\begin{proposition}\label{p3.determinant property}
  If \mxex and \mxexs are $N\times N$ matrices with columns $\mxex_{\bullet i}$ and $\mxexs_{\bullet i}$, $i \in \n$, respectively, then
  \begin{align*}
    \det(\mxex_{\bullet 1} & +\epsilon \mxexs_{\bullet 1},\dots,\mxex_{\bullet N}+\epsilon \mxexs_{\bullet N}) = \underbrace{\det(\mxex_{\bullet 1},\dots,\mxex_{\bullet N})}_{\det(\mxex)} \\
    &+ \epsilon \sum_{i=1}^N\det(\mxex_{\bullet 1},\dots,\mxexs_{\bullet i},\dots,\mxex_{\bullet N}) +O(\epsilon^2).
  \end{align*}
\end{proposition}
\begin{proof}
  The result is an immediate consequence of the additive property of determinants.
\end{proof}

As shown in the following corollary, we can find the roots of the equation $\det \big(\ltm{\mxe} + \epsilon\ltm{\mxk}\big) = 0$, combining the results of Theorem~\ref{t3.method to find the perturbed roots} and Proposition~\ref{p3.determinant property}.

\begin{corollary}\label{c3.roots of the perturbed matrix}
  The number $r_\epsilon=r - \epsilon \delta + O(\epsilon^2) $, where
  \begin{equation*}
    \delta = \frac{\sum_{j=1}^N \det\big(\mxe[\bullet 1],\dots,\mxk[\bullet j],\dots,\mxe[\bullet N]\big)}
                    {\sum_{j=1}^N \det\big(\mxe[\bullet 1],\dots,\dm[1]{\mxe[\bullet j]},\dots,\mxe[\bullet N]\big)},
  \end{equation*}
  is a simple root of the determinant $\det \big(\ltm{\mxe} + \epsilon\ltm{\mxk}\big)=0$.
\end{corollary}
\begin{proof}%
  According to Proposition~\ref{p3.determinant property},
  \begin{align*}
    \det \big(\ltm{\mxe} &+ \epsilon \ltm{\mxk}\big)
        = \det \ltm{\mxe} \\
        &+ \epsilon \sum_{j=1}^N \det\big(\ltm{\mxe[\bullet 1]},\dots,\ltm{\mxk[\bullet j]},\dots,\ltm{\mxe[\bullet N]} \big) \\
        &+ O(\epsilon^2).
  \end{align*}
  Note that $\det \ltm{\mxe} $ is an analytic function in $r$ and its derivative is defined as
  \begin{align*}
    \frac{d}{ds}\det\ltm{\mxe} &= \sum_{j=1}^N
                               \det\big(\ltm{\mxe[\bullet 1]},\dots,\ltm{\dm[1]{\mxe[\bullet j]}},\dots,
                               \ltm{\mxe[\bullet N]} \big).
  \end{align*}
  Since $r$ is a simple eigenvalue of \ltm{\mxe}, by the definition of the multiplicity of a root of an analytic function, it holds that $\frac{d}{ds}\det\ltm{\mxe} \left.\right|_{s=r} \neq 0$ (see \cite{krantz-HCV}). In addition, the function $\sum_{j=1}^N \det\big(\ltm{\mxe[\bullet 1]},\dots,\ltm{\mxk[\bullet j]},\dots,\ltm{\mxe[\bullet N]} \big)$ is also analytic in the neighborhood of $r$. The result is then immediate from Theorem~\ref{t3.method to find the perturbed roots}.
\end{proof}

According to Corollary~\ref{c3.roots of the perturbed matrix}, the eigenvalue $r_\epsilon$ of the matrix $\ltm{\mxe} + \epsilon \ltm{\mxk}$ is simple. Consequently, the dimension of the nullspace of each matrix ${\mxe} + \epsilon \mxk$ is equal to one. We apply Theorem~\ref{t3.eigenvectors that correspond to the eigenvalue zero} to find the eigenvectors of the matrix ${\mxe} + \epsilon \mxk$, that correspond to its eigenvalue $r_\epsilon$. Before that though, we do the following simplification. From Eqs.~\eqref{e3.taylor series matrix E}--\eqref{e3.taylor series matrix K} we have the Taylor expansion
\begin{align*}
  \ltm{\mxe} + \epsilon \ltm{\mxk}
            &= \sum_{n=0}^\infty \frac{(s-r)^n}{n!} \Big( \dm{\mxe} + \epsilon \dm{\mxk} \Big).
  \intertext{Evaluating this at the point $r_\epsilon = r - \epsilon \delta + O(\epsilon^2)$, we obtain}
  \ltm[r_\epsilon]{\mxe} + \epsilon \ltm[r_\epsilon]{\mxk}
                    & = \mxe - \epsilon \delta \dm[1]{\mxe} + \epsilon  \mxk + O(\epsilon^2 \um)\\
                    & = \mxe + \epsilon \big(\mxk - \delta \dm[1]{\mxe} \big) +O(\epsilon^2 \um),
\end{align*}
where we denote by \um the matrix with all its elements equal to one.% Since the same analysis applies for all different values \rootpmix{i}, $i=2,\dots,N$, .

\begin{theorem}\label{t3.eigenvectors perturbed system}
  A right eigenvector of matrix $\mxe + \epsilon \big( \mxk - \delta \dm[1]{\mxe} \big)$ that corresponds to its eigenvalue $r_\epsilon$ is
  \begin{equation*}
    \vect{w}= \vect{t} -\epsilon \delta \dm[1]{\vect{t}} + \epsilon \vect{k},
  \end{equation*}
  where \vect{t} is a right eigenvector of \mxe defined as in Corollary~\ref{c3.eigenvectors of E} and \dm[1]{\vect{t}} is its derivative. Moreover, \vect{k} is an $N\times 1$ vector with coordinates
  \begin{align*}
    k_j &= (-1)^{m+j}\sum_{k=1}^{N-1}\det \Bigg(\Big(\ind{\mxe}{\n\setminus\{j\}}{\n\setminus\{m\}}\Big)_{\bullet 1},\dots,\\
    &\Big(\ind{\mxk}{\n\setminus\{j\}}{\n\setminus\{m\}}\Big)_{\bullet k},
            \dots,\Big(\ind{\mxe}{\n\setminus\{j\}}{\n\setminus\{m\}}\Big)_{\bullet N-1} \Bigg), \qquad j \in \n,
  \end{align*}
  where the choice of $m\in \n$ is explained in Corollary~\ref{c3.eigenvectors of E}.
\end{theorem}
\begin{proof}
  According to Remark~\ref{r3.observation on the eigenvectors} and Corollary~\ref{c3.eigenvectors of E}, there exists an $m \in \n$ such that the vector \vect{t} with coordinates
  \begin{equation*}
    t_j = (-1)^{m+j} \det \ind{\mxe}{\n\setminus\{j\}}{\n\setminus\{m\}}, \qquad j \in \n,
  \end{equation*}
  is a right eigenvector of matrix \mxe. We prove that a right eigenvector that corresponds to the matrix $\mxe + \epsilon \big(\mxk - \delta \dm[1]{\mxe}\big)$ is \vect{w} with coordinates
  \begin{equation*}
    w_j = (-1)^{m+j} \det \big(\mxe + \epsilon \big( \mxk - \delta \dm[1]{\mxe} \big)\big)_{\n\setminus\{m\}}^{\n\setminus\{j\}}, \qquad j \in \n.
  \end{equation*}
  Using Proposition~\ref{p3.determinant property}, the above equation simplifies to
  \begin{small}
  \begin{align*}
    w_j =& \ (-1)^{m+j} \det \Big(\mxe + \epsilon \big(\mxk - \delta \dm[1]{\mxe}\big)\Big)_{\n\setminus\{m\}}^{\n\setminus\{j\}}\\
        =& \ (-1)^{m+j} \det\ind{\mxe}{\n\setminus\{j\}}{\n\setminus\{m\}} \\
         &+ \epsilon (-1)^{m+j}\sum_{k=1}^{N-1}
          \det \Bigg(\Big(\ind{\mxe}{\n\setminus\{j\}}{\n\setminus\{m\}}\Big)_{\bullet 1},\dots, \\
         &\Big(\ind{\big( \mxk - \delta \dm[1]{\mxe} \big)}{\n\setminus\{j\}}{\n\setminus\{m\}}\Big)_{\bullet k},
          \dots,\Big(\ind{\mxe}{\n\setminus\{j\}}{\n\setminus\{m\}}\Big)_{\bullet N-1} \Bigg)\\
        =& \ (-1)^{m+j} \ind{\mxe}{\n\setminus\{j\}}{\n\setminus\{m\}} \\
         &- \epsilon (-1)^{m+j}\delta \sum_{k=1}^{N-1} \det \Bigg(\Big(\ind{\mxe}{\n\setminus\{j\}}{\n\setminus\{m\}}\Big)_{\bullet 1},\dots, \\ &\Big(\ind{\dm[1]{\mxe}}{\n\setminus\{j\}}{\n\setminus\{m\}}\Big)_{\bullet k},
          \dots,\Big(\ind{\mxe}{\n\setminus\{j\}}{\n\setminus\{m\}}\Big)_{\bullet N-1} \Bigg)\\
         &+\ \epsilon (-1)^{m+j} \sum_{k=1}^{N-1}\det \Bigg(\Big(\ind{\mxe}{\n\setminus\{j\}}{\n\setminus\{m\}}\Big)_{\bullet 1},\dots,\\
         &\Big(\ind{\mxk}{\n\setminus\{j\}}{\n\setminus\{m\}}\Big)_{\bullet k},\dots,\Big(\ind{\mxe}{\n\setminus\{j\}}{\n\setminus\{m\}}\Big)_{\bullet N-1} \Bigg)\\
        =& \  t_j -\epsilon \delta \dm[1]{t}_j + \epsilon k_j,
  \end{align*}
  \end{small}
  where $\dm[1]{t}_j = \left.\frac{d}{ds}\fun{t_j}\right|_{s=r}$ and $k_j =(-1)^{m+j}\sum_{k=1}^{N-1}\det \allowbreak \Bigg(\Big(\ind{\mxe}{\n\setminus\{j\}}{\n\setminus\{m\}}\Big)_{\bullet 1},\allowbreak \dots,\Big(\ind{\mxk}{\n\setminus\{j\}}{\n\setminus\{m\}}\Big)_{\bullet k},\dots,
  \Big(\ind{\mxe}{\n\setminus\{j\}}{\n\setminus\{m\}}\Big)_{\bullet N-1} \Bigg)$.

  Observe that \vect{t} is not identically equal to zero, because it is an eigenvector of \mxe. Thus, the vector \vect{w} is also not identically equal to zero. Therefore, according to Remark~\ref{r3.observation on the eigenvectors}, \vect{w} is an eigenvector of the matrix $\mxe + \epsilon \big( \mxk - \delta \dm[1]{\mxe}\big)$, which completes the proof.
\end{proof}

\section{Proofs}\label{appendix B}
\begin{proof}[of Theorem~\ref{t3.determinant matrix E}]
  To prove the theorem, we need formulas that result from the properties of the determinants. We define the sets $F_i=\{1,\dots,i\}$ and $L_i=\{i,\dots,N\}$, where $F_0=L_{N+1}=\emptyset$. Using the additive property of determinants and by expanding in minors on the first row, we obtain for $i \in \n$,
  \begin{align*}
    \det \ltm{\mxe}^{L_{i}}_{L_{i}} =& \ltgs \lambda_{i} \det\Big( \big(\mxprobsreal\circ\mxtrans\big)^{\{i\}}_{L_{i}}, \ltm{\mxe}^{L_{i+1}}_{L_{i}} \Big) \notag\\
                                     &+ \lambda_{i} \det\Big( \big(\mxprobsdummy\circ\mxtrans\big)^{\{i\}}_{L_{i}}, \ltm{\mxe}^{L_{i+1}}_{L_{i}} \Big)  \notag \\
                                     &+ (s-\lambda_i) \det \ltm{\mxe}^{L_{i+1}}_{L_{i+1}}. %\label{e3.resursion for matrix E, 1}
  \end{align*}
  Suppose now that $V=\{i_1,\dots,i_n\}$ and $W=\{j_1,\dots,j_k\}$ are two non\-/overlapping ($V\cap W =\emptyset$) collections of $n$ and $k$ elements from \n, respectively, with $1\leq n+k \leq N-1$. Furthermore, we choose $j$ such that $j>\max\{l:l\in V\cup W\}$. Then, the determinant of the $(N+1-j+n+k)$-dimension square matrix $\Big(\big(\mxprobsdummy\circ\mxtrans\big)^V_{V\cup W\cup L_{j}}\Join \allowbreak\big(\mxprobsreal\circ\mxtrans\big)^W_{V\cup W\cup L_{j}}, \ltm{\mxe}^{L_{j}}_{V\cup W\cup L_{j}}\Big)$  satisfies,
  \begin{small}
  \begin{align*}
    \det & \Big(\big(\mxprobsdummy\circ\mxtrans\big)^V_{V\cup W\cup L_{j}}\Join \big(\mxprobsreal\circ\mxtrans\big)^W_{V\cup W\cup L_{j}}, \notag \\
            &\ltm{\mxe}^{L_{j}}_{V\cup W\cup L_{j}}\Big)\notag \\
           =& \ltgs \lambda_j \det \Big(\big(\mxprobsdummy\circ\mxtrans\big)^V_{V\cup W\cup L_{j}}\Join \big(\mxprobsreal\circ\mxtrans\big)^{W\cup \{j\}}_{V\cup W\cup L_{j}}, \notag \\
            &\ltm{\mxe}^{L_{j+1}}_{V\cup W\cup L_{j}}\Big) \notag \\
            &+ \lambda_j \det \Big(\big(\mxprobsdummy\circ\mxtrans\big)^{V\cup \{j\}}_{V\cup W\cup L_{j}}\Join \big(\mxprobsreal\circ\mxtrans\big)^{W}_{V\cup W\cup L_{j}}, \notag \\
            &\ltm{\mxe}^{L_{j+1}}_{V\cup W\cup L_{j}}\Big) \notag \\
            &+ (s-\lambda_j) \det \Big(\big(\mxprobsdummy\circ\mxtrans\big)^{V}_{V\cup W\cup L_{j+1}}\Join  \notag \\
            & \big(\mxprobsreal\circ\mxtrans\big)^{W}_{V\cup W\cup L_{j+1}},\ltm{\mxe}^{L_{j+1}}_{V\cup W\cup L_{j+1}}\Big). %\label{e3.resursion for matrix E, 2}
  \end{align*}
  \end{small}
  Note that $\det \ltm{\mxe} = \det \ltm{\mxe}^{L_{1}}_{L_{1}}$. The theorem is proven by applying recursively the above formulas.
\end{proof}

\begin{proof}[of Theorem~\ref{t3.adjoint matrix of E}]
  It is known that
  \begin{equation*}
    \ltm{\mxadj[ij]} = (-1)^{i+j} \det \ltm{\mxe}^{\n\setminus\{i\}}_{\n\setminus\{j\}}.
  \end{equation*}
  The case $i=j$ is merely an application of Theorem~\ref{t3.determinant matrix E}, where instead of state space \n we have $\n\setminus\{i\}$. Therefore,
  \begin{align*}
    \ltm{\mxadj[ii]} =& \sum_{S\subset \n\setminus \{i\}}\las{S}\fun{\zeta^{S^c}}\det \big(\mxprobsdummy\circ\mxtrans\big)^S_S \\
                     +& \sum_{k=1}^{N-1} \ltgs[k] \sum_{\substack{\Gamma\subset \n\setminus \{i\}\\  \mid \Gamma\mid =k}}  \sum_{\substack{S\subset \n\setminus \{i\}\\ S\supset \Gamma}}
                        \las{S}  \fun{\zeta^{S^c}} \\
                      & \times \det \Big( \big(\mxprobsdummy\circ\mxtrans\big)^{S\setminus \Gamma}_S  \Join \big(\mxprobsreal\circ\mxtrans\big)^\Gamma_S\Big).
  \end{align*}
  When $i\neq j$, we need to separate the two cases $i<j$ and $i>j$. We first deal with the case $i<j$. We then have,
  \begin{small}
  \begin{align*}
    \ltm{\mxadj[ij]} =& (-1)^{i+j} \ltgs \lambda_j \\
                      & \times \det \Big( \ltm{\mxe}^{F_{j-1}\setminus\{i\}}_{\n\setminus\{j\}}, \big(\mxprobsreal\circ\mxtrans\big)^{\{j\}}_{\n\setminus\{j\}},\ltm{\mxe}^{L_{j+1}}_{\n\setminus\{j\}}\Big)\\
                     +& (-1)^{i+j} \lambda_j \\
                      & \times \det \Big( \ltm{\mxe}^{F_{j-1}\setminus\{i\}}_{\n\setminus\{j\}}, \big(\mxprobsdummy\circ\mxtrans\big)^{\{j\}}_{\n\setminus\{j\}},\ltm{\mxe}^{L_{j+1}}_{\n\setminus\{j\}}\Big).
  \end{align*}
  \end{small}
  We find \ltm{\mxadj[ij]} by expanding the determinants that appear above in minors of their first row. For this reason, it is important to know what is the position of the elements $\ltm{\mxe[n,n]} = \ltgsc{nn} p_{nn}\lambda_n+s-\lambda_n$, $n \in \n\setminus\{i,j\}$, in the above reduced matrix. Note that the elements \ltm{\mxe[n,n]} with $n=i+1,\dots,j-1$, are on the diagonal of matrix \ltm{\mxe}. However, when $j\neq i+1$ they drop to the lower-diagonal of the matrices
  \begin{align*}
    \Big( \ltm{\mxe}^{F_{j-1}\setminus\{i\}}_{\n\setminus\{j\}}, \big(\mxprobsreal\circ\mxtrans\big)^{\{j\}}_{\n\setminus\{j\}},\ltm{\mxe}^{L_{j+1}}_{\n\setminus\{j\}}\Big),
  \intertext{and}
    \Big( \ltm{\mxe}^{F_{j-1}\setminus\{i\}}_{\n\setminus\{j\}}, \big(\mxprobsdummy\circ\mxtrans\big)^{\{j\}}_{\n\setminus\{j\}},\ltm{\mxe}^{L_{j+1}}_{\n\setminus\{j\}}\Big).
  \end{align*}
  It is immediately obvious that if this displacement takes place, it will result in a change of sign for the determinants. For this reason, we split the columns of the latter matrices in the subsets $F_{i-1}$, $T$, $\{j\}$ and $L_{j+1}$, where $T=\{i+1,\dots,j-1\}$. We fix some $m\in \n\setminus\{i,j\}$ and we separate the following cases:
  \begin{enumerate}
    \item $m \in F_{i-1}$. For every two non\-/overlapping collections of $n$ and $k$ elements from $F_{m-1}$, say $V=\{i_1,\dots,i_n\}$ and $W=\{j_1,\dots,j_k\}$, with $1\leq n+k \leq m-1$, it holds that
        \begin{small}
        \begin{align*}
             \det &\Big( \big(\mxprobsdummy\circ\mxtrans\big)^V_\Omega\Join \big(\mxprobsreal\circ\mxtrans\big)^W_\Omega,\ltm{\mxe}^{(L_m\cap F_{i-1})\cup T}_\Omega,\\
                  & \big(\mxprobsreal\circ\mxtrans\big)^{\{j\}}_\Omega,\ltm{\mxe}^{L_{j+1}}_\Omega\Big) \\
                 =& \ltgs \lambda_{m} \det \Big( \big(\mxprobsdummy\circ\mxtrans\big)^V_\Omega\Join \big(\mxprobsreal\circ\mxtrans\big)^{W \cup \{m\}}_\Omega,\\
                  & \ltm{\mxe}^{(L_{m+1}\cap F_{i-1})\cup T}_\Omega, \big(\mxprobsreal\circ\mxtrans\big)^{\{j\}}_\Omega,\ltm{\mxe}^{L_{j+1}}_\Omega\Big)\\
                 +& \lambda_{m} \det \Big( \big(\mxprobsdummy\circ\mxtrans\big)^{V \cup \{m\}}_\Omega\Join \big(\mxprobsreal\circ\mxtrans\big)^{W }_\Omega,\\
                  & \ltm{\mxe}^{(L_{m+1}\cap F_{i-1})\cup T}_\Omega, \big(\mxprobsreal\circ\mxtrans\big)^{\{j\}}_\Omega,\ltm{\mxe}^{L_{j+1}}_\Omega\Big)\\
                 +& (s-\lambda_{m}) \det \Big( \big(\mxprobsdummy\circ\mxtrans\big)^{V}_{\Omega\setminus \{m\}} \Join \big(\mxprobsreal\circ\mxtrans\big)^{W }_{\Omega\setminus \{m\}},\\
                  &\ltm{\mxe}^{(L_{m+1}\cap F_{i-1})\cup T}_{\Omega\setminus \{m\}}, \big(\mxprobsreal\circ\mxtrans\big)^{\{j\}}_{\Omega\setminus \{m\}},\ltm{\mxe}^{L_{j+1}}_{\Omega\setminus \{m\}}\Big),
        \end{align*}
        \end{small}
        and,
        \begin{small}
        \begin{align*}
             \det &\Big( \big(\mxprobsdummy\circ\mxtrans\big)^V_\Omega\Join \big(\mxprobsreal\circ\mxtrans\big)^W_\Omega,\ltm{\mxe}^{(L_m\cap F_{i-1})\cup T}_\Omega,\\
                  & \big(\mxprobsdummy\circ\mxtrans\big)^{\{j\}}_\Omega,\ltm{\mxe}^{L_{j+1}}_\Omega\Big) \\
                 =& \ltgs \lambda_{m} \det \Big( \big(\mxprobsdummy\circ\mxtrans\big)^V_\Omega\Join \big(\mxprobsreal\circ\mxtrans\big)^{W \cup \{m\}}_\Omega,\\
                  & \ltm{\mxe}^{(L_{m+1}\cap F_{i-1})\cup T}_\Omega, \big(\mxprobsdummy\circ\mxtrans\big)^{\{j\}}_\Omega,\ltm{\mxe}^{L_{j+1}}_\Omega\Big)\\
                 +& \lambda_{m} \det \Big( \big(\mxprobsdummy\circ\mxtrans\big)^{V \cup \{m\}}_\Omega\Join \big(\mxprobsreal\circ\mxtrans\big)^{W }_\Omega,\\
                  & \ltm{\mxe}^{(L_{m+1}\cap F_{i-1})\cup T}_\Omega, \big(\mxprobsdummy\circ\mxtrans\big)^{\{j\}}_\Omega,\ltm{\mxe}^{L_{j+1}}_\Omega\Big)\\
                 +& (s-\lambda_{m}) \det \Big( \big(\mxprobsdummy\circ\mxtrans\big)^{V}_{\Omega\setminus \{m\}} \Join \big(\mxprobsreal\circ\mxtrans\big)^{W }_{\Omega\setminus \{m\}},\\
                  &\ltm{\mxe}^{(L_{m+1}\cap F_{i-1})\cup T}_{\Omega\setminus \{m\}}, \big(\mxprobsdummy\circ\mxtrans\big)^{\{j\}}_{\Omega\setminus \{m\}},\ltm{\mxe}^{L_{j+1}}_{\Omega\setminus \{m\}}\Big),
        \end{align*}
        \end{small}
        where $\Omega={V\cup W\cup(L_m\cap F_{i-1})\cup\{i\}\cup T\cup L_{j+1}}$.
    \item $m \in T$ with $T\neq \emptyset$ (note that $T\neq \emptyset$ when $j\neq i+1$). For every two non\-/overlapping collections of $n$ and $k$ elements from $F_{m-1}\setminus \{i\}$, say $V=\{i_1,\dots,i_n\}$ and $W=\{j_1,\dots,j_k\}$, with $1\leq n+k \leq m-2$, it holds that
        \begin{small}
        \begin{align*}
            \det &\Big( \big(\mxprobsdummy\circ\mxtrans\big)^V_\Omega\Join \big(\mxprobsreal\circ\mxtrans\big)^W_\Omega,\ltm{\mxe}^{L_m\cap T}_\Omega, \\
                 & \big(\mxprobsreal\circ\mxtrans\big)^{\{j\}}_\Omega,\ltm{\mxe}^{L_{j+1}}_\Omega\Big) \\
                =& \ltgs \lambda_{m} \det \Big( \big(\mxprobsdummy\circ\mxtrans\big)^V_\Omega\Join \big(\mxprobsreal\circ\mxtrans\big)^{W \cup \{m\}}_\Omega,\\
                 &\ltm{\mxe}^{L_{m+1}\cap T}_\Omega, \big(\mxprobsreal\circ\mxtrans\big)^{\{j\}}_\Omega,\ltm{\mxe}^{L_{j+1}}_\Omega\Big)\\
                +& \lambda_{m} \det \Big( \big(\mxprobsdummy\circ\mxtrans\big)^{V \cup \{m\}}_\Omega\Join \big(\mxprobsreal\circ\mxtrans\big)^{W }_\Omega,\\
                 &\ltm{\mxe}^{L_{m+1}\cap T}_\Omega, \big(\mxprobsreal\circ\mxtrans\big)^{\{j\}}_\Omega,\ltm{\mxe}^{L_{j+1}}_\Omega\Big)\\
                -&(s-\lambda_{m}) \det \Big( \big(\mxprobsdummy\circ\mxtrans\big)^{V}_{\Omega\setminus \{m\}} \Join \big(\mxprobsreal\circ\mxtrans\big)^{W }_{\Omega\setminus \{m\}},\\
                 &\ltm{\mxe}^{L_{m+1}\cap T}_{\Omega\setminus \{m\}}, \big(\mxprobsreal\circ\mxtrans\big)^{\{j\}}_{\Omega\setminus \{m\}},\ltm{\mxe}^{L_{j+1}}_{\Omega\setminus \{m\}}\Big),
        \end{align*}
        \end{small}
        and
        \begin{small}
        \begin{align*}
            \det &\Big( \big(\mxprobsdummy\circ\mxtrans\big)^V_\Omega\Join \big(\mxprobsreal\circ\mxtrans\big)^W_\Omega,\ltm{\mxe}^{L_m\cap T}_\Omega, \\
                 & \big(\mxprobsdummy\circ\mxtrans\big)^{\{j\}}_\Omega,\ltm{\mxe}^{L_{j+1}}_\Omega\Big) \\
                =& \ltgs \lambda_{m} \det \Big( \big(\mxprobsdummy\circ\mxtrans\big)^V_\Omega\Join \big(\mxprobsreal\circ\mxtrans\big)^{W \cup \{m\}}_\Omega,\\
                 &\ltm{\mxe}^{L_{m+1}\cap T}_\Omega, \big(\mxprobsdummy\circ\mxtrans\big)^{\{j\}}_\Omega,\ltm{\mxe}^{L_{j+1}}_\Omega\Big)\\
                +& \lambda_{m} \det \Big( \big(\mxprobsdummy\circ\mxtrans\big)^{V \cup \{m\}}_\Omega\Join \big(\mxprobsreal\circ\mxtrans\big)^{W }_\Omega,\\
                 &\ltm{\mxe}^{L_{m+1}\cap T}_\Omega, \big(\mxprobsdummy\circ\mxtrans\big)^{\{j\}}_\Omega,\ltm{\mxe}^{L_{j+1}}_\Omega\Big)\\
                -&(s-\lambda_{m}) \det \Big( \big(\mxprobsdummy\circ\mxtrans\big)^{V}_{\Omega\setminus \{m\}} \Join \big(\mxprobsreal\circ\mxtrans\big)^{W }_{\Omega\setminus \{m\}},\\
                 &\ltm{\mxe}^{L_{m+1}\cap T}_{\Omega\setminus \{m\}}, \big(\mxprobsdummy\circ\mxtrans\big)^{\{j\}}_{\Omega\setminus \{m\}},\ltm{\mxe}^{L_{j+1}}_{\Omega\setminus \{m\}}\Big),
        \end{align*}
        \end{small}
    where $\Omega={V\cup W\cup\{i\}\cup(L_m\cap T)\cup L_{j+1}}$.
    \item $m \in L_{j+1}$. For every two non\-/overlapping collections of $n$ and $k$ elements from $F_{m-1}\setminus \{i\}$, say $V=\{i_1,\dots,i_n\}$ and $W=\{j_1,\dots,j_k\}$, with $1\leq n+k \leq m-2$, it holds that
        \begin{small}
        \begin{align*}
          \det &\Big( \big(\mxprobsdummy\circ\mxtrans\big)^V_\Omega\Join \big(\mxprobsreal\circ\mxtrans\big)^W_\Omega,\ltm{\mxe}^{L_{m}}_\Omega\Big)\\
            =& \ltgs \lambda_{m} \det \Big( \big(\mxprobsdummy\circ\mxtrans\big)^V_\Omega\Join \big(\mxprobsreal\circ\mxtrans\big)^{W \cup \{m\}}_\Omega,\\
             & \ltm{\mxe}^{L_{m+1}}_\Omega\Big)\\
            +& \lambda_{m} \det \Big( \big(\mxprobsdummy\circ\mxtrans\big)^{V\cup \{m\}}_\Omega\Join \big(\mxprobsreal\circ\mxtrans\big)^{W }_\Omega,\\
             & \ltm{\mxe}^{L_{m+1}}_\Omega\Big)\\
            +&(s-\lambda_{m}) \det \Big( \big(\mxprobsdummy\circ\mxtrans\big)^{V}_{\Omega\setminus \{m\}} \Join \big(\mxprobsreal\circ\mxtrans\big)^{W }_{\Omega\setminus \{m\}},\\
             &\ltm{\mxe}^{L_{m+1}}_{\Omega\setminus \{m\}}\Big),
        \end{align*}
        \end{small}
    where $\Omega=V\cup W\cup L_{m}$.
  \end{enumerate}

  Using the above formulas to evaluate all the involved determinants, we find that
  \begin{small}
  \begin{align*}
    \ltm{\mxadj[ij]} =& (-1)^{i+j} \ltgs \sum_{k=0}^{N-2} \ltgs[k] \sum_{\substack{\Gamma \subset \n\setminus\{i,j\}\\ \mid \Gamma\mid =k;\\S\subset \n\setminus\{i,j\}\\ S\supset \Gamma;\\
                            R\subset S\cap T}}(-1)^{\mid R\mid} \las{S\cup \{j\}}\\
                      & \times  \fun{\zeta^{S^c}} \det \Big( \big(\mxprobsdummy\circ\mxtrans\big)^{S\setminus \Gamma}_{S\cup \{i\}} \Join
                            \big(\mxprobsreal\circ\mxtrans\big)^{\Gamma \cup \{j\}}_{S\cup \{i\}}  \Big)\\
                     +& (-1)^{i+j} \sum_{k=0}^{N-2} \ltgs[k] \sum_{\substack{\Gamma \subset \n\setminus\{i,j\}\\  \mid \Gamma\mid =k; \\ S\subset \n\setminus\{i,j\}\\ S\supset \Gamma;\\
                            R\subset S\cap T}} (-1)^{\mid R\mid} \las{S\cup \{j\}} \fun{\zeta^{S^c}} \\
                      & \times  \det \Big( \big(\mxprobsdummy\circ\mxtrans\big)^{(S\setminus \Gamma)\cup \{j\}}_{S\cup \{i\}} \Join
                            \big(\mxprobsreal\circ\mxtrans\big)^{\Gamma}_{S\cup \{i\}}  \Big),\\
  \end{align*}
  \end{small}
  which holds even when $T=\emptyset$.

  We assume now that $i>j$, and we have to calculate
  \begin{small}
  \begin{align*}
    \ltm{\mxadj[ij]} =& (-1)^{i+j} \ltgs \lambda_j \\
                      & \times \det \Big( \ltm{\mxe}^{F_{j-1}}_{\n\setminus\{j\}}, \big(\mxprobsreal\circ\mxtrans\big)^{\{j\}}_{\n\setminus\{j\}},\ltm{\mxe}^{L_{j+1}\setminus\{i\}}_{\n\setminus\{j\}}\Big)\\
                     +& (-1)^{i+j} \lambda_j \\
                      & \times \det \Big( \ltm{\mxe}^{F_{j-1}}_{\n\setminus\{j\}}, \big(\mxprobsdummy\circ\mxtrans\big)^{\{j\}}_{\n\setminus\{j\}},\ltm{\mxe}^{L_{j+1}\setminus\{i\}}_{\n\setminus\{j\}}\Big).
  \end{align*}
  \end{small}
  In this case, $T=\{j+1,\dots,i-1\}$. When $T\neq \emptyset$, the elements $\ltm{\mxe[n,n]} = \ltgsc{nn} p_{nn}\lambda_n+s-\lambda_n$, with $n=j+1,\dots,i-1$, which are on the diagonal of matrix \ltm{\mxe}, move to the upper-diagonal of the matrices
  \begin{align*}
    \Big( \ltm{\mxe}^{F_{j-1}}_{\n\setminus\{j\}},\big(\mxprobsreal\circ\mxtrans\big)^{\{j\}}_{\n\setminus\{j\}},\ltm{\mxe}^{L_{j+1}\setminus\{i\}}_{\n\setminus\{j\}}\Big),
  \intertext{and}
    \Big( \ltm{\mxe}^{F_{j-1}}_{\n\setminus\{j\}},\big(\mxprobsdummy\circ\mxtrans\big)^{\{j\}}_{\n\setminus\{j\}},\ltm{\mxe}^{L_{j+1}\setminus\{i\}}_{\n\setminus\{j\}}\Big)
  \end{align*}

  The formula is exactly the same, with $T=\{i+1,\dots,j-1\}$. Thus, gathering all the above, for $i\neq j$
  \begin{small}
  \begin{align*}
    \ltm{\mxadj[ij]} =& (-1)^{i+j} \sum_{k=1}^{N-1} \ltgs[k] \sum_{\substack{\Gamma \subset \n\setminus\{i,j\}\\ \mid \Gamma\mid =k-1;\\S\subset \n\setminus\{i,j\}\\ S\supset \Gamma;\\
                            R\subset S\cap T_{ij}}}(-1)^{\mid R\mid} \las{S\cup \{j\}} \fun{\zeta^{S^c}} \\
                      & \times  \det \Big( \big(\mxprobsdummy\circ\mxtrans\big)^{S\setminus \Gamma}_{S\cup \{i\}} \Join
                            \big(\mxprobsreal\circ\mxtrans\big)^{\Gamma \cup \{j\}}_{S\cup \{i\}}  \Big)\\
                     +& (-1)^{i+j} \sum_{k=0}^{N-2} \ltgs[k] \sum_{\substack{\Gamma \subset \n\setminus\{i,j\}\\  \mid \Gamma\mid =k; \\ S\subset \n\setminus\{i,j\}\\ S\supset \Gamma;\\
                            R\subset S\cap T_{ij}}} (-1)^{\mid R\mid} \las{S\cup \{j\}} \fun{\zeta^{S^c}} \\
                      & \times  \det \Big( \big(\mxprobsdummy\circ\mxtrans\big)^{(S\setminus \Gamma)\cup \{j\}}_{S\cup \{i\}} \Join
                            \big(\mxprobsreal\circ\mxtrans\big)^{\Gamma}_{S\cup \{i\}}  \Big),\\
  \end{align*}
  \end{small}
  where $m_{ij}=\min\{i,j\}$, $M_{ij}=\max\{i,j\}$ and $T_{ij}=\{m_{ij}+1,\dots,M_{ij}-1\}$.
\end{proof}

\begin{proof}[of Theorem~\ref{t3.numerator LT workload state i}]
Observe that
\begin{align*}
  s\vup \ltm{\mxadj}\uv_i =& s\sum_{l=1}^N u_l \ltm{\mxadj[li]}
                         = s\sum_{\substack{l=1\\ l\neq i}}^N u_l \ltm{\mxadj[li]} +s u_i \ltm{\mxadj[ii]}.
\end{align*}
Using the definition of \ltm{\mxadj[ij]}, $\forall i,j \in \n$, and Theorem~\ref{t3.adjoint matrix of E}, the result is straightforward.
\end{proof}

\begin{proof}[of Proposition~\ref{p3.LT replace workload state i}]
In this case, the determinant $\det \ltm{\mxe}$ (see Theorem~\ref{t3.determinant matrix E}) takes the form
\begin{align}
  \det \ltm{\mxe} &= \sum_{S\subset \n}\las{S}\fun{\zeta^{S^c}}\det \big(\mxprobsdummy\circ\mxtrans\big)^S_S \notag\\
                  +& \sum_{k=1}^N \Bigg( \frac{\fun{q}}{\fun{p}} \Bigg)^k  \sum_{\substack{\Gamma\subset \n\\  \mid \Gamma\mid =k}}  \sum_{\substack{S\subset \n\\ S\supset \Gamma}}
                    \las{S}  \fun{\zeta^{S^c}} \notag\\
                  & \times \det \Big( \big(\mxprobsdummy\circ\mxtrans\big)^{S\setminus \Gamma}_S  \Join \big(\mxprobsreal\circ\mxtrans\big)^\Gamma_S\Big), \label{e3.denominator with only ph customers replace}
\end{align}
and the numerator of \ltwc{i} (see Theorem~\ref{t3.numerator LT workload state i}) becomes
\begin{small}
\begin{align}
s\vup& \ltm{\mxadj}\uv_i = su_i \sum_{k=0}^{N-1} \Bigg( \frac{\fun{q}}{\fun{p}} \Bigg)^k
                        \sum_{\substack{\Gamma\subset \n\setminus\{i\}\\  \mid \Gamma\mid =k; \\ S\subset \n\setminus\{i\}\\ S\supset \Gamma}}
                              \las{S}  \fun{\zeta^{S^c}} \notag\\
                      & \times \det \Big( \big(\mxprobsdummy\circ\mxtrans\big)^{S\setminus \Gamma}_S \Join \big(\mxprobsreal\circ\mxtrans\big)^\Gamma_S  \Big) \notag\\
                     +& s\sum_{\substack{l=1\\ l\neq i}}^N u_l (-1)^{l+i} \sum_{k=1}^{N-1} \Bigg( \frac{\fun{q}}{\fun{p}} \Bigg)^{k}
                         \sum_{\substack{\Gamma\subset \n\setminus\{l,i\}\\  \mid \Gamma\mid =k-1;\\S\subset \n\setminus\{l,i\}\\ S\supset \Gamma;\\ R\subset S\cap T_{li}}}
                        (-1)^{\mid R\mid}  \las{S\cup \{i\}} \notag\\
                      & \times\fun{\zeta^{S^c}}
                        \det \Big( \big(\mxprobsdummy\circ\mxtrans\big)^{S\setminus \Gamma}_{S\cup \{l\}} \Join \big(\mxprobsreal\circ\mxtrans\big)^{\Gamma \cup \{i\}}_{S\cup \{l\}}  \Big)
                       \notag\\
                     +& s\sum_{\substack{l=1\\ l\neq i}}^N u_l (-1)^{l+i} \sum_{k=0}^{N-2}\Bigg( \frac{\fun{q}}{\fun{p}} \Bigg)^k
                         \sum_{\substack{\Gamma\subset \n\setminus\{l,i\}\\  \mid \Gamma\mid =k;\\S\subset \n\setminus\{l,i\}\\ S\supset \Gamma;\\ R\subset S\cap T_{li}}}
                        (-1)^{\mid R\mid} \las{S\cup \{i\}} \notag\\
                      & \times \fun{\zeta^{S^c}}
                        \det \Big( \big(\mxprobsdummy\circ\mxtrans\big)^{(S\setminus \Gamma)\cup \{i\}}_{S\cup \{l\}} \Join \big(\mxprobsreal\circ\mxtrans\big)^{\Gamma}_{S\cup \{l\}}  \Big).
                        \label{e3.numerator with only ph customers replace}
\end{align}
\end{small}

Observe that both the denominator \eqref{e3.denominator with only ph customers replace} and the numerator \eqref{e3.numerator with only ph customers replace} of \ltwc{i}, are rational functions with denominators the polynomial \fun{p} raised to some power. To simplify as much as possible the expression of \ltwc{i}, we multiply \eqref{e3.denominator with only ph customers replace} and \eqref{e3.numerator with only ph customers replace} with $\big( \fun{p} \big)^r$, where $r \in \n$ is the highest possible power of \fun{p} that is involved in the formulas. It is immediately obvious that $r \leq K$. Therefore, we multiply both \eqref{e3.denominator with only ph customers replace} and \eqref{e3.numerator with only ph customers replace} with $\big(\fun{p} \big)^r$

When multiplied with $\big( \fun{p} \big)^r$, the denominator of \ltwc{i} becomes
\begin{align}
  \big( \fun{p} & \big)^r \det \ltm{\mxe} = \notag\\
                  &  \big( \fun{p} \big)^r \sum_{S\subset \n}\las{S}\fun{\zeta^{S^c}}\det \big(\mxprobsdummy\circ\mxtrans\big)^S_S  \notag \\
                 +& \sum_{k=1}^N \big( \fun{q} \big)^k \big( \fun{p} \big)^{r-k}  \sum_{\substack{\Gamma\subset \n\\  \mid \Gamma\mid =k}}  \sum_{\substack{S\subset \n\\ S\supset \Gamma}} \las{S}  \fun{\zeta^{S^c}} \notag \\
                  & \times \det \Big( \big(\mxprobsdummy\circ\mxtrans\big)^{S\setminus \Gamma}_S  \Join \big(\mxprobsreal\circ\mxtrans\big)^\Gamma_S\Big).
                                 \label{e3.simplified denominator with only ph customers replace}
\end{align}

The term $\big( \fun{p} \big)^r \sum_{S\subset \n}\las{S}\fun{\zeta^{S^c}}\det \big(\mxprobsdummy\circ\mxtrans\big)^S_S$ is a polynomial of degree $rM + N$. The coefficient of $s^{rM+N}$ is found when we set $S = \emptyset$, and it is equal to 1. On the other hand, the second term of the right hand side of \eqref{e3.simplified denominator with only ph customers replace} is a polynomial of degree at most $n+(r-1)M+N-1$ (the highest order of $s$ is found when $\left| S\right| = 1$). Since $n\leq M-1$, it is immediately obvious that $\big( \fun{p} \big)^r \det \ltm{\mxe}$ is a polynomial of degree $N+rM$, thus it has exactly $N+rM$ roots. From Theorem~\ref{t3.solution for vector u}, we know that exactly $N-1$ of its roots have positive real part and that zero is also a root. We denote these roots as $\rootp{1}=0$, and \rootp{k}, $k=2,\dots,N$, and we assume them to be simple. We denote the remaining $rM$ roots with negative real part as $-\rootnden{j}$, $j=1,\dots,rM$. Consequently, the denominator of \ltwc{i} is written as
\begin{equation}\label{e3.denominator as a product replace}
  \big( \fun{p} \big)^r \det \ltm{\mxe} = s \prod_{k=2}^{N} (s-\rootp{k}) \prod_{j=1}^{rM} (s+\rootnden{j}).
\end{equation}
Similarly, the numerator of \ltwc{i} becomes
\begin{small}
\begin{align}
\big( \fun{p} & \big)^r s \vup \ltm{\mxadj} \uv_i \notag \\
                     =& su_i \big( \fun{p} \big)^{r}
                        \sum_{\substack{S\subset \n\setminus\{i\}}}
                              \las{S}  \fun{\zeta^{S^c}} \det \big(\mxprobsdummy\circ\mxtrans\big)^{S}_S \notag\\
                     +& su_i \sum_{k=1}^{N-1} \big( \fun{q} \big)^k \big( \fun{p} \big)^{r-k}
                        \sum_{\substack{\Gamma\subset \n\setminus\{i\}\\  \mid \Gamma\mid =k; \\ S\subset \n\setminus\{i\}\\ S\supset \Gamma}}
                              \las{S}  \fun{\zeta^{S^c}} \notag\\
                      & \times \det \Big( \big(\mxprobsdummy\circ\mxtrans\big)^{S\setminus \Gamma}_S \Join \big(\mxprobsreal\circ\mxtrans\big)^\Gamma_S  \Big) \notag\\
                     +& s\sum_{\substack{l=1\\ l\neq i}}^N u_l (-1)^{l+i} \sum_{k=1}^{N-1} \big( \fun{q} \big)^{k} \big( \fun{p} \big)^{r-k} \notag \\
                      & \times  \sum_{\substack{\Gamma\subset \n\setminus\{l,i\}\\  \mid \Gamma\mid =k-1;\\S\subset \n\setminus\{l,i\}\\ S\supset \Gamma;\\ R\subset S\cap T_{li}}}
                        (-1)^{\mid R\mid} \las{S\cup \{i\}} \fun{\zeta^{S^c}}\notag\\
                      & \times \det \Big( \big(\mxprobsdummy\circ\mxtrans\big)^{S\setminus \Gamma}_{S\cup \{l\}} \Join \big(\mxprobsreal\circ\mxtrans\big)^{\Gamma \cup \{i\}}_{S\cup \{l\}}  \Big)
                       \notag\\
                     +& s\sum_{\substack{l=1\\ l\neq i}}^N u_l (-1)^{l+i} \sum_{k=0}^{N-2} \big( \fun{q} \big)^k \big( \fun{p} \big)^{r-k} \notag \\
                      & \times \sum_{\substack{\Gamma\subset \n\setminus\{l,i\}\\  \mid \Gamma\mid =k;\\S\subset \n\setminus\{l,i\}\\ S\supset \Gamma;\\ R\subset S\cap T_{li}}}
                        (-1)^{\mid R\mid} \las{S\cup \{i\}} \fun{\zeta^{S^c}}\notag\\
                      & \times \det \Big( \big(\mxprobsdummy\circ\mxtrans\big)^{(S\setminus \Gamma)\cup \{i\}}_{S\cup \{l\}} \Join \big(\mxprobsreal\circ\mxtrans\big)^{\Gamma}_{S\cup \{l\}}  \Big).
                        \label{e3.simplified numerator with only ph customers replace}
\end{align}
\end{small}

It is easy to verify that $\big( \fun{p} \big)^r s \vup \ltm{\mxadj} \uv_i$ is also a polynomial of degree $rM+N$. The coefficient of $s^{rM+N}$ is equal to $u_i$ and it is determined by the term $su_i \big( \fun{p} \big)^{r} \allowbreak\sum_{S\subset \n\setminus\{i\}} \las{S}  \fun{\zeta^{S^c}}\det \big(\mxprobsdummy\circ\mxtrans\big)^{S}_S  $ for $S=\emptyset$.
We know from Theorem~\ref{t3.solution for vector u}, that the vector \vup is such that the numbers \rootp{k}, $k \in \n$, are also roots of the numerator of \ltwc{i}. We denote the rest $rM$ roots of the numerator as $-\rootnnum{i,j}$, $j=1,\dots,rM$. Therefore, the numerator of \ltwc{i} is written as
\begin{equation}\label{e3.numerator as a product replace}
  \big( \fun{p} \big)^r s \vup \ltm{\mxadj} \uv_i = u_i s \prod_{k=2}^{N} (s-\rootp{k}) \prod_{j=1}^{rM} (s+\rootnnum{i,j}).
\end{equation}
Combining \eqref{e3.denominator as a product replace} and \eqref{e3.numerator as a product replace}, the result is immediate.
\end{proof}

\begin{proof}[of Theorem~\ref{t3.vector u in the mixture model}]
  Since $\ltm[0]{\mxk}$ is an $N \times N$ zero matrix, it is evident that $\rootpmix{1}=0$ is an eigenvalue of the matrix $\ltm{\mxhmix} +s \im -\mxrates$ (see Eq.~\eqref{e3.definition perturbed matrix e}). According to Corollary~\ref{c3.roots of the perturbed matrix}, the numbers \rootpmix{i}, $i=2,\dots,N$, are also simple eigenvalues of this matrix. Thus, according to Theorem~\ref{t3.solution for vector u}, there are no other roots of the equation $\det \big(\ltm{\mxe} + \epsilon \ltm{\mxk}\big)=0$ with non\-/negative real part besides the values \rootpmix{i}, $i \in \n$. \ \\

  \noindent For the second part of proof we have the following. Using Theorem~\ref{t3.eigenvectors perturbed system}, we can evaluate $N-1$ column vectors $\vect{w}_{\epsilon,i}$ such that
  \begin{equation*}
    \big( \ltm[\rootpmix{i}]{\mxhmix} +\rootpmix{i} \im -\mxrates \big) \vect{w}_{\epsilon,i} = 0, \qquad i=2,\dots,N.
  \end{equation*}
  Since $\rootpmix{i} \neq 0$, $i=2,\dots,N$, post-multiplying equation~\eqref{e3.eq1 for trasform vector mixture model} with $s=\rootpmix{i}$ by $\vect{w}_{\epsilon,i}$, we obtain
  \begin{equation*}
    \vupmix \vect{w}_{\epsilon,i} =0, \qquad i=2,\dots,N.
  \end{equation*}
  To derive the remaining equation, we take the derivative of equation \eqref{e3.eq1 for trasform vector mixture model} with respect to $s$, yielding
  \begin{equation*}
    \vltwmix \big( \ltm{\dm[1]{\mxhmix}} + \im \big) + \vltwmixder\big( \ltm{\mxhmix} +s \im -\mxrates \big) = \vupmix.
  \end{equation*}
  Setting $s=0$ we get
  \begin{equation*}
    \vltwmix[0] \big( \ltm[0]{\dm[1]{\mxhmix}} + \im \big) + \vltwmixder[0]\big( \mxtrans - \im \big)\mxrates = \vupmix.
  \end{equation*}
  Post-multiplying by $\mxrates^{-1}\uv$ gives
  \begin{align*}
    &\vltwmix[0] \big( \ltm[0]{\dm[1]{\mxhmix}} + \im \big)\mxrates^{-1}\uv + \vltwmixder[0]\big( \mxtrans - \im \big)\mxrates \mxrates^{-1}\uv \\ &= \vupmix \mxrates^{-1}\uv.
  \end{align*}
  Finally, using $(\mxtrans - \im) \uv=0$, $\ltm[0]{\dm[1]{\mxhmix}} = - \mxmeans \mxrates +\epsilon \big(\mean[p]-\mean[h]\big)\mxprobsreal \circ\mxtrans \mxrates$ and $\vltwmix[0]=\boldsymbol\pi$ (where the latter follows from \eqref{e3.eq1 for trasform vector mixture model} with $s=0$ and the normalization equation \eqref{e3.eq2 for trasform vector mixture model}), the above can be simplified to
  \begin{equation*}
    \boldsymbol\pi\left(\mxrates^{-1} - \mxmeans\right) \uv + \epsilon (\mean[p]-\mean[h])\boldsymbol\pi \mxprobsreal \circ\mxtrans \uv = \vupmix \mxrates^{-1}\uv.
  \end{equation*}
  The uniqueness of the solution follows from the general theory of Markov chains that under the condition of stability, there is a unique stationary distribution and thus also a unique solution \vltwmix to the equations \eqref{e3.eq1 for trasform vector mixture model} and \eqref{e3.eq2 for trasform vector mixture model}. This completes the proof.
\end{proof}

\begin{proof}[of Proposition~\ref{p3.series expansion of the corrected replace transform vector}]
Recall that $r$ is the maximum power of \fun{p} that appears in the formulas. Therefore, to use perturbation analysis, we multiply both $\det \ltm{\mxemix}$ and $s\vupmix \ltm{\mxadjmix}\uv$ with $\big(\fun{p}\big)^r$. So, if we set
\begin{gather}
\fun{\xi_{rM+N-1}}   = \sum_{k=1}^N k \big( \fun{q} \big)^{k-1} \big( \fun{p} \big)^{r-k+1} \sum_{\substack{\Gamma\subset \n\\  \mid \Gamma\mid =k; \\ S\subset \n\\ S\supset \Gamma}}
                            \las{S}\notag \\
                        \times \fun{\zeta^{S^c}} \det \Big( \big(\mxprobsdummy\circ\mxtrans\big)^{S\setminus \Gamma}_S  \Join \big(\mxprobsreal\circ\mxtrans\big)^\Gamma_S\Big),\label{e3.function xi}
\end{gather}
then,
\begin{align}
    \big(\fun{p}\big)^r\det \ltm{\mxemix}  &= \big(\fun{p}\big)^r\det \ltm{\mxe} \notag \\
                                            +& \epsilon s \big(\mean[p]\ltptes - \mean[h]\lthtes\big) \fun{\xi_{rM+N-1}} \notag \\
                                            +& O(\epsilon^2).\notag
\end{align}
Note that the polynomial \fun{\xi_{rM+N-1}} is of degree at most $rM+N-1$, and the coefficient of $s^{rM+N-1}$ is equal to $\gamma=\sum_{i=1}^N \lambda_i
\det \big(\mxprobsreal\circ\mxtrans\big)^{\{i\}}_{\{i\}} = \sum_{i=1}^N \lambda_i q_{ii}^{(2)} p_{ii}$. Theorem~\ref{t3.vector u in the mixture model} guarantees that the polynomial $\big(\fun{p}\big)^r\times\allowbreak \det \ltm{\mxemix}$ has exactly $N-1$ roots with positive real part and it also has $\rootpmix{1}=0$. The roots with positive real part are of the form $\rootpmix{k} = \rootp{k} - \epsilon \delta_k +O(\epsilon^2)$, $k=2,\dots,N$, where
\begin{small}
\begin{align}
  \delta_k  &= \frac{\big(\mean[p]\ltptes[\rootp{k}] - \mean[h]\lthtes[\rootp{k}]\big) \fun[\rootp{k}]{\xi_{rM+N-1}}}
              {\prod_{\substack{l=2\\ l\neq k}}^{N} (\rootp{k}-\rootp{l}) \prod_{j=1}^{rM} (\rootp{k}+\rootnden{j})} \notag \\
            &= \frac{\big(\mean[p]\ltptes[\rootp{k}] - \mean[h]\lthtes[\rootp{k}]\big) \fun[\rootp{k}]{\xi_{rM+N-1}} \ltwc[\rootp{k}]{i}}
              {u_i \prod_{\substack{l=2\\ l\neq k}}^{N} (\rootp{k}-\rootp{l}) \prod_{j=1}^{rM} (\rootp{k}+\rootnnum{i,j})}.\label{e3.definition 1 for dk}
\end{align}
\end{small}
Thus, if we set
\begin{align}
  \fes =& \frac{\big(\mean[p]\ltptes - \mean[h]\lthtes\big) \fun{\xi_{rM+N-1}}\ltwc{i}}{u_i \prod_{k=2}^{N} (s-\rootp{k}) \prod_{j=1}^{rM} (s+\rootnnum{i,j})} \notag \\
        &- \sum_{k=2}^N\frac{\delta_k}{s - \rootp{k}},\label{e3.function d mixture model}
\end{align}
the denominator of \ltwcmix{,i} multiplied by $ \big(\fun{p}\big)^r$ can be written as
\begin{align}
      \big( \fun{p} & \big)^r\det \ltm{\mxemix} \notag \\
     =& s \prod_{j=1}^{rM} (s+\rootnden{j})\prod_{k=2}^{N} (s-\rootp{k} + \epsilon \delta_k + O(\epsilon^2)) \notag \\
      & \times  \big( 1 + \epsilon \fes + O(\epsilon^2) \big).\label{e3.denominator mixture}
\end{align}
Note that the function \fes is well defined in the positive half plane due to the definition \eqref{e3.definition 1 for dk} of $\delta_k$, $k=2,\dots,N$. Similarly, if we set
\begin{align}
  & \fun{\xi_{i,l,rM+N-2}} \notag \\
 =& \indfun[\{l=i\}] \sum_{k=1}^{N-1} k \big( \fun{q} \big)^{k-1} \big( \fun{p} \big)^{r-k+1}
                        \sum_{\substack{\Gamma\subset \n\setminus\{i\}\\  \mid \Gamma\mid =k; \\ S\subset \n\setminus\{i\}\\ S\supset \Gamma}}   \las{S}  \notag \\
  & \times \fun{\zeta^{S^c}}  \det \Big( \big(\mxprobsdummy\circ\mxtrans\big)^{S\setminus \Gamma}_S \Join \big(\mxprobsreal\circ\mxtrans\big)^\Gamma_S  \Big) \notag\\
 +& \indfun[\{l\neq i\}] \Bigg[(-1)^{l+i} \sum_{k=1}^{N-1} k \big( \fun{q} \big)^{k-1} \big( \fun{p} \big)^{r-k+1} \notag \\
  & \times  \sum_{\substack{\Gamma\subset \n\setminus\{l,i\}\\  \mid \Gamma\mid =k-1}}
            \sum_{\substack{S\subset \n\setminus\{l,i\}\\ S\supset \Gamma;\\ R\subset S\cap T_{li}}}
                        (-1)^{\mid R\mid} \las{S\cup \{i\}} \fun{\zeta^{S^c}} \notag \\
  & \times \det \Big( \big(\mxprobsdummy\circ\mxtrans\big)^{S\setminus \Gamma}_{S\cup \{l\}} \Join \big(\mxprobsreal\circ\mxtrans\big)^{\Gamma \cup \{i\}}_{S\cup \{l\}}  \Big) \notag \\
  &+ (-1)^{l+i} \sum_{k=1}^{N-2} k \big( \fun{q} \big)^{k-1} \big( \fun{p} \big)^{r-k+1} \notag \\
  & \times  \sum_{\substack{\Gamma\subset \n\setminus\{l,i\}\\  \mid \Gamma\mid =k}}
            \sum_{\substack{S\subset \n\setminus\{l,i\}\\ S\supset \Gamma;\\ R\subset S\cap T_{li}}}
                        (-1)^{\mid R\mid} \las{S\cup \{i\}} \fun{\zeta^{S^c}} \notag \\
  & \times \det \Big( \big(\mxprobsdummy\circ\mxtrans\big)^{(S\setminus \Gamma)\cup \{i\}}_{S\cup \{l\}} \Join \big(\mxprobsreal\circ\mxtrans\big)^{\Gamma}_{S\cup \{l\}}  \Big)\Bigg]
                                    ,\label{e3.function xi by state}
\end{align}
and
\begin{align}
 & \fun{\xi'_{i,l, rM+N-1}} \notag \\
=& \indfun[\{l=i\}]\Bigg[ \big( \fun{p} \big)^{r}\sum_{\substack{S\subset \n\setminus\{i\}}} \las{S}  \fun{\zeta^{S^c}} \det \big(\mxprobsdummy\circ\mxtrans\big)^{S}_S  \notag \\
 &+ \sum_{k=1}^{N-1} \big( \fun{q} \big)^k \big( \fun{p} \big)^{r-k} \sum_{\substack{\Gamma\subset \n\setminus\{i\}\\  \mid \Gamma\mid =k}}
 \sum_{\substack{S\subset \n\setminus\{i\}\\ S\supset \Gamma}}
      \las{S}  \fun{\zeta^{S^c}} \notag\\
 & \times \det \Big( \big(\mxprobsdummy\circ\mxtrans\big)^{S\setminus \Gamma}_S \Join \big(\mxprobsreal\circ\mxtrans\big)^\Gamma_S  \Big) \Bigg] \notag\\
+& \indfun[\{l\neq i\}] \Bigg[ (-1)^{l+i} \sum_{k=1}^{N-1} \big( \fun{q} \big)^{k} \big( \fun{p} \big)^{r-k} \notag \\
 & \times  \sum_{\substack{\Gamma\subset \n\setminus\{l,i\}\\  \mid \Gamma\mid =k-1}} \sum_{\substack{S\subset \n\setminus\{l,i\}\\ S\supset \Gamma;\\ R\subset S\cap T_{li}}}
                        (-1)^{\mid R\mid} \las{S\cup \{i\}} \fun{\zeta^{S^c}}\notag\\
 & \times \det \Big( \big(\mxprobsdummy\circ\mxtrans\big)^{S\setminus \Gamma}_{S\cup \{l\}} \Join \big(\mxprobsreal\circ\mxtrans\big)^{\Gamma \cup \{i\}}_{S\cup \{l\}}  \Big) \notag\\
+& (-1)^{l+i} \sum_{k=0}^{N-2} \big( \fun{q} \big)^k \big( \fun{p} \big)^{r-k} \notag \\
 & \times \sum_{\substack{\Gamma\subset \n\setminus\{l,i\}\\  \mid \Gamma\mid =k}} \sum_{\substack{S\subset \n\setminus\{l,i\}\\ S\supset \Gamma;\\ R\subset S\cap T_{li}}}
                        (-1)^{\mid R\mid} \las{S\cup \{i\}} \fun{\zeta^{S^c}}\notag\\
 & \times \det \Big( \big(\mxprobsdummy\circ\mxtrans\big)^{(S\setminus \Gamma)\cup \{i\}}_{S\cup \{l\}} \Join \big(\mxprobsreal\circ\mxtrans\big)^{\Gamma}_{S\cup \{l\}}  \Big)
                                \Bigg],\label{e3.function xi prime by state}
\end{align}
then
\begin{align*}
  \big(\fun{p_m} &\big)^r s \vupmix \ltm{\mxadjmix}\uv_i  = \big(\fun{p_m}\big)^r s \vup \ltm{\mxadj}\uv_i \notag \\
                           +& \epsilon s \Bigg[ \sum_{l=1}^N z_l \fun{\xi'_{i,l, rM+N-1}} \notag \\
                            &+ s \big(\mean[p]\ltptes - \mean[h]\lthtes\big) \sum_{l=1}^N u_l \fun{\xi_{i,l,rM+N-2}}\Bigg]\\
                           +& O(\epsilon^2)  . \notag
\end{align*}
Note that the polynomial $\sum_{l=1}^N z_l \fun{\xi'_{i,l, rM+N-1}}$ is of degree $rM+N-1$, and the coefficient of $s^{rM+N-1}$ is $z_i$. Analogously, the polynomial $s \sum_{l=1}^N u_l \fun{\xi_{i,l,rM+N-2}}$ is of degree at most $rM+N-1$, and the coefficient of $s^{rM+N-1}$ is equal to $\beta = \sum_{l=1}^N u_l \Big( \indfun[\{l=i\}] \sum_{\substack{j=1 \\ j \neq i}}^N  \lambda_j q^{(2)}_{jj}p_{jj} \allowbreak+ \indfun[\{l \neq i\}] (-1)^{l+i} \lambda_i q_{li}^{(2)} p_{li} \Big)$. The first part is for $S=\Gamma=\{j\}$, and the second part for $S=\Gamma=\emptyset$. Theorem~\ref{t3.vector u in the mixture model} guarantees that the roots \rootpmix{k}, $k \in \n$, are also roots of the numerator of \ltwcmix{,i}. Therefore, applying perturbation analysis to $\big(\fun{p_m}\big)^r s \vupmix \ltm{\mxadjmix}\uv_i$ results in an equivalent definition for each $\delta_k$, $k=2,\dots,N$, as
\begin{align}
  \delta_k  =& \frac{\rootp{k} \big(\mean[p]\ltptes[\rootp{k}] - \mean[h]\lthtes[\rootp{k}]\big) \sum_{l=1}^N u_l \fun[\rootp{k}]{\xi_{i,l,rM+N-2}}}
              {u_i \prod_{\substack{l=2\\ l\neq k}}^{N} (\rootp{k}-\rootp{l}) \prod_{j=1}^{rM} (\rootp{k}+\rootnnum{i,j})} \notag \\
             &+ \frac{\sum_{l=1}^N z_l \fun[\rootp{k}]{\xi'_{i,l, rM+N-1}} }
              {u_i \prod_{\substack{l=2\\ l\neq k}}^{N} (\rootp{k}-\rootp{l}) \prod_{j=1}^{rM} (\rootp{k}+\rootnnum{i,j})}  .\label{e3.definition 2 for dk}
\end{align}
Now, if we set
\begin{align}
  \fef =& \frac{s \big(\mean[p]\ltptes - \mean[h]\lthtes\big) \sum_{l=1}^N u_l \fun{\xi_{i,l,rM+N-2}}}
                  {u_i \prod_{k=2}^{N} (s-\rootp{k}) \prod_{j=1}^{rM} (s+\rootnnum{i,j})} \notag \\
       +&\frac{\sum_{l=1}^N z_l \fun{\xi'_{i,l, rM+N-1}} }
                  {u_i \prod_{k=2}^{N} (s-\rootp{k}) \prod_{j=1}^{rM} (s+\rootnnum{i,j})} \notag \\
       -& \sum_{k=2}^N\frac{\delta_k}{s - \rootp{k}},\label{e3.function n mixture model}
\end{align}
the numerator of \ltwcmix{,i} multiplied by $ \big(\fun{p_m}\big)^r$ can be written as
\begin{align}
    \big(\fun{p_m} & \big)^r s \vupmix \ltm{\mxadjmix}\uv_i \notag \\
     =&  u_i s \prod_{j=1}^{rM} (s+\rootnnum{i,j}) \prod_{k=2}^{N} \big(s-\rootp{k} + \epsilon \delta_k + O(\epsilon^2)\big) \notag \\
      & \times \big( 1 + \epsilon \fef + O(\epsilon^2) \big) .\label{e3.numerator mixture}
\end{align}
Note that the function \fef is well defined in the positive half plane due to the definition \eqref{e3.definition 2 for dk} of $\delta_k$, $k=2,\dots,N$. Combining \eqref{e3.denominator mixture} and \eqref{e3.numerator mixture}, we obtain
\begin{small}
\begin{align}
  \ltwcmix{,i} =& \frac{u_i \prod_{j=1}^{rM} (s+\rootnnum{i,j})}{\prod_{j=1}^{rM} (s+\rootnden{j})}\cdot
                  \frac{1 + \epsilon \fef + O(\epsilon^2)}{1 + \epsilon \fes + O(\epsilon^2)} \notag \\
               =& \ltwc{i}\big( 1 + \epsilon \fef + O(\epsilon^2) \big)\big( 1 - \epsilon \fes + O(\epsilon^2) \big)  \notag \\
               =& \ltwc{i} + \epsilon \ltwc{i} \big(\fef - \fes\big) + O(\epsilon^2) \notag \\
               =& \ltwc{i} + \epsilon \frac1{u_i}\ltwc{i} \Bigg( \frac{\sum_{l=1}^N z_l \fun{\xi'_{i,l, rM+N-1}}}{\prod_{k=2}^{N} (s-\rootp{k}) \prod_{j=1}^{rM} (s+\rootnnum{i,j})} \notag \\
                &+ \big(\mean[p]\ltptes - \mean[h]\lthtes\big) \frac{s \sum_{l=1}^N u_l \fun{\xi_{i,l,rM+N-2}}}{\prod_{k=2}^{N} (s-\rootp{k}) \prod_{j=1}^{rM} (s+\rootnnum{i,j})} \notag \\
                &- \big(\mean[p]\ltptes - \mean[h]\lthtes\big)\ltwc{i} \notag \\
                & \times \frac{\fun{\xi_{rM+N-1}}}{\prod_{k=2}^{N} (s-\rootp{k}) \prod_{j=1}^{rM} (s+\rootnnum{i,j})} \Bigg) + O(\epsilon^2) \notag \\
               =& \ltwc{i} + \epsilon \frac1{u_i}\ltwc{i} \Bigg[ \Bigg(z_i + \sum_{k=2}^N \frac{\alpha_{i,k}}{s-\rootp{k}} + \sum_{j=1}^{rM} \frac{\alpha'_{i,j}\cdot \rootnnum{i,j}}{s+\rootnnum{i,j}}\Bigg) \notag \\
                &+ \big(\mean[p]\ltptes - \mean[h]\lthtes\big)  \Bigg(\beta + \sum_{k=2}^N \frac{\beta_{i,k}}{s-\rootp{k}} \notag \\
                &+ \sum_{j=1}^{rM} \frac{\beta'_{i,j} \cdot \rootnnum{i,j}}{s+\rootnnum{i,j}}\Bigg) \notag \\
                &- \big(\mean[p]\ltptes - \mean[h]\lthtes\big)\ltwc{i}  \Bigg(\gamma + \sum_{k=2}^N \frac{\gamma_{i,k}}{s-\rootp{k}} \notag \\
                & + \sum_{j=1}^{rM} \frac{\gamma'_{i,j}\cdot \rootnnum{i,j}}{s+\rootnnum{i,j}} \Bigg) \Bigg]  + O(\epsilon^2),
\end{align}
\end{small}
where the last equality comes from simple fraction decomposition under the assumption that the roots $-\rootnnum{i,j}$, $j=1,\dots,rM$, are simple. The coefficients $\alpha_k, \beta_k, \gamma_k$, $k=2,\dots,N$, and $\alpha'_j, \beta'_j, \gamma'_j$, $j=1,\dots,rM$, are as follows
\begin{align}
  \alpha_{i,k}  &= \frac{\sum_{l=1}^N z_l \fun[\rootp{k}]{\xi'_{i,l, rM+N-1}}}
                      {\prod_{\substack{l=2\\ l\neq k}}^{N} (\rootp{k}-\rootp{l}) \prod_{j=1}^{rM} (\rootp{k}+\rootnnum{i,j})}, \label{e3.coeff 1} \\
  \beta_{i,k}   &= \frac{\rootp{k} \sum_{l=1}^N u_l \fun[\rootp{k}]{\xi_{i,l,rM+N-2}}}
                      {\prod_{\substack{l=2\\ l\neq k}}^{N} (\rootp{k}-\rootp{l}) \prod_{j=1}^{rM} (\rootp{k}+\rootnnum{i,j})}, \label{e3.coeff 2}\\
  \gamma_{i,k}  &= \frac{\fun[\rootp{k}]{\xi_{rM+N-1}}}
                      {\prod_{\substack{l=2\\ l\neq k}}^{N} (\rootp{k}-\rootp{l}) \prod_{j=1}^{rM} (\rootp{k}+\rootnnum{i,j})}, \label{e3.coeff 3}\\
  \alpha'_{i,j} &= \frac{\sum_{l=1}^N z_l \fun[-\rootnnum{i,j}]{\xi'_{i,l, rM+N-1}}}
                      {\rootnnum{i,j}\prod_{k=2}^{N} (-\rootnnum{i,j}-\rootp{k}) \prod_{\substack{l=1\\ l\neq j}}^{rM}(-\rootnnum{i,j}+\rootnnum{i,l})}, \label{e3.coeff 4}\\
  \beta'_{i,j}  &= \frac{-\sum_{l=1}^N u_l \fun[-\rootnnum{i,j}]{\xi_{i,l,rM+N-2}}}
                      {\prod_{k=2}^{N} (-\rootnnum{i,j}-\rootp{k}) \prod_{\substack{l=1\\ l\neq j}}^{rM}(-\rootnnum{i,j}+\rootnnum{i,l})}, \label{e3.coeff 5}\\
  \gamma'_{i,j} &= \frac{\fun[-\rootnnum{i,j}]{\xi_{rM+N-1}}}
                      {\rootnnum{i,j}\prod_{k=2}^{N} (-\rootnnum{i,j}-\rootp{k}) \prod_{\substack{l=1\\ l\neq j}}^{rM}(-\rootnnum{i,j}+\rootnnum{i,l})}. \label{e3.coeff 6}
\end{align}
The above results hold when all roots $-\rootnnum{i,j}$, $j=1,\dots,rM$, are simple. Suppose now that only $\sigma$ of the roots are distinct and that the multiplicity of root $-\rootnnum{i,j}$, $j=1,\dots,\sigma$, is $r_{i,j}$, such that $\sum_{j=1}^\sigma r_{i,j} = rM$. In this case,
\begin{align}
  \ltwcmix{,i}  &= \ltwc{i} + \epsilon \frac1{u_i}\ltwc{i} \Bigg[ \Bigg(z_i + \sum_{k=2}^N \frac{\alpha_{i,k}}{s-\rootp{k}} \notag \\
                &+ \sum_{j=1}^{\sigma}\sum_{l=1}^{r_{i,j}} \frac{\alpha''_{i,j,l} \cdot (\rootnnum{i,j})^{r_{i,j}-l+1}}{(s+\rootnnum{i,j})^{r_{i,j}-l+1}}\Bigg) \notag \\
                &+ \big(\mean[p]\ltptes - \mean[h]\lthtes\big)\Bigg( \beta
                 + \sum_{k=2}^N \frac{\beta_{i,k}}{s-\rootp{k}} \notag \\
                &+ \sum_{j=1}^{\sigma}\sum_{l=1}^{r_{i,j}} \frac{\beta''_{i,j,l} \cdot (\rootnnum{i,j})^{r_{i,j}-l+1}}{(s+\rootnnum{i,j})^{r_{i,j}-l+1}}\Bigg) \notag\\
                &- \big(\mean[p]\ltptes - \mean[h]\lthtes\big)\ltwc{i} \Bigg( \gamma
                 + \sum_{k=2}^N \frac{\gamma_{i,k}}{s-\rootp{k}} \notag \\
                &+ \sum_{j=1}^{\sigma}\sum_{l=1}^{r_{i,j}} \frac{\gamma''_{i,j,l} \cdot (\rootnnum{i,j})^{r_{i,j}-l+1}}{(s+\rootnnum{i,j})^{r_{i,j}-l+1}} \Bigg) \Bigg] %\notag \\
                 + O(\epsilon^2),\label{e3.unordered laplace transform of the mixture workload in state i}
\end{align}
where $\alpha_{i,k}$, $\beta_{i,k}$ and $\gamma_{i,k}$, $k=2,\dots,N$, are defined through \eqref{e3.coeff 1}--\eqref{e3.coeff 3}. For each $j=1,\dots,\sigma$, the coefficients $\alpha''_{i,j,p}$, $p=1,\dots,r_{i,j}$, are the unique solution to the following linear system of $r_{i,j}$ equations
\begin{small}
\begin{align}
  \frac{d}{ds^n}&\left.\Bigg[\sum_{l=1}^N z_l \fun{\xi'_{i,l, rM+N-1}}\Bigg]\right|_{s=-\rootnnum{i,j}} \notag \\
              =& \frac{d}{ds^n}\Bigg[ \prod_{k=2}^{N} (s-\rootp{k}) \prod_{\substack{l=1\\ l\neq j}}^{\sigma}(s+\rootnnum{i,l})^{r_{i,l}} \sum_{p=1}^{r_{i,j}} \alpha''_{i,j,p} \notag \\
               & \times \left. (\rootnnum{i,j})^{r_{i,j}-p+1}(s+\rootnnum{i,j})^{p-1} \Bigg]\right|_{s=-\rootnnum{i,j}}, \label{e3.coeff 7}
\end{align}
\end{small}
\noindent for $n=0,\dots,r_{i,j}$. Similarly, for each $j=1,\dots,\sigma$, the coefficients $\beta''_{i,j,p}$ and $\gamma''_{i,j,p}$, $p=1,\dots,r_{i,j}$, are the respective unique solutions to the following two linear system of $r_{i,j}$ equations
\begin{small}
\begin{align}
  \frac{d}{ds^n}&\left.\Bigg[s \sum_{l=1}^N u_l \fun{\xi_{i,l,rM+N-2}}\Bigg]\right|_{s=-\rootnnum{i,j}} \notag \\
               =& \frac{d}{ds^n}\Bigg[ \prod_{k=2}^{N} (s-\rootp{k}) \prod_{\substack{l=1\\ l\neq j}}^{\sigma}(s+\rootnnum{i,l})^{r_{i,l}}  \sum_{p=1}^{r_{i,j}} \beta''_{i,j,p} \notag \\
                & \times \left.(\rootnnum{i,j})^{r_{i,j}-p+1}(s+\rootnnum{i,j})^{p-1} \Bigg]\right|_{s=-\rootnnum{i,j}}, \label{e3.coeff 8} \\
  \frac{d}{ds^n}&\left.\Bigg[\fun{\xi_{rM+N-1}}\Bigg]\right|_{s=-\rootnnum{i,j}} \notag \\
               =&\frac{d}{ds^n}\Bigg[ \prod_{k=2}^{N} (s-\rootp{k}) \prod_{\substack{l=1\\ l\neq j}}^{\sigma}(s+\rootnnum{i,l})^{r_{i,l}} \sum_{p=1}^{r_{i,j}} \gamma''_{i,j,p} \notag \\
                & \times \left.(\rootnnum{i,j})^{r_{i,j}-p+1}(s+\rootnnum{i,j})^{p-1} \Bigg]\right|_{s=-\rootnnum{i,j}}, \label{e3.coeff 9}
\end{align}
\end{small}
\noindent for $n=0,\dots,r_{i,j}$.
\end{proof}

\begin{proof}[of Theorem~\ref{t3.series expansion of mixture workload state i}]
Here, we follow the notation we introduced in Proposition~\ref{p3.series expansion of the corrected replace transform vector}. We denote by \fun{\lt{\theta_i}} the correction term (the coefficient of $\epsilon$) in the expression of \ltwcmix{i}. In order to apply Laplace inversion to \fun{\lt{\theta_i}}, we first reorder the involved terms (see Eq.~\eqref{e3.unordered laplace transform of the mixture workload in state i}) as
\begin{align}
 \fun{\lt{\theta_i}} =& \frac1{u_i}\ltwc{i} \Bigg[ \Bigg(z_i + \beta \big(\mean[p]\ltptes - \mean[h]\lthtes\big) \notag \\
              &- \gamma \big(\mean[p]\ltptes - \mean[h]\lthtes\big)\ltwc{i}\Bigg) \notag \\
             +&\sum_{k=2}^N \frac1{s-\rootp{k}} \Bigg( \alpha_{i,k} + \beta_{i,k}\big(\mean[p]\ltptes - \mean[h]\lthtes\big) \notag \\
              &- \gamma_{i,k} \big(\mean[p]\ltptes - \mean[h]\lthtes\big)\ltwc{i} \Bigg) \notag \\
             +& \sum_{j=1}^{\sigma}\sum_{l=1}^{r_{i,j}} \frac1{(s+\rootnnum{i,j})^{r_{i,j}-l+1}} \Bigg( \alpha''_{i,j,l} + \beta''_{i,j,l} \notag \\
              & \times \big(\mean[p]\ltptes - \mean[h]\lthtes\big)  - \gamma''_{i,j,l} \big(\mean[p]\ltptes \notag \\
              &- \mean[h]\lthtes\big)\ltwc{i} \Bigg) \Bigg] %\notag \\
              . \label{e3.LST power series with corrected replace}
\end{align}
From the above formula it is evident that only the terms in the middle bracket cannot be inverted directly as they are, because of the singularities they seem to have in the positive half plane. Thus, we treat them separately in the next lines. From the two equivalent definitions \eqref{e3.definition 1 for dk} and \eqref{e3.definition 2 for dk} of the perturbation terms $\delta_k$, $k=2,\dots,N$, and the relations \eqref{e3.coeff 1}--\eqref{e3.coeff 3}
we obtain that
 \begin{small}
 \begin{gather*}
   \alpha_{i,k}\ltwc[\rootp{k}]{i} + \beta_{i,k}\big(\mean[p]\ltptes[\rootp{k}] - \mean[h]\lthtes[\rootp{k}]\big)\ltwc[\rootp{k}]{i} \\
                    - \gamma_{i,k} \big(\mean[p]\ltptes[\rootp{k}] - \mean[h]\lthtes[\rootp{k}]\big)\big(\ltwc[\rootp{k}]{i}\big)^2 = 0, \quad k=2,\dots,N.
 \end{gather*}
 \end{small}
The above equations are equivalent to
  \begin{small}
  \begin{align}
   0=& \alpha_{i,k} \int_{x=0}^{\infty} e^{-\rootp{k} x} d\pr(\workload[i]\leq x) + \beta_{i,k} \notag \\
    & \times \Big( \mean[p] \int_{x=0}^{\infty} e^{-\rootp{k} x} d\pr(\workload[i] + \epts \leq x) \notag \\
    &- \mean[h] \int_{x=0}^{\infty} e^{-\rootp{k} x} d\pr(\workload[i] +\ehts \leq x) \Big) \notag \\
    &- \gamma_{i,k} \Big( \mean[p] \int_{x=0}^{\infty} e^{-\rootp{k} x} d\pr(\workload[i] + \workload[i]' + \epts \leq x) \notag \\
    &- \mean[h] \int_{x=0}^{\infty} e^{-\rootp{k} x} d\pr(\workload[i] + \workload[i]' + \ehts \leq x) \Big),  \label{e3.integral relations}
  \end{align}
  \end{small}
  $k=2,\dots,N$. We first show that
  \begin{small}
  \begin{align}
   \mathcal{L}^{-1}&\Bigg( \sum_{k=2}^N \frac1{s-\rootp{k}} \Big( \alpha_{i,k}\ltwc{i} + \beta_{i,k}\big(\mean[p]\ltptes - \mean[h]\lthtes\big)\ltwc{i} \notag \\
                   &- \gamma_{i,k} \big(\mean[p]\ltptes - \mean[h]\lthtes\big)\big(\ltwc{i}\big)^2 \Big) \Bigg) \notag \\
   =&\sum_{k=2}^N \Bigg[ \gamma_{i,k} \Big(\mean[p] \int_{y=x}^{\infty} e^{\rootp{k} (x-y)} d\pr(\workload[i] + \workload[i]' + \epts \leq y) \notag \\
                   &- \mean[h] \int_{y=x}^{\infty} e^{\rootp{k} (x-y)} d\pr(\workload[i] + \workload[i]' + \ehts \leq y) \Big)\notag \\
                   &- \beta_{i,k} \Big( \mean[p] \int_{y=x}^{\infty} e^{\rootp{k} (x-y)} d\pr(\workload[i] + \epts \leq y) \notag \\
                   &-\mean[h] \int_{y=x}^{\infty} e^{\rootp{k} (x-y)} d\pr(\workload[i] +\ehts \leq y) \Bigg) \notag \\
                   &-  \alpha_{i,k} \int_{y=x}^{\infty} e^{\rootp{k} (x-y)} d\pr(\workload[i]\leq y) \Bigg]. \label{e3.laplace inversion to the middle bracket terms}
  \end{align}
  \end{small}
  Since Laplace transforms turn convolutions of functions into their product, using the property $\int_{y=0}^{\infty}\fun[y]{f} dy = \int_{y=0}^{x} \fun[y]{f} dy + \int_{y=x}^{\infty}\fun[y]{f} dy$ and the relations \eqref{e3.integral relations} we obtain
  \begin{small}
  \begin{align*}
    &\mathcal{L}\left\{\sum_{k=2}^N \Bigg[ \gamma_{i,k} \Bigg(\mean[p] \int_{y=x}^{\infty} e^{\rootp{k} (x-y)} d\pr(\workload[i] + \workload[i]' + \epts \leq y)\right. \notag \\
                                   &- \mean[h] \int_{y=x}^{\infty} e^{\rootp{k} (x-y)} d\pr(\workload[i] + \workload[i]' + \ehts \leq y) \Bigg) \notag \\
                                   &-  \beta_{i,k} \Bigg( \mean[p] \int_{y=x}^{\infty} e^{\rootp{k} (x-y)} d\pr(\workload[i] + \epts \leq y) \notag \\
                                   &- \mean[h] \int_{y=x}^{\infty} e^{\rootp{k} (x-y)} d\pr(\workload[i] +\ehts \leq y) \Bigg) \notag \\
                                   &-  \left. \alpha_{i,k} \int_{y=x}^{\infty} e^{\rootp{k} (x-y)} d\pr(\workload[i]\leq y) \Bigg] \right\} \notag\\
    =&\mathcal{L}\left\{\sum_{k=2}^N \Bigg[ -\gamma_{i,k} \Bigg( \mean[p] \int_{y=0}^{x} e^{\rootp{k} (x-y)} d\pr(\workload[i] + \workload[i]' + \epts \leq y) \right. \notag \\
                                   &- \mean[h] \int_{y=0}^{x} e^{\rootp{k} (x-y)} d\pr(\workload[i] + \workload[i]' + \ehts \leq y)  \Bigg) \notag \\
                                   &+ \beta_{i,k} \Bigg( \mean[p] \int_{y=0}^{x} e^{\rootp{k} (x-y)} d\pr(\workload[i] + \epts \leq y) \notag \\
                                   &- \mean[h] \int_{y=0}^{x} e^{\rootp{k} (x-y)} d\pr(\workload[i] +\ehts \leq y) \Bigg) \Bigg] \notag \\
                                   &+\left. \alpha_{i,k} \int_{y=0}^{x} e^{\rootp{k} (x-y)} d\pr(\workload[i]\leq y)\right\} \notag \\
    =&\sum_{k=2}^N \frac1{s-\rootp{k}} \Big( \alpha_{i,k}\ltwc{i} + \beta_{i,k}\big(\mean[p]\ltptes - \mean[h]\lthtes\big)\ltwc{i} \notag \\
                                   &- \gamma_{i,k} \big(\mean[p]\ltptes - \mean[h]\lthtes\big)\big(\ltwc{i}\big)^2 \Big) \Bigg),
  \end{align*}
  \end{small}
  which proves \eqref{e3.laplace inversion to the middle bracket terms}.

To find the tail probabilities that correspond to the terms in the middle bracket of \eqref{e3.LST power series with corrected replace}, we integrate the inverted Laplace transform in Eq.~\eqref{e3.laplace inversion to the middle bracket terms} from $t$ to $\infty$, and we obtain
\begin{small}
\begin{align}
   \sum_{k=2}^N& \Bigg[\gamma_{i,k} \Bigg( \mean[p] \int_{x=t}^{\infty}\int_{y=x}^{\infty} e^{\rootp{k} (x-y)} d\pr(\workload[i] + \workload[i]' + \epts \leq y) dx \notag \\
                                  &- \mean[h] \int_{x=t}^{\infty}\int_{y=x}^{\infty} e^{\rootp{k} (x-y)} d\pr(\workload[i] + \workload[i]' + \ehts \leq y)dx \Bigg)\notag \\
                                  &-  \beta_{i,k} \Bigg( \mean[p] \int_{x=t}^{\infty}\int_{y=x}^{\infty} e^{\rootp{k} (x-y)} d\pr(\workload[i] + \epts \leq y)dx \notag \\
                                  &- \mean[h] \int_{x=t}^{\infty}\int_{y=x}^{\infty} e^{\rootp{k} (x-y)} d\pr(\workload[i] +\ehts \leq y)dx \Bigg) \notag \\
                                  &- \alpha_{i,k} \int_{x=t}^{\infty}\int_{y=x}^{\infty} e^{\rootp{k} (x-y)} d\pr(\workload[i]\leq y)dx \Bigg]\notag \\
   = \sum_{k=2}^N &\Bigg[ \gamma_{i,k} \Bigg( \mean[p] \int_{y=t}^{\infty} e^{-\rootp{k} y}d\pr(\workload[i] + \workload[i]' + \epts \leq y)\int_{x=t}^{y} e^{\rootp{k} x} dx \notag \\
                                  & - \mean[h] \int_{y=t}^{\infty} e^{-\rootp{k} y}d\pr(\workload[i] + \workload[i]' + \ehts \leq y) \int_{x=t}^{y} e^{\rootp{k} x} dx \Bigg) \notag \\
                                  &- \beta_{i,k} \Bigg( \mean[p] \int_{y=t}^{\infty}e^{-\rootp{k} y}d\pr(\workload[i] + \epts \leq y) \int_{x=t}^{y} e^{\rootp{k} x} dx \notag \\
                                  &- \mean[h] \int_{y=t}^{\infty}e^{-\rootp{k} y}d\pr(\workload[i] +\ehts \leq y) \int_{x=t}^{y} e^{\rootp{k} x} dx \Bigg)\notag \\
                                  &-  \alpha_{i,k} \int_{y=t}^{\infty} e^{-\rootp{k} y} d\pr(\workload[i]\leq y) \int_{x=t}^{y} e^{\rootp{k} x} dx \Bigg] \notag \\
   = \sum_{k=2}^N &\Bigg[ \frac{\gamma_{i,k}}{\rootp{k}} \Bigg( \mean[p] \int_{y=t}^{\infty} d\pr(\workload[i] + \workload[i]' + \epts \leq y) \notag \\
                                  &- \mean[h] \int_{y=t}^{\infty} d\pr(\workload[i] + \workload[i]' + \ehts \leq y) \Bigg) \notag \\
                                  &- \frac{\beta_{i,k}}{\rootp{k}} \Bigg( \mean[p] \int_{y=t}^{\infty}d\pr(\workload[i] + \epts \leq y) \notag \\
                                  &-\mean[h]  \int_{y=t}^{\infty}d\pr(\workload[i] +\ehts \leq y) \Bigg)
                                   -  \frac{\alpha_{i,k}}{\rootp{k}} \int_{y=t}^{\infty}  d\pr(\workload[i]\leq y) \notag \\
                                  &- \frac{\gamma_{i,k}}{\rootp{k}} \Bigg( \mean[p] \int_{y=t}^{\infty} e^{-\rootp{k} (y-t)} d\pr(\workload[i] + \workload[i]' + \epts \leq y) \notag \\
                                  &- \mean[h] \int_{y=t}^{\infty} e^{-\rootp{k} (y-t)} d\pr(\workload[i] + \workload[i]' + \ehts \leq y) \Bigg) \notag \\
                                  &+ \frac{\beta_{i,k}}{\rootp{k}} \Bigg( \mean[p] \int_{y=t}^{\infty} e^{-\rootp{k} (y-t)} d\pr(\workload[i] + \epts \leq y) \notag \\
                                  &- \mean[h]  \int_{y=t}^{\infty} e^{-\rootp{k} (y-t)} d\pr(\workload[i] +\ehts \leq y) \Bigg) \notag \\
                                  &+  \frac{\alpha_{i,k}}{\rootp{k}} \int_{y=t}^{\infty} e^{-\rootp{k} (y-t)} d\pr(\workload[i]\leq y) \Bigg] \notag \\
   =\sum_{k=2}^N &\frac1{\rootp{k}} \Bigg[ -\gamma_{i,k} \Bigg( \mean[p] \pr \big(t < \workload[i] + \workload[i]' + \epts < t + \erlang{\rootp{k}}\big) \notag \\
                                  &- \mean[h] \pr \big(t<\workload[i] + \workload[i]' + \ehts < t + \erlang{\rootp{k}} \big) \Bigg) \notag \\
                                  &+ \beta_{i,k} \Bigg( \mean[p] \pr\big(t<\workload[i] + \epts < t + \erlang{\rootp{k}} \big) \notag \\
                                  &-\mean[h] \pr\big(t<\workload[i] + \ehts < t + \erlang{\rootp{k}}\big) \Bigg) \notag \\
                                  &+ \alpha_{i,k} \pr\big(t<\workload[i] < t + \erlang{\rootp{k}} \big) \Bigg] \notag \\
   + \sum_{k=2}^N &\frac1{\rootp{k}} \Bigg[ \gamma_{i,k} \Bigg( \mean[p] \pr (\workload[i] + \workload[i]' + \epts >t ) \notag \\
                                  &- \mean[h] \pr (\workload[i] + \workload[i]' + \ehts >t ) \Bigg) \notag \\
                                  &- \beta_{i,k} \Bigg( \mean[p] \pr (\workload[i] + \epts >t ) - \mean[h] \pr (\workload[i] + \ehts >t ) \Bigg) \notag \\
                                  &- \alpha_{i,k} \pr (\workload[i] >t ) \Bigg]. \label{e3.middle terms}
\end{align}
\end{small}

By using now the property $\mathcal{L}^{-1}\{\frac{a^{n+1}}{(s+a)^{n+1}}\} = \frac1{n!}a^{n+1}t^n \allowbreak \times e^{-a t}$, $t\geq 0$, of the inverse Laplace transform, we see that the terms $\frac{(\rootnnum{i,j})^{r_{i,j}-l+1}}{(s+\rootnnum{i,j})^{r_{i,j}-l+1}}$ in Eq.~\eqref{e3.LST power series with corrected replace} correspond to the Laplace transform of an \erlang[r_{i,j}-l+1]{\rootnnum{i,j}} r.v. Combining all the above, the result in immediate, which completes the proof of the theorem.
\end{proof}

\begin{proof}[of Proposition~\ref{p3.extended continuity theorem}]
  In Proposition~\ref{p3.series expansion of the corrected replace transform vector}, we found that
  \begin{equation*}
    \ltwcmix{,i} = \ltwc{i} + \epsilon \fun{\lt{\theta_i}} + O(\epsilon^2),
  \end{equation*}
  where \fun{\lt{\theta_i}} is the Laplace-Stieltjes transform of the signed measure \fun[t]{\Theta_i} introduced in Theorem~\ref{t3.series expansion of mixture workload state i}. The above equation implies that
  \begin{equation}\label{e3.convergence rates for laplace transforms}
    \frac{\ltwcmix{,i} - \ltwc{i}}{\epsilon} = \fun{\lt{\theta_i}} + o(1).
  \end{equation}
  We set $n=\frac1{\epsilon}$ and we define the sequence of functions
  \begin{equation*}
    \fun{\lt{v}_n}: = \frac1{\epsilon}\big( \ltwcmix{,i} - \ltwc{i} \big),
  \end{equation*}
  where $\fun{\lt{v}_n}$ is the Laplace-Stieltjes transform of the measure $\fun[t]{V_n} = \frac1{\epsilon}\big( \pr(\workloadmix[,i] > t) - \pr(\workload[i] > t) \big)$. By using \eqref{e3.convergence rates for laplace transforms}, we obtain that $\fun{\lt{v}_n} \rightarrow \fun{\lt{\theta_i}}$, for all $s>0$ as $n \rightarrow \infty$ (or equivalently $\epsilon \rightarrow 0$). Thus, it follows from the Extended Continuity Theorem (see Theorem XIII.2 \cite{feller-IPTIA}) that $\frac{\pr(\workloadmix[,i] > t) - \pr(\workload[i] > t)}{\epsilon} \rightarrow \fun[t]{\Theta_i}$, which completes the proof.
\end{proof}

\end{document}